\documentclass{amsart}
\usepackage{graphicx,amsmath,amssymb,amsfonts,amscd,color}
\usepackage[all,cmtip,arrow,2cell]{xy}
\usepackage{mathrsfs}
\usepackage[colorlinks,
             linkcolor=black,
             citecolor=blue,
             bookmarksnumbered=true]{hyperref}
\usepackage{appendix}
\newtheorem{thm}{Theorem}[section]
\newtheorem{conj}[thm]{Conjecture}
\newtheorem{cor}[thm]{Corollary}
\newtheorem{lem}[thm]{Lemma}
\newtheorem{prop}[thm]{Proposition}
\newtheorem{propdef}[thm]{Proposition-Definition}

\newtheorem*{note}{Note}
\theoremstyle{definition}
\newtheorem{defn}[thm]{Definition}
\theoremstyle{definition}
\newtheorem{eg}[thm]{Example}
\theoremstyle{remark}
\newtheorem{rem}[thm]{Remark}
\theoremstyle{remark}
\newtheorem*{notation}{Notation}
\newtheorem*{acknowledgement}{Acknowledgement}
\numberwithin{equation}{section}

\newcommand{\set}[1]{\left\{#1\right\}}
\newcommand{\Ind}{\mathbf{Ind}}
\newcommand{\Sind}{\mathbf{Sind}}
\newcommand{\Fun}{\mathbf{Fun}}
\newcommand{\Lan}{\mathbf{Lan}}

\newcommand{\F}{\mathcal F}

\newcommand{\W}{\mathcal W}

\newcommand{\A}{\mathcal{A}}
\newcommand{\B}{\mathcal{B}}
\newcommand{\C}{\mathcal{C}}
\newcommand{\D}{\mathcal{D}}
\newcommand{\N}{\mathcal{N}}
\renewcommand{\P}{\mathcal{P}}
\newcommand{\Q}{\mathcal{Q}}
\newcommand{\R}{\mathcal{R}}
\newcommand{\J}{\mathcal{J}}
\renewcommand{\O}{\mathcal{O}}

\renewcommand{\AA}{\mathbb A}
\newcommand{\RR}{\mathbb R}
\newcommand{\QQ}{\mathbb Q}
\newcommand{\KK}{\mathbb K}
\newcommand{\ZZ}{\mathbb Z}
\newcommand{\NN}{\mathbb N}
\renewcommand{\SS}{\mathbb S}
\newcommand{\CC}{\mathbb C}
\newcommand{\AlgCat}{\mathbf{AlgCat}}
\newcommand{\algcat}{\mathscr{A}\!\mathit{lg}\mathscr{C}\!\mathit{at}}
\newcommand{\ATh}{\mathbf{ATh}}

\newcommand{\ath}{\mathscr{A\!T}\!\mathit{h}}

\newcommand{\LTh}{\mathbf{LTh}}
\newcommand{\Np}{\bE\Alg_{\mathbf{npd}}}

\newcommand{\To}{\longrightarrow}

\newcommand{\bB}{\mathbf{B}}
\newcommand{\bC}{\mathbf{C}}
\newcommand{\bD}{\mathbf{D}}
\newcommand{\bE}{\mathbf{E}}
\newcommand{\bK}{\mathbf{K}}
\newcommand{\bN}{\mathbf{N}}
\newcommand{\bQ}{\mathbf{Q}}
\newcommand{\bR}{\mathbf{R}}
\newcommand{\bSE}{\mathbf{SE}}
\newcommand{\bT}{\mathbf{T}}
\newcommand{\bS}{\mathbf{S}}
\newcommand{\bG}{\mathbf{G}}

\newcommand{\bW}{\mathbf{W}}
\newcommand{\bH}{\mathbf{H}}
\newcommand{\bF}{\mathbf{F}}
\newcommand{\bcinf}{\mathbf{\cinf}}
\newcommand{\bcomega}{\mathbf{\comega}}
\newcommand{\bcom}{\mathbf{Com}}
\newcommand{\bscom}{\mathbf{SCom}}

\newcommand{\even}{{\underline{0}}}
\newcommand{\odd}{{\underline{1}}}
\newcommand{\cinf}{{C^\infty}}
\newcommand{\comega}{{C^\omega}}

\newcommand{\Set}{\mathbf{Set}}

\newcommand{\End}{\mathbf{End}}

\newcommand{\Alg}{\mathbf{Alg}}
\newcommand{\comalg}{\mathbf{ComAlg}}

\newcommand{\comsalg}{\mathbf{SComAlg}}

\newcommand{\ev}{\mathrm{ev}}
\newcommand{\op}{\mathrm{op}}

\newcommand{\red}{\mathrm{red}}
\newcommand{\Hom}{\operatorname{Hom}}
\newcommand{\Ker}{\operatorname{Ker}}
\newcommand{\Spec}{\mathbf{Spec}}
\renewcommand{\Im}{\operatorname{Im}}
\newcommand{\del}{\partial}

\def\Rad{\mathfrak{Rad}}
\def\RAD{\mathfrak{R}}
\def\oinft{ \mbox{\,\put(-1.5,0){$\bigcirc$}\put(-0.2,1.05){\hbox{\tiny$\infty$}}}\mspace{16mu}}
\def\oi{\mspace{3mu}\resizebox{0.28cm}{!}{$\oinft$}\mspace{3mu}}
\def\oinfty{\raisebox{.24ex}{$\oi$}}
\def\oinftpy{\mspace{10mu}\mbox{\,\put(-4.73,0){$\odot$}\put(-3.33,0){\hbox{$\circ$}}}\mspace{10mu}}
\def\fweil{\mbox{formal Weil}}
\def\fweilk{\mbox{formal Weil $\KK$-algebra}}

\def\longlongrightarrow{-\!\!\!-\!\!\!-\!\!\!-\!\!\!-\!\!\!-\!\!\!\longrightarrow}

\newcommand{\Adj}[4]{\xymatrix@1{#2 \ar@<-0.5ex>[r]_-{#4} & #3 \ar@<-0.5ex>[l]_-{#1}}}
\newcommand{\Adjlong}[4]{\xymatrix@C=2cm{#2 \ar@<-0.65ex>[r]_-{#4} & #3 \ar@<-0.65ex>[l]_-{#1}}}
\begin{document}

\title{On Theories of Superalgebras of Differentiable Functions}%
\author{David Carchedi}
\address{D. Carchedi \hspace{16pt}\mbox{Max Planck Institute for Mathematics, Bonn, Germany. }
}
\author{Dmitry Roytenberg}
\address{D. Roytenberg\hspace{5pt}
Department of Mathematics,
Utrecht University,
The Netherlands.}







\date{\today}
\maketitle
\begin{abstract}
This is the first in a series of papers laying the foundations for a differential graded approach to derived differential geometry (and other geometries in characteristic zero). In this paper, we study theories of supercommutative algebras for which infinitely differentiable functions can be evaluated on elements. Such a theory is called a \emph{super Fermat theory.} Any category of superspaces and smooth functions has an associated such theory. This includes both real and complex supermanifolds, as well as algebraic superschemes. In particular, there is a super Fermat theory of \emph{$\bcinf$-superalgebras}. $\bcinf$-superalgebras are the appropriate notion of supercommutative algebras in the world of $\bcinf$-rings, the latter being of central importance both to synthetic differential geometry and to all existing models of derived smooth manifolds. A super Fermat theory is a natural generalization of the concept of a Fermat theory introduced by E. Dubuc and A. Kock. We show that any Fermat theory admits a canonical \emph{superization,} however not every super Fermat theory arises in this way. For a fixed super Fermat theory, we go on to study a special subcategory of algebras called near-point determined algebras, and derive many of their algebraic properties.
\end{abstract}



\tableofcontents
\section{Introduction.}
The purpose of this paper is to introduce \emph{super Fermat theories}. This theory will form the basis of our approach to differential graded models for derived manifolds. Super Fermat theories are theories of supercommutative algebras in which, in addition to evaluating polynomials on elements, one can evaluate infinitely differentiable functions. In particular, they provide a unifying framework to study the rings of functions of various flavors of smooth superspaces, e.g. regular functions on algebraic superschemes, smooth functions on smooth supermanifolds, and holomorphic functions on complex supermanifolds. The basic idea is to take seriously the notion that every type of geometry must have its own intrinsic version of commutative algebra associated to it (with the classical theory of commutative rings corresponding to algebraic geometry). Of central importance is the example of $\bcinf$-superalgebras which are the appropriate notion of supercommutative algebras in the world of $\bcinf$-rings, the commutative algebras associated to differential geometry.

A $\bcinf$-ring is a commutative $\RR$-algebra which, in addition to the binary operations of addition and multiplication, has an $n$-ary operation for each smooth function $$f:\RR^n \to \RR,$$ subject to natural compatibility. They were introduced by W. Lawvere in his Chicago lectures on \emph{categorical dynamics}, but first appeared in the literature in \cite{ReyWra} and \cite{dubuc1}. Their inception lies in the development of models for synthetic differential geometry \cite{cinfsch,smoothfun,loc1,arch,smzar,germint,msia,sdg}; however, recently they have played a pivotal role in developing models for \emph{derived} differential geometry \cite{spivak,joyce,borisovnoel}.

In \cite{1forms}, E. Dubuc and A. Kock introduce \emph{Fermat theories}, which provide a unifying framework for the algebraic study of polynomials using commutative rings, and the algebraic study of smooth functions using $\bcinf$-rings. Fermat theories are, in a precise way, theories of rings of infinitely differentiable functions. Recall that for a smooth function $f\left(x,z_1,\ldots,z_n\right)$ on $\RR^{n+1},$ there exists a unique smooth function $\frac{\Delta f}{\Delta x}\left(x,y,z_1,\ldots,z_n\right)$ on $\RR^{n+2}$ -- \textit{the difference quotient} -- such that for all $x$ and $y,$
\begin{equation}\label{eq:Fermtprop}
f\left(x,\mathbf{z}\right)-f\left(y,\mathbf{z}\right)=\left(x-y\right)\cdot \frac{\Delta f}{\Delta x}\left(x,y,\mathbf{z}\right).
\end{equation}
Indeed, for $x' \ne y',$ $$\frac{\Delta f}{\Delta x}\left(x',y',\mathbf{z}\right)=\frac{f\left(x',\mathbf{z}\right)-f\left(y',\mathbf{z}\right)}{x'-y'}$$ but for $x'=y',$
\begin{equation}\label{eq:partial}
\frac{\Delta f}{\Delta x}\left(x',x',\mathbf{z}\right)= \frac{\partial f}{\partial x} \left(x'\right),
\end{equation}
a result known as \emph{Hadamard's Lemma}. In fact, if one took (\ref{eq:partial}) as a \emph{definition} of the partial derivative, all of the classical rules for differentiation could be derived from (\ref{eq:Fermtprop}) using only algebra. The key insight of Dubuc and Kock in \cite{1forms} is that equation (\ref{eq:Fermtprop}) makes sense in a more general setting, and one can consider algebraic theories extending the theory of commutative rings whose operations are labeled by functions satisfying a generalization of (\ref{eq:Fermtprop}) called the \emph{Fermat property} (since Fermat was the first to observe that it holds for polynomials). For examples and non-examples of Fermat theories, see Section \ref{sec:fermex}. Many important properties of $\bcinf$-rings hold for any Fermat theory. For example, if $\bE$ is a Fermat theory, $\A$ an $\bE$-algebra, and $I \subset \A$ an ideal of the underlying ring, then $\A/I$ has the canonical structure of an $\bE$-algebra. Moreover, the fact that the theory of $\bcinf$-rings satisfies the Fermat property is a key ingredient in many well-adapted models of synthetic differential geometry.

In this paper, we show that the Fermat property is ideally suited to study supergeometry, as any Fermat theory admits a canonical \emph{superization}. The superization of a Fermat theory is a $2$-sorted algebraic theory extending the theory of supercommutative algebras, and satisfies a modified version of the Fermat property which we call the \emph{super Fermat property.} An important example of a super Fermat theory is the theory of super $\cinf$-rings, which is the  superization of the theory of $\cinf$-rings. However, super Fermat theories are more general than Fermat theories, as not every super Fermat theory arises as a superization.

Super Fermat theories are an essential ingredient to our development of a differential graded approach to derived differential geometry. Of particular importance is the theory of $\bcinf$-superalgebras. Our approach is based upon exploiting the connection between supercommutativity and differential graded algebras. In \cite{dg2}, we define the concept of a differential graded $\bE$-algebra for a super Fermat theory $\bE,$ and develop homological algebra in this setting.

In light of the history of $\bcinf$-rings and their role in synthetic differential geometry, it is natural to believe that super Fermat theories should play a pivotal role in synthetic supergeometry, but we do not pursue this in this paper. It is worth mentioning however, that our notion of superization is different from that of Yetter's \cite{yetter}, as his approach results in a uni-sorted Lawvere theory, and also applies in a more restrictive context; the superization of the theory of $\bcinf$-rings in Yetter's sense embeds diagonally into our $2$-sorted superization. Our theory is also quite different from that of \cite{nish} as his theory concerns itself with  $G^\infty$-supermanifolds, whereas our approach is more in tune with supermanifolds in the sense of \cite{manin}.


\subsection{Organization and main results}
In Section \ref{sec:ffermat}, we begin by reviewing the concept of a Fermat theory introduced in \cite{1forms}. We then introduce the concept of a \emph{reduced} Fermat theory, which is a Fermat theory that is, in a precise sense, a ``theory  of functions.'' We go on to show that we can associate to any Fermat theory a reduced Fermat theory in a functorial way; moreover, for every commutative ring $\KK$ there is a \emph{maximal} reduced Fermat theory with $\KK$ as the ground ring (for $\KK=\RR$, we recover the theory $\bcinf$ of smooth functions).

Section \ref{sec:superfermat} introduces the main subject of this paper, the concept of a super Fermat theory. We show that any Fermat theory has associated to it a canonical super Fermat theory called its \emph{superization}, and conversely, any super Fermat theory has an underlying Fermat theory. Moreover, we prove the following:

\begin{thm}(Corollary \ref{cor:superisadj})
The superization functor $$\bS:\mathbf{FTh} \to \mathbf{SFTh}$$ from Fermat theories to super Fermat theories is left adjoint to the underlying functor $$\left(\quad\right)_\even:\mathbf{SFTh} \to \mathbf{FTh}.$$
\end{thm}
We also develop some aspects of supercommutative algebra for algebras over a super Fermat theory and prove a useful property of morphisms of super Fermat theories:

\begin{thm} \ref{thm:morphismprodspres}
Let $F:\bS\to\bS'$ be a morphism of super Fermat theories. Then the induced functor $$F_!:\bS\Alg\to\bS'\Alg$$ preserves finite products.
\end{thm}

In Section \ref{sec:npd}, we begin the study of near-point determined algebras for a super Fermat theory. This is a generalization of the notion near-point determined introduced in \cite{msia} for finitely generated $\cinf$-rings to the setting of not necessarily finitely generated algebras over any super Fermat theory. We then prove that near-point determined algebras are completely determined by their underlying $\KK$-algebra, where $\KK$ is the ground ring of the theory:

\begin{thm} (Corollary \ref{cor:npdralgs})
If $\A$ and $\B$ are $\bE$-algebras and $\B$ is near-point determined, then any $\KK$-algebra morphism $\varphi:\A \to \B$ is a map of $\bE$-algebras.
\end{thm}

This is a broad generalization of the result proven by Borisov in \cite{Borisov} in the case of $\bcinf$-rings. Borisov uses topological methods in his proof, tailored specifically to the case of $\bcinf$-rings. We show that this result holds in a much more general context, and follows by completely elementary algebraic methods.

We go on to define what it means for a super Fermat theory to be \emph{super reduced}, a subtle generalization of the notion of a reduced Fermat theory suitable in the supergeometric context, and show that any free algebra for a super reduced Fermat theory is near-point determined. We end this section by investigating to what extent near-point determined algebras over a super Fermat theory $\bE$ whose ground ring $\KK$ is a field are flat. Near-point determined $\bE$-algebras are a reflective subcategory of $\bE$-algebras, and hence have their own tensor product $\oinftpy$ (coproduct). In particular, we prove the following:

\begin{thm} (Lemma \ref{lem:npfinp} and Lemma \ref{lem:monoprev})
For any near-point determined $\bE$-algebra $A,$ the endofunctor
$$\A \oinftpy \left(\quad\right): \Np \to \Np$$ preserves finite products and monomorphisms, where $\Np$ is the category of near-point determined $\bE$-algebras.
\end{thm}

Finally, in the appendices, we give a detailed introduction to algebraic theories and multi-sorted Lawvere theories, mostly following \cite{algthy,borceux2}, and introduce many of the conventions and notations concerning their use in this paper.

\begin{acknowledgement}
We would like to thank Mathieu Anel, Christian Blohmann, Dennis Borisov, Eduardo Dubuc, Wilberd van der Kallen, Anders Kock, Ieke Moerdijk, Justin Noel, Jan Stienstra, and Peter Teichner for useful conversations. The first author would like to additionally thank the many participants in the ``Higher Differential Geometry'' seminar (formerly known as the ``Derived Differential Geometry'' seminar) at the Max Planck Institute for Mathematics. The second author was supported by the Dutch Science Foundation ``Free Competition'' grant. He would also like to thank the Radboud University of Nijmegen, where part of this work was carried out, for hospitality.
\end{acknowledgement}

\section{Fermat Theories.}\label{sec:ffermat}
\subsection{Examples of Lawvere Theories}
A review of the basics of algebraic theories and multi-sorted Lawvere theories, as well as many notational conventions concerning them, is given in Appendices \ref{sec:theories} and \ref{sec:concrete}.

Before presenting Fermat theories in general, we begin by introducing some instructive motivating examples:
\begin{eg}
Let $\bcom$ be the opposite category of finitely generated free commutative (unital, associative) rings. Up to isomorphism, its objects are of the form $$\ZZ[x_1,\cdots,x_n] \cong \ZZ[x]^{\otimes n}.$$ Since we are in the opposite category, and the tensor product of commutative rings is the coproduct, every object of $\bcom$ is a finite product of the object $\ZZ[x].$ With this chosen generator, $\bcom$ is a Lawvere theory.

It is sometimes useful to take the dual geometric viewpoint. We can consider the category whose objects are finite dimensional affine spaces $\AA_{\ZZ}^n$ over $\ZZ,$ so their morphisms are polynomial functions with integer coefficients. This category is canonically equivalent to $\bcom.$ Indeed, each affine space $\AA_{\ZZ}^n$ corresponds to the ring $\ZZ[x_1,\cdots,x_n],$ and since $\ZZ[x]$ is the free commutative ring on one generator, we have the following string of natural isomorphisms:
\begin{eqnarray*}
\Hom\left(\AA_{\ZZ}^n,\AA_{\ZZ}^m\right) &\cong& \underline{\ZZ[x_1,\cdots,x_n]^m}\\
&\cong& \Hom\left(\ZZ[x],\ZZ[x_1,\cdots,x_n]\right)^m\\
&\cong& \Hom\left(\ZZ[x]^{\otimes m},\ZZ[x_1,\cdots,x_n]\right)\\
&\cong& \Hom\left(\ZZ[x_1,\cdots,x_m],\ZZ[x_1,\cdots,x_n]\right),\\
\end{eqnarray*}
where $\underline{\ZZ[x_1,\cdots,x_n]^m}$ denotes the underlying set of the ring.

Notice that the affine line $\AA_{\ZZ}$ is a commutative ring object in $\bcom$. Indeed, the polynomial $$m\left(x,y\right)=x\cdot y \in \ZZ[x,y]$$ is classified by a morphism $$\ZZ[x] \to \ZZ[x,y]$$ which corresponds to a morphism $$m:\AA_{\ZZ}^2 \to \AA_{\ZZ},$$ which is multiplication. Similarly, the polynomial $$a\left(x,y\right)=x+y$$ induces a map $$a:\AA_{\ZZ}^2 \to \AA_{\ZZ},$$ which is addition. Finally, the ring unit $$1 \in \ZZ[x]$$ induces a map $$u:\AA^0_{\ZZ} \to \AA_{\ZZ}$$ which is the unit map of this ring object.

Let $A$ be a commutative ring. Then it induces a functor
\begin{eqnarray*}
\tilde A:\bcom &\to& \Set\\
\ZZ[x_1,\cdots,x_n] &\mapsto& \Hom\left(\ZZ[x_1,\cdots,x_n],A\right),
\end{eqnarray*}
which is product preserving, hence a $\bcom$-algebra. Moreover, since $\ZZ[x]$ is the free commutative ring on one generator, the underlying set of $\tilde A$ is $\tilde A\left(\ZZ[x]\right) \cong \underline{A},$ the underlying set of $A.$

Conversely, suppose that $\B$ is a $\bcom$-algebra. Then, as it is a finite product preserving functor, and the diagram expressing that an object in a category is a ring object only uses finite products, it follows that the data $$\left(\underline{\B}:=\B\left(\AA_{\ZZ}\right),\B\left(m\right),\B\left(a\right),\B\left(u\right)\right)$$ encodes a commutative ring (in $\Set.$) Moreover, it can be checked that if $$\mu:\B \Rightarrow \B'$$ is a natural transformation between product preserving functors from $\bcom$ to $\Set,$ that $$\mu\left(\AA_{\ZZ}\right):\underline{\B} \to \underline{\B'}$$ is a ring homomorphism, and conversely, if $$\varphi:A \to A'$$ is a ring homomorphism, $$\mu\left(\AA^n_{\ZZ}\right)=\varphi^n:\underline{A}^n \to \underline{A'}^n$$ defines a natural transformation $$\tilde A \Rightarrow \tilde A'.$$ It follows that the category $\bcom\Alg$ is equivalent to the category of commutative rings.

Using the notation (\ref{eq:notta}), one has that $$\bcom\left(n\right)=\ZZ[x_1,\cdots,x_n]$$ and its underlying set is given by $$\underline{\ZZ[x_1,\cdots,x_n]} = \bcom\left(n,1\right) \cong \Hom\left(\AA_\ZZ^n,\AA_\ZZ\right).$$ One may succinctly say that $\bcom$ is the Lawvere theory whose $n$-ary operations are labeled by the elements of $\ZZ[x_1,\ldots,x_n],$ and its algebras are commutative rings.

Notice that for a given commutative ring $A,$ congruences of $A$ are in bijection with ideals. Indeed, given an ideal $I,$ it defines a subring $R\left(I\right)$ of $A \times A$ whose elements are pairs $\left(a,a'\right)$ such that $$a-a' \in I.$$ Conversely, given a congruence $R \rightarrowtail A \times A,$ the subset $$I:=\left\{a \in A\mbox{  }|  \left(a,0\right) \in R.\right\},$$ is an ideal of $A.$
\end{eg}

\begin{eg}
Let $\KK$ be a commutative ring. Then one may consider $\bcom_{\KK}$ to be the opposite category of finitely generated free $\KK$-algebras. Up to isomorphism, its objects are of the form $$\KK[x_1,\cdots,x_n] \cong \KK[x]^{\otimes^n_{\KK}}.$$ Since we are in the opposite category, tensoring over $\KK$ corresponds to taking the product, and this category is a Lawvere theory with generator $$\KK[x]\cong \KK \otimes_{\ZZ} \ZZ[x].$$ We see that this Lawvere theory is a particular instance of Remark \ref{rem:relative}. Algebras for this Lawvere theory are precisely $\KK$-algebras, and congruences are again ideals. One may also view the category $\bcom_\KK$ as the category of finite dimensional affine planes $\AA^n_\KK$ over $\KK.$
\end{eg}

\begin{eg}\label{ex:cinf}
When $\KK=\RR,$ one may view the category $\bcom_\RR$ as the category whose objects are manifolds of the form $\RR^n$ and whose morphisms are polynomial functions. From the geometric view point, it is natural to ask what happens if one allows arbitrary smooth functions instead. The resulting category, which is a full subcategory of the category of smooth manifolds $\textbf{Mfd},$ is a Lawvere theory with generator $\RR$. We denote this Lawvere theory by $\bcinf.$ It is the motivating example for this paper. It may be described succinctly by saying its $n$-ary operations are labeled by elements of $$\bcinf(n,m)=\cinf(\RR^n,\RR^m).$$ Notice that, there is a canonically induced functor $$\bcom_\RR \to \bcinf$$ sending each manifold to itself, and each polynomial to itself viewed as a smooth function. This is a map of Lawvere theories, so there is an induced adjunction
$$\Adjlong{\widehat{(\quad)}}{\bcom_\RR\Alg}{\bcinf\Alg,}{(\quad)_\sharp}$$ where, for a $\bcinf$-algebra $A$, $A_\sharp$ is its underlying $\RR$-algebra, and if $R$ is an $\RR$-algebra, $\widehat{R}$ is its $\cinf$-completion. In particular, $$\widehat{\RR[x_1,\cdots,x_n]} \cong \bcinf\left(\RR^n\right).$$ The functor $(\quad)_\sharp$ is faithful and conservative, therefore one may regard a $\bcinf$-algebra as an $\RR$-algebra with extra structure. This extra structure is encoded by a whole slew of $n$-ary operations, one for each smooth function $$f:\RR^n \to \RR,$$ subject to natural compatibility. For example, if $M$ is a smooth manifold, then it induces a product preserving functor
\begin{eqnarray*}
\bcinf\left(M\right):\bcinf &\to& \Set\\
\RR^n &\mapsto& \Hom\left(M,\RR^n\right).
\end{eqnarray*}
$\bcinf\left(M\right)$ is a $\bcinf$-algebra whose underlying $\RR$-algebra is the ordinary ring of smooth functions $\cinf\left(M\right).$ Given a smooth function $$f:\RR^n \to \RR,$$ it induces an $n$-ary operation $$\bcinf\left(M\right)\left(f\right):\bcinf\left(M\right)^n \to \bcinf\left(M\right),$$ defined by $$\bcinf\left(M\right)\left(f\right)\left(\varphi_1,\cdots,\varphi_n\right)\left(x\right)=f\left(\varphi_1\left(x\right),\cdots,\varphi_n\left(x\right)\right).$$
$\bcinf$-algebras come with their own notion of tensor product (coproduct), and we denote the $\bcinf$-tensor product of $A$ and $B$ by $A \oinfty B.$ Unlike for the ordinary tensor product of $\RR$-algebras, one has for (Hausdorff, second countable) smooth manifolds $M$ and $N,$ the equality \cite{msia}: $$\bcinf\left(M\right) \oinfty \mspace{4mu} \bcinf\left(N\right) \cong \bcinf\left(M \times N\right).$$ Hence, they are ideally suited for the theory of manifolds. At the same time, the theory $\bcinf$-algebras closely resembles the theory of commutative rings, as it enjoys a very nice property, namely the \emph{Fermat property}, which is the subject of the next subsection.
\end{eg}

\subsection{Fermat theories}\label{sec:fermat}

A \emph{Fermat theory} is an extension $\bE$ of $\bcom$ that has an intrinsic notion of derivative obeying the expected rules (the chain rule, the Taylor formula, etc.). Standard notions of differential calculus, such as derivations and differentials, can be defined for $\bE$-algebras. The notion of Fermat theory was introduced and studied by Dubuc and Kock in \cite{1forms}. In what follows, we recall some key definitions and results from that paper.

\subsubsection{The Fermat Property.}

Let $\bE$ be extension of $\bcom,$ that is a Lawvere theory $\bE$ with a map of Lawvere theories $\tau_\bE:\bcom\to\bE$, i.e. an object of the undercategory $\bcom/\LTh.$ The structure map $\tau_\bE$ induces an adjunction
$$\Adj{\tau_\bE^*}{\bcom\Alg}{\bE\Alg}{\tau^\bE_!}.$$
Observe first that $\KK=\bE(0)$, being the free $\bE$-algebra on the empty set, has an underlying ring structure. Categorically, $\bE(0)$ is a finite product preserving $\Set$-valued functor, and composition with $\tau_\bE$ induces a $\bcom$-algebra
$$\bcom \stackrel{\tau_\bE}{\longlongrightarrow} \bE \stackrel{\bE(0)}{\longlongrightarrow} \Set.$$ This $\bcom$-algebra is just $\tau^*_\bE\left(\bE(0)\right),$ and since the underlying set of an algebra for a Lawvere theory is determined by its value as a functor on the generator, and $\tau_E$ preserves generators, $\tau^*_\bE\left(\bE(0)\right)$ has the same underlying set as $\bE(0).$ Now that this is clear, we will abuse notation and denote
$\tau^*_\bE\left(\A\right),$ for an $\bE$-algebra $\A,$ simply by $\A.$
On one hand, since $\bE(0)$ is the initial $\bE$-algebra, there is a unique $\bE$-algebra map from $\KK$ to $\bE(n)$ for each $n$. In particular, it is a map of rings. On the other hand, the unit of the adjunction $\tau^\bE_! \dashv \tau^*_\bE$ is map of rings $$\bcom(n)\to\bE(n).$$ Hence, we have a map of rings
\[
\KK[x_1,\ldots,x_n]=\bcom_\KK(n)=\KK\otimes_\ZZ\bcom(n)\rightarrow\bE(n).
\]
Since this is obviously compatible with compositions, we deduce that the structure map $\tau_\bE:\bcom\rightarrow\bE$ factors through $\bcom_\KK$. So every $\bE$-algebra has an underlying commutative $\KK$-algebra structure. Let us denote the corresponding forgetful functor by
\[
(\quad)_\sharp:\bE\Alg\To\comalg_\KK
\]
and its left adjoint -- the $\bE$-algebra completion -- by
\[
\widehat{(\quad)}:\comalg_\KK\To\bE\Alg.
\]
We shall refer to $\KK$ as the \emph{ground ring} of the theory $\bE$.
\begin{rem}
The name \emph{completion} should be taken with a grain of salt. For example, the theory of $\bcinf$-algebras is an extension of $\bcom$ with ground ring $\RR$. Consider $\CC$ as an $\RR$-algebra. We can present it as $$\CC=\RR[x]/\left(x^2+1\right).$$ It follows that the $\bcinf$-completion of $\CC$ is $$\widehat{\CC}=\bcinf\left(\RR\right)/\left(x^2+1\right),$$ however the function $x^2+1$ is a unit in $\bcinf\left(\RR\right),$ so we have that $\widehat{\CC}=\{0\}$, the terminal algebra.
\end{rem}

\begin{notation}
Denote by $\odot$ the binary coproduct in $\bE\Alg$. Denote the free $\bE$-algebra on generators $x_1,\ldots,x_n$ by
$\KK\{x_1,\ldots,x_n\}$ (or $\KK_\bE\{x_1,\ldots,x_n\}$ when there are several theories around and we need to be clear which one we mean). It is synonymous with $\bE(n)$, but with the generators named explicitly. If $\A\in\bE\Alg$, let
\[
\A\{x_1,\ldots,x_n\}=\A\odot\KK\{x_1,\ldots,x_n\}.
\]
It solves the problem of universally adjoining variables to an $\bE$-algebra.
\end{notation}

\begin{rem}
The $\bE$-algebra completion of $\KK[x_1,\ldots,x_n]$ is $\KK\{x_1,\ldots,x_n\}$.
\end{rem}

\begin{defn}\cite{1forms}
An extension $\bE$ of $\bcom$ is called a \emph{Fermat theory} if  for every $f\in\KK\{x,z_1,\ldots,z_n\}$ there exists a unique $\frac{\Delta f}{\Delta x}\in\KK\{x,y,z_1,\ldots,z_n\}$, called the \textit{difference quotient}, such that
\begin{equation}\label{eqn:fermatproperty}
f(x,\mathbf{z})-f(y,\mathbf{z})=(x-y)\cdot\frac{\Delta f}{\Delta x}(x,y,\mathbf{z})
\end{equation}
where $\mathbf{z}=(z_1,\ldots,z_n)$. A \emph{Fermat theory over $\QQ$} is a Fermat theory whose structure map factors through $\bcom_\QQ$.

Let $\mathbf{FTh}$ (resp. $\mathbf{FTh}_{/\QQ}$, $\mathbf{FTh}_\KK$) denote the full subcategory of $\bcom/\mathbf{LTh}$ (resp. $\bcom_\QQ/\mathbf{LTh}$, $\bcom_\KK/\mathbf{LTh}$) consisting of the Fermat theories (resp. Fermat theories over $\QQ$, Fermat theories with ground ring $\KK$).
\end{defn}

For the rest of this subsection, let $\bE$ denote a Fermat theory with ground ring $\KK$. 

\begin{note}
 The ground ring of a Fermat theory over $\QQ$ always contains $\QQ$ but is generally different from it, so the categories $\mathbf{FTh}_\QQ$ and $\mathbf{FTh}_{/\QQ}$ are different.
 \end{note}

An immediate consequence of the Fermat property is the following:

\begin{prop}\label{prop:hadamard}\cite{1forms}
For any $f\in\KK\{x_1,\ldots,x_n\}$, there exist $$g_i\in\KK\{x_1,\ldots,x_n,y_1,\ldots,y_n\},$$ $i=1,\ldots,n$, such that
\begin{equation}\label{eq:hada}
f\left(x_1,\ldots,x_n\right)-f\left(y_1,\ldots,y_n\right)=\sum_{i=1}^n\left(x_i-y_i\right) \cdot g_i\left(x_1,\ldots,x_n,y_1,\ldots,y_n\right).
\end{equation}
\end{prop}

The following corollary is the cornerstone of the theory of Fermat theories:

\begin{cor}\label{cor:idealsarecongs}\cite{1forms}
Let $\A$ be an $\bE$-algebra, $I\subset\A$ an ideal in the underlying ring. Then the induced equivalence relation on $\A$ ($a\sim b$ modulo $I$ iff $a-b\in I$) is an $\bE$-congruence. Consequently, there is a unique $\bE$-algebra structure on $\A/I$ making the projection $\A\to\A/I$ a map of $\bE$-algebras.
\end{cor}

\begin{proof}
It suffices to show that if $I$ is an ideal of $\A,$ and $$a_1,\ldots,a_n$$ and $$b_1,\ldots,b_n$$ are in $\A$ such that for each $i,$ $$a_i - b_i \in I,$$ then for each $f \in \bE\left(n,1\right),$
$$A\left(f\right)\left(a_1,\ldots,a_n\right) - A\left(f\right)\left(b_1,\ldots,b_n\right) \in I.$$ There exists a unique morphism $$\varphi:\KK\{x_1,\ldots,x_n,y_1,\ldots,y_n\} \to \A$$ sending each $x_i$ to $a_i$ and each $y_i$ to $b_i.$ Note that by \ref{prop:hadamard} there exists $$g_1,\ldots,g_n \in \KK\{x_1,\ldots,x_n,y_1,\ldots,y_n\}$$ such that (\ref{eq:hada}) holds. Notice for any $g \in \KK\{x_1,\ldots,x_n,y_1,\ldots,y_n\},$ $$\varphi\left(g\right)=A\left(g\right)\left(\varphi\left(x_1\right),\ldots,\varphi\left(x_n\right),\varphi\left(y_1\right),\ldots,\varphi\left(y_n\right)\right).$$
It follows that
$$A\left(f\right)\left(a_1,\ldots,a_n\right) - A\left(f\right)\left(b_1,\ldots,b_n\right) = \sum\limits_{i=1}^n \left(a_i-b_i\right) \cdot A\left(g_i\right)\left(a_1,\ldots,a_n,b_1,\ldots,b_n\right) \in I.$$
\end{proof}

\subsubsection{Derivatives.}
Suppose we are given an $f\in\KK\{x^1,\ldots,x^n\}$. Fix an $i\in\{1,\ldots,n\}$, let $x=x^i$ and consider $f$ as an element of $\R\{x\}$ with $\R=\KK\{x^1,\ldots,\hat{x}^i,\ldots,x^n\}$ (the hat indicates omission). By the Fermat property \eqref{eqn:fermatproperty}, there is a unique $\frac{\Delta f}{\Delta x}\in\R\{x,y\}$ such that
\[
f(x)-f(y)=(x-y)\cdot\frac{\Delta f}{\Delta x}(x,y).
\]
Define the \emph{partial derivative of $f$ with respect to $x^i$} to be
\[
\del_if=\frac{\del f}{\del x^i}=\frac{\Delta f}{\Delta x}(x,x)\in\R\{x\}=\KK\{x^1,\ldots,x^n\}.
\]
When $n=1$, we shall also write $f'(x)$ for $\del f/\del x^1$.

\medskip

As expected, the partial derivatives satisfy the chain rule:
\begin{prop}\label{prop:chain}\cite{1forms}
Let $\varphi\in\bE(k,1)$, $f=(f^1,\ldots,f^k)\in\bE(n,k)$. Then for all $i=1,\ldots,n$ we have
\[
\del_i(\varphi\circ f)=\sum_{j=1}^k(\del_if^j)(\del_j\varphi\circ f).
\]
\end{prop}
Here the partial derivatives can be interpreted as operators
\[
\del_i:\bE(n)\To\bE(n)
\]
on the $\bE$-algebra $\bE(n)$, satisfying a ``derivation rule'' for every $k$-ary $\bE$-operation $\varphi$ on $\bE(n)$ \cite{1forms}. In particular, letting $\varphi$ be addition (resp. the multiplication) we get the familiar $\KK$-linearity (resp. Leibniz rule).
\begin{rem}
We have slightly abused notation since the partial derivative operators $\del_i$ are \emph{not} morphisms of $\bE$-algebras.
\end{rem}

\begin{prop}\label{prop:clairaut} \emph{(Clairaut's theorem).}
The partial derivatives commute:
\[
\del_i\del_j=\del_j\del_i\quad\forall i,j.
\]
\end{prop}
\begin{proof}
We shall give the proof in the two-variable case only; the general case is proven in exactly the same way. Let $f=f(x,y)\in\bE(2,1)$. We obviously have
\[
(f(x,y)-f(z,y))-(f(x,w)-f(z,w))=(f(x,y)-f(x,w))-(f(z,y)-f(z,w)).
\]
Applying the Fermat property on both sides we get
\[
(x-z)(g(x,z,y)-g(x,z,w))=(y-w)(h(x,y,w)-h(z,y,w))
\]
for unique difference quotients $g$ and $h$. Applying the Fermat property again, we get
\[
(x-z)(y-w)\phi(x,z,y,w)=(y-w)(x-z)\psi(x,z,y,w)
\]
for unique difference quotients $\phi$ and $\psi$. By uniqueness, $x-z$ and $y-w$ are not zero-divisors, hence
\[
\phi(x,z,y,w)=\psi(x,z,y,w).
\]
Now, setting $x=z$ and $y=w$, we obtain the sought after
\[
\del_2\del_1f(x,y)=\del_1\del_2f(x,y).
\]
\end{proof}

\begin{cor}\label{cor:taylor}\cite{1forms} \emph{(The Taylor formula).}
For any $f\in\KK\{x_1,\ldots,x_p,z_1,\ldots,z_q\}$, $n\geq0$ and multi-indices $\alpha$ and $\beta$, there exist unique $h_\alpha\in\KK\{x_1,\ldots,x_p,z_1,\ldots,z_q\}$ and (not necessarily unique) $g_\beta\in\KK\{x_1,\ldots,x_p,y_1,\ldots,y_p,z_1,\ldots,z_q\}$ such that
\begin{equation}\label{eqn:taylor}
f(\mathbf{x}+\mathbf{y},\mathbf{z})=\sum_{|\alpha|=0}^nh_\alpha(\mathbf{x},\mathbf{z})\mathbf{y}^\alpha+\sum_{|\beta|=n+1}\mathbf{y}^\beta g_\beta(\mathbf{x},\mathbf{y},\mathbf{z}).
\end{equation}
Furthermore, if $\KK\supset\QQ$, we have
\[
h_\alpha(\mathbf{x},\mathbf{z})=\frac{\del^\alpha_xf(\mathbf{x},\mathbf{z})}{\alpha!},
\]
the usual Taylor coefficients.
\end{cor}

\subsubsection{Examples and non-examples.}\label{sec:fermex}
\begin{eg}
The theory $\bcom$ of commutative algebras is itself a Fermat theory, with ground ring $\ZZ$. It is the \emph{initial} Fermat theory. Similarly, $\bcom_\KK$ is the initial Fermat theory over $\KK.$
\end{eg}

\begin{eg}\label{eg:cinf}
 The theory $\bcinf$ of $\cinf$-algebras, with $\bcinf(n,m)=\cinf(\RR^n,\RR^m)$, the set of real smooth functions is a Fermat theory. This is Example \ref{ex:cinf}. The category $\cinf$ is the full subcategory of smooth manifolds spanned by those of the form $\RR^n,$ and the ground ring of this theory is $\RR$.
\end{eg}

\begin{eg} The theory $\bcomega$ of real analytic algebras, with\\ $\bcomega(n,m)=\comega(\RR^n,\RR^m)$, the set of real analytic functions is a Fermat theory. It is equivalent to the full subcategory of real analytic manifolds spanned by those of the form $\RR^n.$ The ground ring is again $\RR$.
\end{eg}

\begin{eg} The theory $\bH$ of complex holomorphic algebras, with $\bH(n,m)=H(\CC^n,\CC^m)$, the set of complex holomorphic (entire) functions, is a Fermat theory. The category is equivalent to the full subcategory of complex manifolds spanned by those of the form $\CC^n.$ The ground ring of this theory is $\CC$.
\end{eg}

\begin{eg} The theory $\bH^\RR$ of real holomorphic algebras, with $\bH^\RR(n,m)=H(\RR^n,\RR^m)$, the set of those entire functions which are invariant under complex conjugation, is a Fermat theory. The ground ring is $\RR$.
\end{eg}

\begin{eg}
Let $\KK$ be an integral domain, $\Bbbk$ its field of fractions. Let the theory $\bR^\KK$ consist of global rational functions, i.e. rational functions with coefficients in $\KK$ having no poles in $\Bbbk$. It is a Fermat theory with ground ring $\Bbbk$.
\end{eg}

\begin{eg} The theory $\bcinf^\CC$ with $\bcinf^\CC(n,m)=\bcinf(\CC^n,\CC^m)$, the set of functions which are smooth when viewed as functions from $\RR^{2n}$ to $\RR^{2m}$, is \emph{not} a Fermat theory, as the Fermat property for complex-valued functions implies the Cauchy-Riemann equations. Likewise, the theory $\bcomega^\CC$, defined analogously, is not a Fermat theory.
\end{eg}

\begin{eg}\label{eg:Ck} The theory $\mathbf{C^k}$ of $k$ times continuously differentiable real functions is \emph{not} a Fermat theory for any $0\leq k<\infty$: given an $f$ of class $C^k$, the difference quotient appearing in \eqref{eqn:fermatproperty} is only of class $C^{k-1}$.
\end{eg}

\begin{eg} As shown in \cite{1forms}, if $\bE$ is a Fermat theory and $\A$ is any $\bE$-algebra, the theory $\bE_\A$ of $\bE$-algebras over $\A$ is also Fermat, with $\A$ as the ground ring. This gives many examples of Fermat theories.
\end{eg}

We have proper inclusions of Fermat theories
\[
\bcom_\RR\subsetneq\bH^\RR\subsetneq\bcomega\subsetneq\bcinf,\quad\bcom_\RR\subsetneq\bcom_\CC\subsetneq\bH\quad
\mathrm{and}\quad\bH^\RR\subsetneq\bH,
\]
making various diagrams commute.

\begin{rem}
As $\CC$ is neither a $\bcinf$- nor $\bcomega$-algebra, and nor is $\RR$ an $\bH$-algebra, putting superscripts instead of subscripts in our notation for the theories $\bcinf^\CC$, $\bcomega^\CC$ and $\bH^\RR$ avoids possible confusion. However, notice that $\CC$ \emph{is} an $\bH_\RR$-algebra, and $\bH^\RR_\CC=\bH$.
\end{rem}

\begin{eg} Fermat theories have associated geometries. For instance, if $X$ is a smooth (resp. real analytic, complex) manifold, its structure sheaf $\O_X$ is actually a sheaf of $\bcinf$- (resp. $\bcomega$-, $\bH$-) algebras. If $X$ is real analytic, $X_\CC$ its complexification, the $\bcomega$-algebra structure on $\O_X$ does not extend to $\O_{X_\CC}$, but the underlying $\bH^\RR$-structure does extend to an $\bH^\RR_\CC=\bH$-algebra structure on $\O_{X_\CC}$.
\end{eg}

\subsubsection{Evaluations at $\KK$-points.}

Let $\bE$ be an extension of $\bcom.$ For the initial $\bE$-algebra $\KK$, we can think of its $\bE$-algebra structure as a collection of evaluation maps
\begin{equation}\label{eqn:ev}
\ev^n:\KK\{x_1,\ldots,x_n\}\To\Set(\KK^n,\KK),\quad n\geq0,
\end{equation}
where, for $f=f(x_1,\ldots,x_n)\in\KK\{x_1,\ldots,x_n\}$ and $p=(p_1,\ldots,p_n)\in\KK^n$ we denote
\[
f(p)=f(p_1,\ldots,p_n)=\ev(f)(p)=\ev_p(f)\in\KK.
\]
Notice that $\ev$ is in fact a map of $\bE$-algebras; in other words, $$\ev_p:\KK\{x_1,\ldots,x_n\}\to\KK$$ is a morphism of $\bE$-algebras for each $p$. Furthermore, $\Set(\KK^n,\KK)$ has a point-wise $\bE$-algebra structure making $\ev^n$ into a $\bE$-algebra map.

The following proposition can be proved in several ways; the short proof below was suggested by E. Dubuc.

\begin{prop}\label{prop:homsintoK}
Let $\bE$ be a Fermat theory, $\KK=\bE(0)$. Then, given an arbitrary $\bE$-algebra $\A$, any $\KK$-algebra homomorphism $\phi:\A\to\KK$ is a morphism of $\bE$-algebras.
\end{prop}
\begin{proof}
Let $I=\Ker{\phi}$. On one hand, by Corollary \ref{cor:idealsarecongs}, there is a unique $\bE$-algebra structure on $\A/I$ making the projection $\pi:\A\to\A/I$ an $\bE$-algebra homomorphism. On the other hand, notice that $\phi$ is surjective (since it preserves units), hence it factors as $\phi=\tilde{\phi}\circ\pi$, where $\tilde{\phi}:\A/I\to\KK$ is a $\KK$-algebra \emph{isomorphism}. However, since $\KK$ is initial both as a $\KK$-algebra and as a $\bE$-algebra, the $\bE$-algebra structure on $\KK$ is uniquely determined by its $\KK$-algebra structure; therefore, $\tilde{\phi}$ is also an isomorphism of $\bE$-algebras. It follows immediately that $\phi$ is an $\bE$-algebra homomorphism.
\end{proof}

\begin{cor}
Any $\KK$-algebra homomorphism
\[
P:\KK\{x_1,\ldots,x_n\}\To\KK
\]
is of the form $\ev_p$ for some $p\in\KK^n$.
\end{cor}
\begin{proof}
$P$ is in fact an $\bE$-algebra homomorphism by Proposition \ref{prop:homsintoK}. The conclusion follows by observing that $\KK\{x_1,\ldots,x_n\}$ is the free $\bE$-algebra on $n$ generators.
\end{proof}




\subsubsection{Reduced Fermat theories.}
As the following examples show, the evaluation maps \eqref{eqn:ev} need not be injective.

\begin{eg}\label{eg:finitering} Let $\bE=\bcom_\KK$, where $\KK$ is a \emph{finite} ring, with elements labeled $k_1,\ldots,k_N$. Then the polynomial
\[
p(x)=(x-k_1)\cdots(x-k_N)
\]
evaluates to $0$ on every $k\in\KK$, and yet is itself non-zero (being a monic polynomial of degree $N$).
\end{eg}

It is easy to see that this phenomenon cannot occur for $\bcom_\KK$ with $\KK$ containing $\QQ$. However, the next example illustrates that it \emph{can} occur even for theories over $\QQ$.

\begin{eg}\label{eg:germsnonpd}
Consider the $\bcinf$-algebra $\KK=\cinf(\RR)_0$ of germs of smooth functions of one variable. It can be presented as the quotient of $\cinf(\RR)$ by the ideal $m_0^g$ consisting of those smooth functions $f(x)$ which vanish on some neighborhood of $0$. Let $\bE=\bcinf_\KK$. It follows (cf. \cite{msia}, p. 49) that $\KK\{y\}$ is the quotient of $\cinf(\RR^2)$ by the ideal $m^{tub}_{x=0}$ consisting of those functions $f(x,y)$ which vanish on some \emph{tubular} neighborhood of the $y$-axis (i.e. a set of the form $(-\epsilon,\epsilon)\times\RR$ for some $\epsilon>0$). If $g$ is the germ at the origin of some function $g(x)$ and $[f]$ is the class modulo $m^{tub}_{x=0}$ of some function $f(x,y)$, then $\ev_g([f])$ is the germ of $f(x,g(x))$ at the origin.

Now, let $f(x,y)$ be a smooth function whose vanishing set contains \emph{some} neighborhood of the $y$-axis but does not contain any tubular neighborhood. Then the class $[f]\in\KK\{y\}$ is non-zero, and yet $\ev_g([f])=0$ for all $g$.
\end{eg}

\begin{defn}\label{def:reduced}
A Fermat theory $\bE$ is called \emph{reduced} if all the evaluation maps \eqref{eqn:ev} are injective.
\end{defn}

Thus, reduced theories are ``theories of differentiable functions'' in the sense that $n$-ary operations are labeled by functions from $\KK^n$ to $\KK$. For instance, $\bcom_\KK$ is reduced for $\KK\supset\QQ$, as are the theories $\bcinf$, $\bcomega$, $\bH$ and $\bH_\RR$. Non-reduced theories, such as the ones in Examples \ref{eg:finitering} and \ref{eg:germsnonpd}, can be viewed as pathological in some sense. We are now going to describe a functorial procedure of turning any Fermat theory into a reduced one.

\begin{defn}\label{def:reduction}
Given a Fermat theory $\bE$, define $\bE_\red$ by setting $$\bE_\red(n,1)=\Im\ev^n=\bE(n,1)/\Ker(\ev^n).$$
\end{defn}

\begin{rem}
As a category, one can describe $\bE_\red$ as the opposite category of the full subcategory of $\bE\Alg$ on algebras of the form $\bE\left(n\right)/Ker(\ev^n).$ As a consequence, one has that the free $\bE_\red$-algebra on $n$-generators is $\bE\left(n\right)/Ker(\ev^n).$
\end{rem}

\begin{prop}\label{prop:reduction}
If $\bE$ is a Fermat theory, $\bE_\red$ is a reduced Fermat theory. Moreover, the assignment $\bE\mapsto\bE_\red$ is functorial and is left adjoint to the inclusion
\[
\mathbf{FTh}_\red\hookrightarrow\mathbf{FTh}
\]
of the full subcategory of reduced Fermat theories. The same holds with $\mathbf{FTh}$ replaced with $\mathbf{FTh}_{/\QQ}$ or $\mathbf{FTh}_\KK$
\end{prop}
\begin{proof}
First, $\bE_\red$ is in fact a theory since $\ev=\{\ev^n\}_{n\in\NN}$ is a map of algebraic theories from $\bE$ to $\End_\KK$ (see Example \ref{eg:end}), and $\bE_\red=\Im(\ev),$ as in Remark \ref{rem:image}; obviously, $\bE_\red$ is reduced.

To see that $\bE_\red$ remains a Fermat theory observe that, in the Fermat property \eqref{eqn:fermatproperty} for $\bE$, if $f\in \Ker(\ev^{n+1})$, then $\frac{\Delta f}{\Delta x}\in\Ker(\ev^{n+2})$, and vice versa.

The functoriality and adjointness are clear (the latter follows immediately from the presentation $\bE_\red=\bE/\Ker(\ev)$).

Finally, the last statement is simply the observation that reduction does not change the ground ring.
\end{proof}

\begin{eg}
If $\bE$ is as in Example \ref{eg:germsnonpd}, its reduction $\bE_\red$ can be described as follows. Let $\RR\{x,y_1,\ldots,y_n\}=\cinf(\RR^{n+1})$, the free $\bcinf$-algebra on $n+1$ generators. Then $\bE(n)=\RR\{x,y_1,\ldots,y_n\}/m^{tub}_{x=0}$ while $\bE_\red(n)=\RR\{x,y_1,\ldots,y_n\}/m^g_{x=0}$, where $m^{tub}_{x=0}$ is the ideal of functions vanishing in some tubular neighborhood of the hyperplane $x=0$, while $m^g_{x=0}$ consists of functions vanishing in some (not necessarily tubular) neighborhood of $x=0$. Clearly, $m^{tub}_{x=0}\subsetneq m^g_{x=0}$; in fact, $m^g_{x=0}$ is the \emph{germ-determined closure} of $m^{tub}_{x=0}$, hence, viewed as $\bcinf$-algebras, $\bE_\red(n)$ is the germ-determined quotient of $\bE(n)$. Specifically, $\bE_\red(n)$ consists of germs of smooth functions on $\RR^{n+1}$ at the hyperplane $x=0$. In fact, the free $\bE_\red$-algebras $\bE_\red(n)$ are formed by adjoining variables to $\KK=\cinf(\RR)_0$ using the coproduct in the category of \emph{germ-determined} $\bcinf$-algebras. We refer to \cite{msia} for the appropriate definitions and discussion.
\end{eg}

Let us conclude this section by constructing, for any ring $\KK$, the \emph{maximal} reduced Fermat theory $\bF(\KK)$ with $\KK$ as the ground ring. Indeed, the $\KK$-algebra structure on $\KK$ amounts to a map of Lawvere theories $\bcom_\KK\to\End_\KK$, and any reduced Fermat theory with ground ring $\KK$ is a subtheory of $\End_\KK$, as we have seen. $\bF(\KK)$ will be the maximal Fermat subtheory of $\End_\KK$. More precisely, we have

\begin{defn}
Let $f:\KK^n\to\KK$ be a function. Given a $k=1,\ldots,n$, we say that $f$ is \emph{differentiable in the $k$th variable} if there is a unique function $\frac{\Delta f}{\Delta x^k}:\KK^{n+1}\to\KK$ such that
\[
f(\ldots,x,\ldots)-f(\ldots,y,\ldots)=(x-y)\cdot\frac{\Delta f}{\Delta x^k}(\ldots,x,y,\ldots).
\]
This $g_k$ is then called the \emph{difference quotient} of $f$ with respect to the $k$th variable. Say that $f$ is \emph{differentiable} if it is differentiable in all the variables. Given $N>1$, say that $f$ is $N$ times differentiable if $f$ is differentiable and all its difference quotients are $N-1$ times differentiable. Finally, say that $f$ is \emph{smooth} if it is $N$ times differentiable for all $N$.
\end{defn}

\begin{propdef}
Let $\bF(\KK)(n,1)$ consist of all smooth functions $$f:\KK^n\to\KK.$$ Then $\bF$ is the maximal reduced Fermat theory with ground ring $\KK$.
\end{propdef}
\begin{proof}
To see that $\bF(\KK)$ is a subtheory of $\End_\KK$, just observe that the superposition of smooth functions is again smooth (the chain rule!). By construction, $\bF(\KK)$ is a reduced Fermat theory with ground ring $\KK$ and for any other Fermat theory $\bE$ with ground ring $\KK$, the structure map $\bE\to\End_\KK$ for the $\bE$-algebra structure on $\KK$ factors through $\bF(\KK)$.
\end{proof}

\begin{eg}\label{eg:maxfermat}
$\bF(\RR)=\bcinf$, while $\bF(\CC)=\bH$. We do not know what $\bF(\QQ)$ is but it certainly contains $\bR^\QQ$ ($=\bR^\ZZ$).
\end{eg} 
\subsection{Nilpotent extensions}



\begin{defn}
Let $\KK$ be a commutative ring, $\A$ a $\KK$-algebra. An \emph{extension} of $\A$ (over $\KK$) is a surjective $\KK$-algebra homomorphism $\pi:\A'\to\A.$ A \emph{split extension} of $\A$ (over $\KK$) is an extension $\pi:\A'\to\A$ together with a section (splitting) $\iota$ of $\pi$ which is also a $\KK$-algebra homomorphism. Such a split extension is called a \emph{split nilpotent extension} if additionally, the kernel $\Ker\left(\pi\right)$ is a nilpotent ideal.
\end{defn}

\begin{rem}
Warning: This notion of \emph{extension} should not be confused with the notion of extension used in Galois theory!
\end{rem}

\begin{rem}
Notice that any $\KK$-algebra map $$\A \to \KK$$ is automatically surjective and split in a canonical way; a section is provided by the unique $\KK$-algebra homomorphism $$\KK \to \A.$$
\end{rem}

\begin{defn}
A \emph{Weil $\KK$-algebra} is a nilpotent extension of $\KK$ (over $\KK$) which is finitely generated as a $\KK$-module. A \emph{formal Weil $\KK$-algebra} is a nilpotent extension of $\KK$ (over $\KK$).
\end{defn}

\begin{rem}\label{rem:local}
When $\KK$ is a \emph{field}, any $\fweilk$ $\A'$ has an underlying $\KK$-algebra of the form $\KK\oplus \mathfrak{m},$ with $\mathfrak{m}$ an nilpotent maximal ideal. Moreover, from the direct sum decomposition, every element of $\A'$ can be expressed uniquely as $a=k+m$ with $k \in \KK$ and $m$ nilpotent. If $k \ne 0,$ $a$ is a unit. If $k=0,$ $a$ is nilpotent. Hence, $\A'$ is a local $\KK$-algebra with the maximal ideal $\mathfrak{m}$, and residue field $\KK$. 
\end{rem}

\begin{rem}
By Nakayama's lemma, if $\KK$ is a field, it follows that $\A$ is a Weil $\KK$-algebra if and only if there exists a surjection $$\pi:\A \to \KK$$ whose kernel is a finitely generated $\KK$-module.
\end{rem}

Let $\bE$ be a Fermat theory with ground ring $\KK$.

\begin{prop}\label{prop:nilpext}
Let $\A\in\bE\Alg$, and $$\xymatrix{\A' \ar@<+0.65ex>^-{\pi}[r] & \A_\sharp \ar@/^0.8pc/^{\iota} [l]}$$ be any split nilpotent extension of $\A_\sharp$ in $\comalg_\KK$. Then there is a unique $\bE$-algebra structure on $\A'$, consistent with its commutative algebra structure and making both the projection $$\pi:\A'\to\A$$ and the splitting $\iota:\A\to\A'$ into $\bE$-algebra maps. Furthermore, for any $\bE$-algebra $\B$, we have
\begin{enumerate}
\item Any $\KK$-algebra map $\Psi:\A'\to\B$ such that the precomposition $$\psi=\Psi\circ\iota:\A\to\B$$ is a map of $\bE$-algebras, is a map of $\bE$-algebras;
\item Any $\KK$-algebra map $\Phi:\B\to\A'$ such that the composition $\phi=\pi \circ \Phi:\B\to\A$ is a split map of $\bE$-algebras, with splitting $\sigma:\A\to\B$ such that $$\Phi\circ\sigma=\iota,$$ is a map of $\bE$-algebras.
\end{enumerate}
\end{prop}
\begin{proof}
Any element of $\A'$ can be written uniquely as a sum $a'=\iota(a)+\tilde{a}$ with $a\in\A$ and $\tilde{a}\in\N=\Ker\pi$. Let $n$ be the nilpotence degree of $\N$ (so $\N^{n+1}=0$). The Taylor expansion (Corollary \ref{cor:taylor}) now provides a unique evaluation of any operation in $\bE$ on any tuple of elements of $\A'$. More precisely, let $f\in\bE(k,1)$ and
\[
a'_1=\iota(a_1)+\tilde{a}_1,\ldots,a'_k=\iota(a_k)+\tilde{a}_k\in\A'.
\]
Use the Taylor formula \eqref{eqn:taylor} to write
\[
f(\mathbf{x}+\mathbf{y})=\sum_{|\alpha|=0}^nh_\alpha(\mathbf{x})\mathbf{y}^\alpha+\sum_{|\beta|=n+1}\mathbf{y}^\beta g_\beta(\mathbf{x},\mathbf{y}),
\]
and define
\[
f(a'_1,\ldots,a'_k)=f(\mathbf{a'})=\sum_{|\alpha|=0}^n\iota(h_\alpha(\mathbf{a}))\mathbf{\tilde{a}}^\alpha.
\]
Since the Taylor expansion is compatible with compositions (the generalized chain rule), this defines an $\bE$-algebra structure on $\A'$. It is clearly compatible with its $\KK$-algebra structure and makes both $\iota$ and $\pi$ into $\bE$-algebra homomorphisms.

Now let us prove the properties $(1)$ and $(2)$ of this structure. For $(1)$, observe first that the direct image ideal $\Psi_*(\N)$ is also nilpotent of degree $n$. Using our assumptions, we have
\begin{eqnarray*}
\Psi(f(a'_1,\ldots,a'_k))&=&\Psi(\sum_{|\alpha|=0}^n\iota(h_\alpha(\mathbf{a}))\mathbf{\tilde{a}}^\alpha)
=\sum_{|\alpha|=0}^n\psi(h_\alpha(\mathbf{a}))\Psi(\mathbf{\tilde{a}})^\alpha\\
&=&\sum_{|\alpha|=0}^nh_\alpha(\psi(\mathbf{a}))\Psi(\mathbf{\tilde{a}})^\alpha
=f(\psi(\mathbf{a})+\Psi(\mathbf{\tilde{a}}))\\
&=&f(\Psi(a'_1),\ldots,\Psi(a'_k)).
\end{eqnarray*}

To prove (2), let $b\in\B$ and decompose
\[
\Phi(b)=\iota\pi\Phi(b)+\widetilde{\Phi(b)}=\iota\phi(b)+\widetilde{\Phi(b)}
\]
with
\[
\widetilde{\Phi(b)}=\Phi(b)-\iota\pi\Phi(b)\in\N.
\]
We can also decompose
\[
b=\sigma\phi(b)+\tilde{b}
\]
with
\[
\tilde{b}=b-\sigma\phi(b)\in\Ker\phi.
\]
Since $\Phi$ is a $\KK$-algebra homomorphism, we have
\[
\Phi(b)=\Phi\sigma\phi(b)+\Phi(\tilde{b})=\iota\phi(b)+\Phi(\tilde{b}).
\]
Hence,
\[
\widetilde{\Phi(b)}=\Phi(\tilde{b})\in\N.
\]
Now let $f\in\bE(k,1)$, $\mathbf{b}=(b_1,\ldots,b_k)\in\B^k$. Using Taylor's formula and our assumptions, we have
\begin{eqnarray*}
\Phi(f(\mathbf{b}))&=&\Phi(f(\sigma\phi(\mathbf{b})+\mathbf{\tilde{b}}))\\
&=&\Phi(\sum_{|\alpha|=0}^nh_\alpha(\sigma\phi(\mathbf{b}))\mathbf{\tilde{b}}^\alpha+\sum_{|\beta|=n+1}\mathbf{\tilde{b}}^\beta g_\beta(\sigma\phi(\mathbf{b}),\mathbf{\tilde{b}}))\\
&=&\sum_{|\alpha|=0}^n\Phi\sigma\phi(h_\alpha(\mathbf{b}))\Phi(\mathbf{\tilde{b}})^\alpha+\sum_{|\beta|=n+1}\Phi(\mathbf{\tilde{b}})^\beta \Phi(g_\beta(\sigma\phi(\mathbf{b}),\mathbf{\tilde{b}}))\\
&=&\sum_{|\alpha|=0}^n\iota\phi(h_\alpha(\mathbf{b}))\widetilde{\Phi(\mathbf{b})}^\alpha+\sum_{|\beta|=n+1}\widetilde{\Phi(\mathbf{b})}^\beta \Phi(g_\beta(\sigma\phi(\mathbf{b}),\mathbf{\tilde{b}}))\\
&=&\sum_{|\alpha|=0}^n\iota(h_\alpha(\phi(\mathbf{b})))\widetilde{\Phi(\mathbf{b})}^\alpha\\
&=&f(\iota\phi(\mathbf{b})+\widetilde{\Phi(\mathbf{b})})=f(\Phi(\mathbf{b})).
\end{eqnarray*}
\end{proof}

This has the following important consequences:

\begin{cor}\label{cor:weilisgood}
Let $\W$ be a $\fweilk.$ Then there is a unique $\bE$-algebra structure on $\W$ consistent with its $\KK$-algebra structure. Furthermore, it has the following properties:
\begin{enumerate}
\item Given an arbitrary $\bE$-algebra $\A$, any $\KK$-algebra homomorphism $\A\to\W$ or $\W\to\A$ is an $\bE$-algebra homomorphism;
\item The algebraic tensor product $\A\otimes\W$ (over $\KK$) coincides with the coproduct $\A\odot\W$ of $\bE$-algebras;
\item The tensor product of finitely many $\fweil$ $\KK$-algebras is again a $\fweilk.$ Similarly, the tensor product of finitely many Weil algebras is again a Weil algebra.
\end{enumerate}
\end{cor}
\begin{proof}
To see that $\W$ supports a unique $\bE$-algebra structure making any $\KK$-algebra map to or from $\W$ a homomorphism of $\bE$-algebras, we invoke Propositions \ref{prop:nilpext} and \ref{prop:homsintoK}.

To see that $\A\otimes\W$ is the coproduct in $\bE$, observe first that since $A \otimes \left(\quad\right)$ is a functor, $\A\otimes\W$ is a split extension of $\A$. Moreover, its kernel may naturally be identified with $\A \otimes \N=\left(i_\W\right)_*\left(\N\right),$ where $\N$ is the kernel of the extension $\W \to \KK,$ and $$i_\W:\W \to \A \otimes \W$$ is the canonical map. It follows that this kernel has nilpotency degree equal to that of $\N,$ so that $$\A\otimes\W \to \A$$ is a split nilpotent extension. Hence $\A\otimes\W$ supports a unique $\bE$-algebra structure making the canonical inclusions from $\A$ and $\W$ into $\bE$-algebra homomorphisms. Now, suppose we are given an $\bE$-algebra $\B$ and maps of $\bE$-algebras $f:\A\to\B$ and $g:\W\to\B$. Then we get a unique $\KK$-algebra map $f\otimes g:\A\otimes\W\to\B$ extending $f$ and $g$. But then $f\otimes g$ is an $\bE$-algebra map by Proposition \ref{prop:nilpext}.

For $(3)$, notice that if
$$\xymatrix{\W\cong \KK\oplus \mathfrak{m} \ar@<+0.65ex>^-{\pi}[r] & \KK \ar@/^1.2pc/^{\sigma} [l]}$$
and
$$\xymatrix{\W'\cong \KK\oplus \mathfrak{m'} \ar@<+0.65ex>^-{\pi'}[r] & \KK \ar@/^1.2pc/^{\sigma'} [l]}$$
are nilpotent extensions, then
$$\xymatrix@C=2.5cm{\W'\otimes \W \ar@<+0.65ex>^-{\pi'':=\pi' \circ \left(id_{\W'} \otimes \pi\right)}[r] & \KK \ar@/^1.2pc/^{\sigma'':=\left(id_{\W'} \otimes \sigma\right) \circ \sigma'} [l]}$$
is a split extension of $\KK$ with $$\Ker\left(\pi''\right) \cong \mathfrak{m} \oplus \mathfrak{m'} \oplus \left(\mathfrak{m}\otimes\mathfrak{m'}\right),$$ which is finitely generated as a $\KK$-module if both $\mathfrak{m}$ and $\mathfrak{m'}$ are. This implies $(3)$ holds for Weil algebras. For the case of general nilpotent extensions, note that if $$i_\W':\W' \hookrightarrow \W' \otimes \W$$ and $$i_\W:\W \hookrightarrow \W' \otimes \W,$$ are the canonical maps, then $\Ker\left(\pi''\right)=\left(i_\W'\right)_{*}\left(\mathfrak{m'}\right)+\left(i_\W\right)_{*}\left(\mathfrak{m}\right).$ If $\mathfrak{m}$ is nilpotent of degree $m$ and $\mathfrak{m'}$ is nilpotent of degree $n$ and, then it follows that $\Ker\left(\pi''\right)$ is nilpotent of degree $m+n-1.$
\end{proof}

\begin{cor}\label{cor:weilcomp}
For any $\fweilk$ $\W,$ the co-unit $$\widehat{\W_\sharp} \to \W$$ is an isomorphism.
\end{cor}

\begin{proof}
This follows immediately from $(1)$ of \ref{cor:weilisgood}.
\end{proof}

\section{Super Fermat Theories}\label{sec:superfermat}
\subsection{Superalgebras and superizations}\label{sec:supalg}

\begin{defn}\label{defn:supcom}
Let $\KK$ be a commutative ring. A \emph{supercommutative superalgebra over $\KK$} (or supercommutative algebra) is a $\ZZ_2$-graded associative unital $\KK$-algebra $$\A=\{\A_\even,\A_\odd\}$$ such that $\A_\even$ is commutative and for every $a\in\A_\odd$,
\[
a^2=0.
\]
 We say that $a$ is of (Grassman) \emph{parity} $\epsilon$ if $a\in\A_\epsilon$; we say it is \emph{even} (resp. \emph{odd}) if $\epsilon=\even$ (resp. $\epsilon=\odd$). Supercommutative superalgebras over $\KK$ form a category, denoted by $\bscom_\KK\Alg$, whose morphisms are parity-preserving $\KK$-algebra homomorphisms.
\end{defn}

\begin{rem} The definition implies that
\[
a_1a_2=(-1)^{\epsilon_1\epsilon_2}a_2a_1
\]
whenever $a_i\in\A_{\epsilon_i}$, $i=1,2$, justifying the term ``supercommutative''; if $\frac{1}{2}\in\KK$, the converse also holds.
\end{rem}

\begin{rem}
In our formulation, there is no such thing as elements of mixed parities in a super commutative algebra $\A,$ as $\A$ lacks an underlying set. Instead, it has an underlying $\ZZ_2$-graded set, that is a set of two sets $\{\A_\even,\A_\odd\}$. Consequently, if $a \in \A_\even$ and $b \in \A_\odd,$ the expression $a +b$ has no meaning. The advantage of this treatment is that it behaves nicely with respect to the fact that supercommutative algebras are algebras for a $2$-sorted Lawvere theory. This viewpoint is not essential, as the category of super commutative algebras as defined in Definition \ref{defn:supcom} is canonically equivalent to the category non-commutative $\KK$-algebras $\A$ together with a grading $\A=\A_\even \oplus \A_\odd,$ making $\A$ supercommutative, where the morphisms are algebra morphisms respecting the grading. If one would like, one may work entirely within the framework of uni-sorted Lawvere theories as in \cite{yetter}, but this makes things unnecessarily complicated and yields less flexibility.
\end{rem}

There is a forgetful functor
\[
u^*:\bscom_\KK\Alg\To\Set^{\{\even,\odd\}}
\]
to the category of $\ZZ_2$-graded sets, whose left adjoint $u_!$ assigns to a pair of sets $P=(P_\even|P_\odd)$ the free supercommutative superalgebra on the set $P_\even$ of even and the set $P_\odd$ of odd generators.

\begin{defn}\label{dfn:grass}
Given a $\KK$-algebra $\R$, the \emph{Grassmann (or exterior) $\R$-algebra on $n$ generators} is the free supercommutative $\R$-superalgebra on $n$ odd generators. In other words, it is generated as an $\R$-algebra by odd elements $\xi^1,\ldots,\xi^n$ subject to relations
\[
\xi^i\xi^j+\xi^j\xi^i=0.
\]
Denote this algebra by $\Lambda^n_\R$ (or simply $\Lambda^n$ if $\R=\KK$).
\end{defn}

\begin{rem}
$(\Lambda^n_\R)_\even$ is a Weil $\R$-algebra.
\end{rem}

It is easy to see that the free supercommutative $\KK$-superalgebra on $m$ even and $n$ odd generators is nothing but $\Lambda_\R^n$ with $\R=\KK[x^1,\ldots,x^m]$. Denote this algebra by
\[
\KK[x^1,\ldots,x^m;\xi^1,\ldots,\xi^n].
\]
Supercommutative superalgebras over $\KK$ are algebras over a $2$-sorted Lawvere theory $\bscom_\KK$, which we now describe. As a category, $\bscom_\KK$ is equivalent to the opposite of the category of finitely generated supercommutative $\KK$-superalgebras:
\[
\bscom_\KK(m|n)=\KK[x^1,\ldots,x^m;\xi^1,\ldots,\xi^n].
\]
It is generated by the set $\{\even,\odd\}$ of Grassmann parities; the product of $m$ copies of $\even$ and $n$ copies of $\odd$ will be denoted by $(m|n)$. The morphisms are
\[
\bscom_\KK((m|n),(p|q))=\bscom_\KK(m|n)_\even^p\times\bscom_\KK(m|n)_\odd^q,
\]
and the composition is defined by substitution. Notice that the ground ring of $\bscom_\KK$ is $\bscom(0|0)=\KK$.

Observe that $\bcom_\KK$ sits inside $\bscom_\KK$ as the full subcategory of ``purely even'' objects of the form $(m|0)$, $m\in\NN$. Clearly, the embedding
\[
\iota:\bcom_\KK\To\bscom_\KK, \quad m\mapsto(m|0),
\]
is a morphism of algebraic theories, hence induces an adjunction
\[
\Adj{\iota^*}{\bcom_\KK\Alg}{\bscom_\KK\Alg}{\iota_!},
\]
such that $\iota^*\A=\A_\even$, while $\iota_!\A=\{\A,0\}$ (the superalgebra with even part equal to $\A$ and trivial odd part).

We now observe that the $2$-sorted Lawvere theory $\bscom_\KK$ satisfies a natural generalization of the Fermat property:

Suppose that $f$ is an element of $\KK[x,z^1,\ldots,z^m;\xi^1,\ldots,\xi^n].$ Then $f$ can be expressed uniquely in the form $$f\left(x,\mathbf{z},\mathbf{\xi}\right)=\sum\limits_{I \subset \set{1,\ldots,n}} f^I\! \mathbf{\xi}^I,$$ where if $I=\set{i_1,\ldots i_k},$ $$\mathbf{\xi}^I=\xi^{i_1} \ldots \xi^{i_k},$$ with each $f^I\!$  in $\KK[x,z^1,\ldots,z^m].$ (If $f$ is even, $f^I=0$ for all $I$ with odd cardinality, and vice-versa for $f$ odd.) By the Fermat property for $\bcom_\KK,$ for each $$I \subset \set{1,\ldots,n},$$ there is a unique $\frac{\Delta f^I}{\Delta x}\!\left(x,y,\mathbf{z}\right) \in \KK[x,y,z^1,\ldots,z^m],$ such that $$f^I\!\left(x,\mathbf{z}\right)-f^I\!\left(y,\mathbf{z}\right)=\left(x-y\right) \cdot \frac{\Delta f^I}{\Delta x}\!\left(x,y,\mathbf{z}\right).$$ Let $$\frac{\Delta f}{\Delta x}\left(x,y,\mathbf{z},\mathbf{\xi}\right):=\sum\limits_{I \subset \set{1,\ldots,n}} \frac{\Delta f^I}{\Delta x}\!\left(x,y,\mathbf{z}\right) \mathbf{\xi}^I \in \KK[x,y,z^1,\ldots,z^m;\xi^1,\ldots,\xi^n].$$ Then we have that
\begin{equation}
f\left(x,\mathbf{z},\mathbf{\xi}\right)-f\left(y,\mathbf{z},\mathbf{\xi}\right)=\left(x-y\right) \cdot \frac{\Delta f}{\Delta x}\left(x,y,\mathbf{z},\mathbf{\xi}\right).
\end{equation}
Moreover, it is not hard to see that $\frac{\Delta f}{\Delta x}\left(x,y,\mathbf{z},\mathbf{\xi}\right)$ is unique with this property.

Suppose now that the role of $x$ and $y$ are played by odd generators. That is, suppose $f \in \KK[x^1,\ldots,x^m;\eta,\xi^1,\ldots,\xi^n].$  Then since $\eta^2=0,$ $f$ can be uniquely expressed in the form
$$f\left(\mathbf{x},\eta,\mathbf{\xi}\right)=\sum\limits_{I \subset \set{1,\ldots,n}} h^I\! \mathbf{\xi}^I + \eta \cdot \left(\sum\limits_{I \subset \set{1,\ldots,n}} g^I\! \mathbf{\xi}^I\right),$$ with $h^I$ and $g^I$ in $\KK[x^1,\ldots,x^m].$ Let $$h\left(\mathbf{x},\mathbf{\xi}\right):=\sum\limits_{I \subset \set{1,\ldots,n}} h^I\! \mathbf{\xi}^I$$ and
$$g\left(\mathbf{x},\mathbf{\xi}\right):=\sum\limits_{I \subset \set{1,\ldots,n}} g^I\! \mathbf{\xi}^I.$$
Then we have $$f\left(\mathbf{x},\eta,\mathbf{\xi}\right)=h\left(\mathbf{x},\mathbf{\xi}\right) + \eta \cdot g\left(\mathbf{x},\mathbf{\xi}\right).$$
Notice that $h$ is the value of $f$ at $\eta=0$ while $g$ is the (left) partial derivative  $\frac{\del f}{\del \eta}$ of $f$ with respect to $\eta$. Furthermore, we have the following:
\begin{equation}
f\left(\mathbf{x},\eta,\mathbf{\xi}\right)-f\left(\mathbf{x},\theta,\mathbf{\xi}\right)=\left(\eta-\theta\right) \cdot g\left(\mathbf{x},\mathbf{\xi}\right).
\end{equation}
Regarding $g\left(\mathbf{x},\mathbf{\xi}\right)$ as $g\left(\mathbf{x},\eta,\theta,\mathbf{\xi}\right) \in \KK[x^1,\ldots,x^m;\eta,\theta,\xi^1,\ldots,\xi^n],$ we have that
\begin{equation}\label{sfermmm}
f\left(\mathbf{x},\eta,\mathbf{\xi}\right)-f\left(\mathbf{x},\theta,\mathbf{\xi}\right)=\left(\eta-\theta\right) \cdot g\left(\mathbf{x},\eta,\theta,\mathbf{\xi}\right).
\end{equation}
Note however that $g\left(\mathbf{x},\eta,\theta,\mathbf{\xi}\right)$ is \emph{not} unique with this property; one could also use $$g\left(\mathbf{x},\mathbf{\xi}\right)+\left(\eta-\theta\right) \cdot p\left(\mathbf{x},\mathbf{\xi}\right)$$ for any $p$. However, by differentiating (\ref{sfermmm}) with respect to $\eta$ and $\theta,$ one sees immediately that there is a unique such $g\left(\mathbf{x},\eta,\theta,\mathbf{\xi}\right)$ such that $$\frac{\del g}{\del \eta}=\frac{\del g}{\del \theta}=0,$$ in other words there exists a unique $g$ which is only a function of $\mathbf{x}$ and $\mathbf{\xi}$.

This motivates the following definitions:

\begin{defn}
Let $$\tau_{\bS}:\bscom \to \bS$$ be an extension of $\bscom$ as a $2$-sorted Lawvere theory, where implicitly $$\bscom=\bscom_{\ZZ}.$$ Without loss of generality, assume the objects are given by pairs $\left(m|n\right)$ with $m$ and $n$ non-negative integers, such that $\tau_\bS$ is the identity on objects when $\bscom$ is equipped with the usual sorting. Denote by $$\bS\left(0|0\right)=:\KK$$ the initial $\bS$-algebra. The free $\bS$-algebra $\bS\left(m|n\right)$ is called \textit{the free $\bS$-algebra on $m$ even and $n$ odd generators}, and is denoted by $$\KK\{x^1,\ldots x^m;\xi^1,\ldots,\xi^n\}.$$
\end{defn}

\begin{defn}\label{dfn:superfermat}
An extension $\bS$ of $\bscom$ is called a \emph{super Fermat theory} if for every $$f\in\KK\{x,z_1,\ldots,z_n;\xi^1,\ldots,\xi^n\}$$ there exists a unique $\frac{\Delta f}{\Delta x}\in\KK\{x,y,z_1,\ldots,z_n;\xi^1,\ldots,\xi^n\}$ such that
\begin{equation}\label{eqn:sfermatproperty1}
f(x,\mathbf{z},\mathbf{\xi})-f(y,\mathbf{z},\mathbf{\xi})=(x-y)\cdot \frac{\Delta f}{\Delta x}(x,y,\mathbf{z},\mathbf{\xi}),
\end{equation}
and for every $\varphi \in \KK\{x^1,\ldots,x^m;\eta,\xi^1,\ldots,\xi^n\}$ there exists a unique $$\frac{\Delta\varphi}{\Delta\eta} \in \KK\{x^1,\ldots,x^m;\xi^1,\ldots,\xi^n\}$$ such that
\begin{equation}\label{eqn:sfermatproperty2}
\varphi\left(\mathbf{x},\eta,\mathbf{\xi}\right)-\varphi\left(\mathbf{x},\theta,\mathbf{\xi}\right)=\left(\eta-\theta\right) \cdot \frac{\Delta\varphi}{\Delta\eta}\left(\mathbf{x},\mathbf{\xi}\right).
\end{equation}
Denote by $\mathbf{SFTh}$ the full subcategory of $2$-sorted Lawvere theories under $\bscom$ consisting of those which are super Fermat theories.
\end{defn}

\begin{propdef}
Let $\bE$ be a Fermat theory with ground ring $\KK$. There exists an algebraic theory $\bSE$, the \emph{superization} of $\bE$, with the set $\{\even,\odd\}$ of Grassmann parities as sorts, and with operations given by
\[
\bSE((m|n),(p|q))=\bSE(m|n)_\even^p\times\bSE(m|n)_\odd^q,
\]
where
\[
\bSE(m|n)=\bE(m)\otimes_\KK\Lambda^n.
\]
An \emph{$\bE$-superalgebra} is a $\bSE$-algebra.
\end{propdef}
\begin{proof} The only thing to check is that the composition by substitution is well-defined. To this end we observe that, as $\Lambda^n_\even$ is a Weil algebra, $$\bSE(m|n)_\even=\bE(m)\otimes\Lambda^n_\even$$ has a canonical $\bE$-algebra structure by Corollary \ref{cor:weilisgood}. Now let
\[
f=f(x^1,\ldots,x^p;\xi^1,\ldots,\xi^q)=\mathop{\sum_{k\geq0}}_{i_1<\cdots<i_k}f_{i_1\ldots i_k}(x^1,\ldots,x^p)\xi^{i_1}\cdots\xi^{i_k}
\]
be an element of $\bSE(p|q)$ (with $f_{i_1\ldots i_k}\in\bE(p)=\KK\{x^1,\ldots,x^p\}$).\\ Let $g^1,\ldots,g^p\in\bSE(m|n)_\even$, and $\gamma^1,\ldots,\gamma^q\in\bSE(m|n)_\odd.$ It follows that
\[
f(g^1,\ldots,g^p;\gamma^1,\ldots,\gamma^q)=\mathop{\sum_{k\geq0}}_{i_1<\cdots<i_k}f_{i_1\ldots i_k}(g^1,\ldots,g^p)\gamma^{i_1}\cdots\gamma^{i_k}
\]
is a well-defined element of $\bSE(m|n)$, of the same parity as $f$.

\end{proof}


\begin{notation}
Let
\[
\KK\{x^1,\ldots,x^m;\xi^1,\ldots,\xi^n\}
\]
denote the free $\bSE$-algebra on $m$ even generators $x^1,\ldots,x^n$ and $n$ odd generators $\xi^1,\ldots,\xi^n$.
\end{notation}

\begin{prop}\label{prop:superization}
If $\bE$ is a Fermat theory, its superization $\bSE$ is a super Fermat theory.
\end{prop}
\begin{proof}
The proof is nearly identical to the proof of the super Fermat property for $\bscom_\KK,$ so we leave it to the reader.

\end{proof}

\begin{eg}
Let $\bE=\bcinf$ (Example \ref{eg:cinf}). Its superization is the theory $\bS\bcinf$ of $\bcinf$-superalgebras. The free $\bcinf$-superalgebra $\RR\{x^1,\ldots,x^m;\xi^1,\ldots,\xi^n\}$ on $m$ even and $n$ odd generators is known as the \emph{Berezin algebra}. It is often denoted by $\cinf(\RR^{m|n})$ and thought of as the superalgebra of smooth functions on the $(m|n)$-dimensional Euclidean supermanifold $\RR^{m|n}$. Thus, $\bS\bcinf$ is the category of real finite-dimensional Euclidean supermanifolds and parity-preserving smooth maps between them.

Suppose that $\mathcal{M}$ is a smooth supermanifold. It induces a functor
\begin{eqnarray*}
\cinf\left(\mathcal{M}\right):\bS\bcinf &\to& \Set\\
\RR^{m|n} &\mapsto& \Hom\left(\mathcal{M},\RR^{m|n}\right),
\end{eqnarray*}
which preserves finite products. This $\bS\bcinf$-algebra is the \emph{$\bcinf$-superalgebra of smooth functions on $\mathcal{M}$}. By construction, its even elements correspond to smooth functions into $\RR$ in the traditional sense: $$\mathcal{M} \to \RR=\RR^{1|0},$$
whereas its odd elements correspond to smooth functions into the odd line: $$\mathcal{M} \to \RR^{0|1}.$$ The underlying supercommutative $\RR$-algebra of $\cinf\left(\mathcal{M}\right)$ is the global sections of its structure sheaf. More generally, the structure sheaf of any smooth supermanifold $\mathcal{M}$ is in fact canonically a sheaf of $\bcinf$-superalgebras.
\end{eg}

We will now give a more categorical description of superization. Let $\bS$ be any $2$-sorted Lawvere theory. Denote by $\bS_\even$ the full subcategory on the objects of the form $\left(n|0\right).$ Notice that $\bS_\even$ is generated under finite products by $\left(1|0\right),$ so $\bS_\even$ is a Lawvere theory. As the notation suggests, the free $\bS_\even$-algebra on $n$ generators has underlying set
\begin{eqnarray*}
\Hom_{\bS_\even}\left(n,1\right) &=& \Hom_\bS\left(\left(n|0\right),\left(1|0\right)\right)\\
&\cong& \bS\left(\left(n|0\right)\right)_\even.
\end{eqnarray*}
This produces a functor
$$\left(\quad\right)_\even:\LTh_{\set{\even,\odd}} \to \LTh$$ from $2$-sorted Lawvere theories to Lawvere theories.
By abuse of notation there is an induced functor
$$\left(\quad\right)_\even:\bscom/\LTh_{\set{\even,\odd}} \to \bcom/\LTh.$$
Notice that if $\bS$ is a super Fermat theory, then (\ref{eqn:sfermatproperty1}) implies that $\bS_\even$ is a Fermat theory. Hence there is furthermore an induced functor
$$\left(\quad\right)_\even:\mathbf{SFTh} \to \mathbf{FTh}$$
from super Fermat theories to Fermat theories.

On one hand, the superization of $\bcom_\KK$ is obviously $\bscom_\KK$ so we have a map of theories $\bscom_\KK\to\bSE$ induced by the structure map $\bcom_\KK\to\bE$. On the other hand, we also have a fully faithful embedding $\bE\to\bSE$ sending $m$ to $(m|0)$. The diagram
\begin{equation}\label{eq:pushout}
\xymatrix{\bcom_\KK \ar[r] \ar[d] & \bE \ar[d]\\
\bscom_\KK \ar[r] & \bS\bE}
\end{equation}
commutes. Therefore, we have a map of theories
\[
\phi:\bC=\bscom_\KK\coprod_{\bcom_\KK}\bE\To\bSE
\]

\begin{prop}
The map $\phi$ is an isomorphism.
\end{prop}
\begin{proof}
Suppose that $\bT$ is an algebraic theory fitting into a commutative diagram of morphisms of theories
$$\xymatrix{\bcom_\KK \ar[r] \ar[d] & \bE \ar[d]^-{\theta}\\
\bscom_\KK \ar[r]^-{\varphi} & \bT.}$$
Denote by $\left(N|M\right)$ the image $\varphi\left(n|m\right)$ in $\bT.$ (Since $\theta$ and $\varphi$ do not necessarily preserve generators, these need not be unique objects.) The functor $\theta$ induces a map of $\bE$-algebras
$$\bE\left(n\right) \to \theta^*\bT\left(N|0\right),$$ where $\theta^*\bT\left(N|0\right)$ denotes the underlying $\bE$-algebra corresponding to $\left(N|0\right)$ under the identification of $\bT^{op}$ with finitely generated $\bT$-algebras. With similar notational conventions, $\varphi$ induces a map of supercommutative algebras $$\bscom_\KK\left(0|m\right)=\Lambda^m \to \varphi^*\bT\left(0|M\right).$$ Since there are canonical $\bT$-algebra maps from $\bT\left(N|0\right)$  and $\bT\left(0|M\right)$ to $\bT\left(N|M\right),$ there is an induced map of supercommutative algebras $$\bE\left(n\right) \otimes \Lambda^m \to \bT\left(N|M\right)_\sharp.$$ These algebra maps assemble into a finite product preserving functor $\bS\bE \to \bT$ making the diagram commute. It is easy to see that this functor is unique with this property.
\end{proof}

As a corollary, we get a categorical description of superization:

\begin{cor}\label{cor:superisadj}
The functor $$\bS:\mathbf{FTh} \to \mathbf{SFTh}$$ is left adjoint to $$\left(\quad\right)_\even:\mathbf{SFTh} \to \mathbf{FTh}.$$
\end{cor}

\begin{rem}
For any super Fermat theory $\bF,$ the obvious diagram
$$\xymatrix{\bcom_\KK \ar[r] \ar[d] & \bF_\even \ar[d]\\
\bscom_\KK \ar[r] & \bF}$$
commutes, and the induced map $\bS\bF_\even \to \bF$ is the co-unit of the adjunction. The unit of $\bS \dashv \left(\quad\right)_\even$ is always an isomorphism. Hence $\mathbf{FTh}$ is a coreflective subcategory of $\mathbf{SFTh}.$ In particular, $\bS$ is full and faithful.
\end{rem}

\begin{cor}
An $\bE$-superalgebra is a superalgebra with an additional $\bE$-algebra structure on its even part; a morphism of $\bE$-superalgebras is a morphism of superalgebras whose even component is a morphism of $\bE$-algebras.
\end{cor}

\begin{rem}
For \emph{any} extension $\bE$ of $\bcom$ (not necessarily Fermat), we could simply \emph{define} $\bSE$ to be the pushout $\bC$ (\ref{eq:pushout}). However, this notion would not be very useful since one would generally have too few interesting examples of $\bSE$-algebras unless $\bE$ was Fermat. For instance, if $\bE=\mathbf{C^k}$ for some $k<\infty$ (Example \ref{eg:Ck}), even the Grassmann algebras $\Lambda^n$ are not $\bE$-superalgebras for $n$ sufficiently larger than $k$.
\end{rem}

One can draw the same conclusions from this ``super'' Fermat property as we did from the Fermat property \eqref{eqn:fermatproperty}. For instance, we have

\begin{thm}
Let $\bS$ be a super Fermat theory, $\A\in\bS\Alg$, $I=\{I_\even,I_\odd\}$ a homogeneous ideal. Then $I$ induces an $\bS$-congruence on $\A$, so that the superalgebra $\A/I$ is canonically an $\bS$-algebra and the projection $\A\to\A/I$ is an $\bS$-algebra map.
\end{thm}

\begin{rem}
The ground ring $\KK=\bS(0|0)$ of a super Fermat theory $\bS$ is generally a superalgebra, whereas for $\bS=\bSE$ it is an algebra (i.e. has trivial odd part). This indicates that not all super Fermat theories arise as superizations of Fermat theories.
\end{rem}

\begin{eg}\label{eg:relativesuper}
Let $\A\in\bSE\Alg$. Then the theory $\bSE_\A$ of $\A$-algebras enjoys the super Fermat property. If $\A$ has trivial odd part, $\bSE_\A=\bS(\bE_\A)$; otherwise, $\bSE_\A$ is not the superization of any Fermat theory.
\end{eg}

Lastly, we comment on two ways of turning a superalgebra into an algebra. We already mentioned the inclusion of theories
\[
\iota:\bE\To\bSE,\quad m\mapsto(m|0)
\]
inducing the adjunction $(\iota^*\vdash\iota_!)$:
\[
\Adj{\iota^*}{\bE\Alg}{\bSE\Alg}{\iota_!},
\]
with $\iota^*$ taking an $\bE$-superalgebra $\A$ to its even part $\A_\even$, while $\iota_!$ takes an $\bE$-algebra $\B$ to the $\bE$-superalgebra $\{\B,0\}$.

Now, observe that the inclusion $\iota$ has a right adjoint
$$
\pi:\bSE\to\bE,
$$
defined on objects by $\pi(m|n)=m$ for all $n$, and on morphisms by setting all the odd generators to $0$. Since $\pi$ is a right adjoint, it preserves products and is, therefore, a morphism of algebraic theories (though not of Lawvere theories, as it fails to preserve generators). Hence it induces an adjunction $\pi^*\vdash\pi_!$:
\[
\Adj{\pi^*}{\bSE\Alg}{\bE\Alg}{\pi_!}
\]
Moreover, $\pi^*$ and $\iota_!$ are naturally isomorphic (Remark \ref{rem:rightadj}), so the inclusion $\iota_!$ of algebras into superalgebras has also a left adjoint, $\pi_!$, sending a superalgebra $\A$ to the algebra $\A_{\mathrm{rd}}=\A/(\A_\odd)=\A_\even/(\A_\odd)^2$. Here, $(\A_\odd)$ denotes the homogeneous ideal generated by $\A_\odd$, namely, $(\A_\odd)_\odd=\A_\odd$, $(\A_\odd)_\even=\A_\odd^2$.

Observe that $\A_{\mathrm{rd}}$ is generally different from $\A_{\mathrm{red}}$ obtained by setting all the nilpotents to $0$, since $\A_\even$ may contain nilpotent elements which are not products of odd elements. Moreover, although each element in $\A_\odd$ is nilpotent, the ideal $\A_\odd$ is not, unless $\A$ is finitely generated as an $\A_\even$-algebra: otherwise, one can have non-vanishing products of arbitrarily many different odd elements.
\subsection{Nilpotent Extensions of Superalgebras}
The concepts of split nilpotent extensions, and of Weil algebras, generalize readily to the setting of supercommutative algebras:

\begin{defn}
Let $\KK$ be a supercommutative ring, $\A$ a $\KK$-algebra. A \emph{split nilpotent extension} of $\A$ is a $\KK$-algebra $\A'$ together with a surjective $\KK$-algebra homomorphism $\pi:\A'\to\A$ such that $N=\Ker\pi$ is a nilpotent ideal, and a section (splitting) $\sigma$ of $\pi$ which is also a $\KK$-algebra homomorphism.

A \emph{(super) Weil $\KK$-algebra} is an extension of $\KK$ which is finitely generated as a $\KK$-module.
\end{defn}

\begin{rem}
When $\KK$ is a purely even algebra (for instance, a field), the phrase ``Weil $\KK$-algebra'' is ambiguous as it could either mean a Weil algebra when viewing $\KK$ as a commutative algebra, a Weil algebra viewing $\KK$ as a supercommutative algebra. We shall always mean the latter, and if we need to distinguish, we will call the former a \emph{purely even} Weil $\KK$-algebra. In this context, by Nakayama's lemma, any Weil algebra is a split nilpotent extension, but the converse is false.
\end{rem}

\begin{rem}\label{rem:localsuper}
When $\KK$ is a field, any nilpotent extension $\A'$ of $\KK$ has an underlying $\KK$-algebra of the form $\KK\oplus \mathfrak{m},$ with $\mathfrak{m}$ an nilpotent maximal ideal, and $\A'$ is a local $\KK$-algebra with unique maximal ideal $\mathfrak{m},$ and residue field $\KK.$ Moreover, $\mathfrak{m}$ must contain $\A_\odd,$ since $\KK$ is purely even.
\end{rem}

Let $\bE$ be a super Fermat theory with ground ring $\KK$. Proposition \ref{prop:nilpext}, and its proof readily generalizes to the supercommutative case:

\begin{prop}\label{prop:nilpextsup}
Let $\A\in\bE\Alg$, $\pi:\A' \to \A_\sharp$ any split nilpotent extension of $\A_\sharp$ in $\comsalg_\KK$. Then there is a unique $\bE$-algebra structure on $\A'$, consistent with its supercommutative algebra structure and making both the projection $$\pi:\A'\to\A$$ and the splitting $\sigma:\A\to\A'$ into $\bE$-algebra maps. Furthermore, for any $\bE$-algebra $\B$, we have
\begin{enumerate}
\item Any $\KK$-algebra map $\Psi:\A'\to\B$ such that the precomposition $$\psi=\Psi\circ\iota:\A\to\B$$ is a map of $\bE$-algebras, is a map of $\bE$-algebras;
\item Any $\KK$-algebra map $\Phi:\B\to\A'$ such that the composition $\phi=\pi \circ \Phi:\B\to\A$ is a split map of $\bE$-algebras, with splitting $\sigma:\A\to\B$ such that $$\Phi\circ\sigma=\iota,$$ is a map of $\bE$-algebras.
\end{enumerate}
\end{prop}

\begin{cor}\label{cor:weilisgoodsup}
Let $\W$ be a split nilpotent extension of $\KK$ in $\comsalg_\KK.$ Then there is a unique $\bE$-algebra structure on $\W$ consistent with its super $\KK$-algebra structure. Furthermore, it has the following properties:
\begin{enumerate}
\item Given an arbitrary $\bE$-algebra $\A$, any super $\KK$-algebra homomorphism $\A\to\W$ or $\W\to\A$ is an $\bE$-algebra homomorphism;
\item The algebraic tensor product $\A\otimes\W$ (over $\KK$) coincides with the coproduct $\A\odot\W$ of $\bE$-algebras;
\item The tensor product of finitely many (super) Weil algebras is again a Weil algebra.
\end{enumerate}
\end{cor}

\begin{cor}\label{cor:cococococo}
For any $\fweilk$ $\W,$ the co-unit $$\widehat{\W_\sharp} \to \W$$ is an isomorphism.
\end{cor}

\begin{rem}\label{rem:locallynil}
Proposition \ref{prop:nilpextsup}, Corollary \ref{cor:weilisgoodsup}, and Corollary \ref{cor:cococococo} (as well as Proposition \ref{prop:nilpext}, Corollary \ref{cor:weilisgood}, and Corollary \ref{cor:weilcomp}) remain valid for a larger class of examples. One can define a \emph{locally} nilpotent extension in the same way as a nilpotent extension, with the role of nilpotent ideals generalized to locally nilpotent ideals. Recall that an ideal $I$ is locally nilpotent if every finitely generated subideal of $I$ is nilpotent. This is equivalent to asking for each element of the ideal $I$ to be nilpotent. The reason for this is that the operations of $\bE$ are \emph{finitary}; therefore, to evaluate an operation on a finite tuple of elements, we need only use the Taylor expansion up to the nilpotence order of the subideal generated by their nilpotent parts, rather than of the whole ideal, which may be infinite.

An important example of a locally nilpotent but not globally nilpotent extension is an infinitely generated Grassmann algebra. 
\end{rem}

\subsection{Some constructions.}\label{sec:constr}
\subsubsection{Ideals.} Recall that, given a supercommutative superalgebra $\A$ and homogeneous ideals $I_1,\ldots,I_r$ of $\A$, we can form their sum $\Sigma_k I_k$, product $\prod_k I_k$ and intersection $\bigcap_k I_k$, which are again ideals of $\A$. Two ideals $I,J\subset\A$ are called \emph{coprime} if $I+J=(1)$.

The following is a standard fact from commutative algebra (see \cite{atmcd}, the proof given there carries over verbatim to the super case).

\begin{prop} \label{prop:coprimeideals}
Let $\A_k=\A/I_k$, $k=1,\ldots,r$, let $\phi_k:\A\to\A_k$ be the canonical projections and let
\[
\phi=(\phi_1,\ldots,\phi_r):\A\To\prod_k\A_k.
\]
\begin{enumerate}
\item If the ideals $I_1,\ldots, I_r$ are mutually coprime, then $\prod_k I_k=\bigcap_k I_k$;
\item The homomorphism $\phi$ is surjective iff the $I_k$'s are mutually coprime;
\item In any case, $\Ker\phi=\bigcap_k I_k$.
\end{enumerate}
\end{prop}

Let $\phi:\A\to\B$ is a homomorphism of superalgebras and $I\subset\A$ is a (homogeneous) ideal, we can form its \emph{direct image} $\phi_*I\subset\B$ as the ideal $(\phi(I))$ generated by the image of $I$ under $\phi$. It consists of finite sums of homogeneous elements of the form $b\phi(a)$ with $b\in\B$, $a\in I$ (in other words,
\[
\phi_*I=\B\otimes_\A I
\]
as a $\B$-module).

\begin{prop}\label{prop:dirimage}
\begin{enumerate}
\item $\phi_*$ preserves arbitrary sums and finite products of ideals;
\item if $I,J\subset\A$ are coprime, so are $\phi_*I$ and $\phi_*J$.
\end{enumerate}
\end{prop}
\begin{proof}
Left to the reader.
\end{proof}

\medskip

Now let $\bS$ be a super Fermat theory with ground ring $\KK$. Recall that $\bS$-congruences on $\bS$-algebras are the same thing as homogeneous ideals in the underlying superalgebras. Let $\A\in\bS\Alg$ and let $P=(P_\even|P_\odd)$, where $P_\even\subset\A_\even$, $P_\odd\subset\A_\odd$ are subsets; let $(P)$ denote the homogeneous ideal generated by $P$. Observe that the quotient $\A/(P)$ is the coequalizer of the pair of maps
\[
\bS(P)\rightrightarrows\A,
\]
where the top map sends each generator $x_p$ to the corresponding element $p\in\A$, while the bottom one sends each $x_p$ to $0$. The following is then immediate:

\begin{prop}\label{prop:completionofpresentations}
Let $F:\bS\to\bS'$ be a map of super Fermat theories,
\[
\Adj{F^*}{\bS\Alg}{\bS'\Alg}{F_!}
\]
the corresponding adjunction, $\A\in\bS\Alg$, $I\subset\A$ a homogeneous ideal. Then
\[
F_!(\A/I)=F_!\A/u_*I,
\]
where
\[
u:\A\To F^*F_!\A
\]
is the unit of the adjunction.
\end{prop}

\subsubsection{Completions, coproducts and change of base.} Applying the above proposition to the special case of the structure map $\bscom_\KK\to\bS$ we get:

\begin{cor}
Let $\A\in\bscom_\KK\Alg$, $I\subset\A$ a homogeneous ideal, $\hat{\A}\in\bS\Alg$ the $\bS$-algebra completion of $\A$. Then
\[
\widehat{\A/I}=\hat{\A}/\hat{I},
\]
where $\hat{I}=u_*I$ for $u:\A\to(\hat{\A})_\sharp$ the unit of the adjunction.
\end{cor}

Let now $\B\in\bS\Alg$. The map $u:\KK\to\B$ induces a map of theories $\bS\to\bS_\B$; the corresponding adjunction takes the form
\[
\Adjlong{(\quad)\circ u}{\bS\Alg}{\bS_\B\Alg}{\B\odot_\KK(\quad)}=\B/\bS\Alg
\]
where the right adjoint is simply precomposition with $u$, i.e. it is the functor assigning the underlying $\bS$-algebra, while the left adjoint is the change of base. Notice that the unit of the adjunction is the canonical inclusion into the coproduct:
\[
\iota:\A\To\B\odot\A.
\]

\begin{cor} If $\A\in\bS\Alg$, $I\subset\A$ a homogeneous ideal, then
\[
\B\odot(\A/I)=(\B\odot\A)/\iota_*I
\]
\end{cor}

\begin{cor}
Let $\A_i\in\bS\Alg$, $I_i\subset\A_i$ homogeneous ideals, $i=1,2$. Then
\[
(\A_1/I_1)\odot(\A_2/I_2)=(\A_1\odot\A_2)/(\iota_{1,*}I_1+\iota_{2,*}I_2),
\]
where $\iota_i:\A_i\to\A_1\odot\A_2$ are the canonical inclusions. In particular, if
\[
A_i=\bS(P_i)/I_i,
\]
are presentations, then
\[
\A_1\odot\A_2=\bS(P_1\amalg P_2)/(\iota_{1,*}I_1+\iota_{2,*}I_2).
\]
\end{cor}

\subsubsection{Products.}

As is the case for all algebraic theories, products of $\bS$-algebras are computed ``pointwise'', i.e. on underlying sets. What is quite remarkable is that, in sharp contrast with general algebraic theories, finite products are preserved by the left adjoints of algebraic morphisms between categories of algebras over (super) Fermat theories.

\begin{thm}\label{thm:morphismprodspres}
Let $F:\bS\to\bS'$ be a morphism of super Fermat theories. Then $F_!:\bS\Alg\to\bS'\Alg$ preserves finite products.
\end{thm}
\begin{proof}
Let $\A_1,\ldots,\A_r$ be $\bS$-algebras, and
\[
\A=\prod_k\A_k
\]
their product. Pick any presentation of $\A$, i.e. a pair of sets $P=(P_\even,P_\odd)$ and a surjective homomorphism
\[
\phi:\bS(P)\To\A.
\]
Composing with the canonical projections, we get surjective homomorphisms
\[
\phi_k:\bS(P)\To\A_k,\quad k=1,\ldots,r,
\]
so that
\[
\phi=(\phi_1,\ldots,\phi_r).
\]
Let $I_k=\Ker\phi_k$ for each $k$, and $I=\Ker\phi$, so that
\[
\A_k=\bS(P)/I_k \quad \mathrm{and}\quad\A=\bS(P)/I.
\]
Since $\phi$ is surjective, the $I_k$'s are mutually coprime by Proposition \ref{prop:coprimeideals}. Therefore,
\[
I=\bigcap_k I_k=\prod_k I_k.
\]
Now apply $F_!$. We have
\[
\F_!\A_k=\bS'(P)/u_*I_k
\]
by Proposition \ref{prop:completionofpresentations}, and the ideals $u_*I_k$ are mutually coprime by Proposition \ref{prop:dirimage}. Therefore, the map
\[
\psi=(F_!\phi_1,\ldots,F_!\phi_r):\bS'(P)\To\prod_kF_!\A_k
\]
is surjective by Proposition \ref{prop:coprimeideals} and its kernel is
\[
\Ker\psi=\bigcap_k u_*I_k=\prod_k u_*I_k=u_*(\prod_k I_k)=u_*I
\]
by Propositions \ref{prop:coprimeideals} and \ref{prop:dirimage}. Therefore,
\[
\prod_kF_!\A_k=\bS'(P)/u_*I=F_!\A=F_!(\prod_k\A_k)
\]
and $\psi=F_!\phi$.
\end{proof}

\begin{cor}\label{cor:complepresprod}

\begin{enumerate}
\item For any super Fermat theory $\bS$ with ground ring $\KK$, the completion functor
\[
\widehat{(\quad)}:\bscom_\KK\Alg\To\bS\Alg
\]
preserves finite products;
\item for any $\bS$-algebra $\B$, the change of base functor
\[
\B\odot(\quad):\bS\Alg\To\bS_\B\Alg
\]
preserves finite products.
\end{enumerate}
\end{cor}

\begin{rem}
One can easily see by repeating the above arguments (or by restriction) that the results of this subsection remain valid for morphisms of Fermat theories (not super), and more generally, for morphisms of algebraic theories (over $\bcom$) between Fermat theories and super Fermat theories.
\end{rem}

\begin{rem}
In general, for a morphism $F:\bT\to\bT'$ of algebraic theories, the left adjoint $F_!$ seldom preserves products. For instance, the free $\bT$-algebra functor $\Set\to\bT\Alg$ almost never does.
\end{rem}

\begin{rem}
Although they preserve finite products, left adjoints of algebraic morphisms of categories of algebras over (super) Fermat theories generally fail to preserve other finite limits. For instance, consider the structure map $\bcom_\RR\to\bcinf$ and the corresponding $\bcinf$-completion functor
\[
\widehat{(\quad)}:\bcom_\RR\Alg\To\bcinf\Alg.
\]
The equalizer of the shift by $1$ map
\[
\phi:\RR[x]\To\RR[x],\quad x\mapsto x+1
\]
and the identity is $\RR$, while the equalizer of
\[
\hat\phi:\RR\{x\}\To\RR\{x\}
\]
and the identity is isomorphic to $\cinf(S^1)$: there are no non-constant periodic polynomials, but lots of periodic smooth functions.
\end{rem}







\subsubsection{Localizations}
We end this section by giving a brief account of localization in the context of super Fermat theories. If $\Sigma \subset \A$ is any subset of a $\bE$-algebra $\A,$ with $\bE$ a (super) Fermat theory, one can form the \emph{localization} of $\A$ with respect to $\Sigma,$ $\A\set{\Sigma^{-1}}.$ There is a canonical $\bE$-algebra map $$l:\A \to \A\set{\Sigma^{-1}}$$ which satisfies the following universal property:

Given any $\bE$-algebra $\B,$ any $\bE$-algebra map $$\varphi:\A \to \B$$ which send every element of $\Sigma$ to a unit extends uniquely to a $\bE$-algebra map $$\A\set{\Sigma^{-1}} \to \B.$$ Even in the case where $\Sigma$ is multiplicatively closed, the localization of $\A$ with respect to $\Sigma$ can not usually be computed by the methods customary to commutative algebra. However, the universal properties of $\A\set{\Sigma^{-1}}$ give rise to a canonical presentation. Let $$\A\set{\Sigma}=\A \odot \KK\set{\Sigma}=\A\set{\left(x_s\right)_{s\in\Sigma}},$$ where $\KK\set{\Sigma}$ is the free $\bE$-algebra on $|\Sigma|$-generators (or the free $\bE$-algebra on $|\Sigma_\even|$-even generators and $|\Sigma_\odd|$-odd generators, in the super case). Then $$\A\set{\Sigma^{-1}}=\A\set{\left(x_s\right)_{s\in\Sigma}}/\left(\left(1-s\cdot x_s\right)\right),$$ where $\left(\left(1-s\cdot x_s\right)\right)$ is the ideal generated by all elements of the form $1-s\cdot x_s,$ for some $s\in \Sigma.$

\begin{rem}
Of course, if $\Sigma_\odd$ is non-empty, $\A\set{\Sigma^{-1}}=\{0\}$.
\end{rem}

For certain Fermat theories (besides those of the form $\bcom_\KK$), e.g. the theory of $\bcinf$-algebras, other descriptions of localizations are possible. For example, if $f \in \bcinf\left(\RR^n\right),$ and $\Sigma=\set{f},$ then $$\bcinf\left(\RR^n\right)\set{f^{-1}}\cong \bcinf\left(U\right),$$ where $$U=f^{-1}\left(\RR/\set{0}\right)$$ (c.f. \cite{msia}).

\section{Near-point Determined Algebras}\label{sec:npd}
In this section, we introduce for a (super) Fermat theory $\bE$ its subcategory of near-point determined algebras. We then go on to prove many of their pleasant properties.
\subsection{Radicals}
In this subsection, let $\bE$ be either a Fermat theory or a super Fermat theory, and let $\bQ$ be a full subcategory of $\bE\Alg$.
\begin{defn}
Given an $\A\in\bE\Alg$ and $\Q\in\bQ$, a \emph{$\Q$-point of $\A$} is a homomorphism $p:\A\to\Q$; a \emph{$\bQ$-point} of $\A$ is a $\Q$-point for some $\Q\in\bQ$.
\end{defn}

\begin{defn}
An ideal $P$ of a $\bE$-algebra $\A$ is said to be a \emph{$\bQ$-ideal} if $P$ is the kernel of some $\bQ$-point $p.$ Denote the set of $\bQ$-ideals by $\Spec_\bQ\left(\A\right).$ Let $I$ be an arbitrary ideal of $\A$. Define the $\bQ$-radical of $I$ to be the ideal
$$\Rad_\bQ\left(I\right)=\mathop{\bigcap_{P \supseteq I}}_{P\in\Spec_\bQ\left(\A\right)}P.$$ We call $\Rad_\bQ\left(0\right)$ the $\bQ$-radical of $\A,$ and will denote it by $\RAD_\bQ\left(\A\right).$
\end{defn}

\begin{rem}
Assume that $\bE=\bcom.$ If $\bQ$ is the subcategory of integral domains, then a $\bQ$-ideal is the same as a prime ideal, so $\Spec_\bQ\left(\A\right)$ is the prime spectrum and the $\bQ$-radical of $\A$ is the same thing as the nilradical. When $\bQ$ is the subcategory of fields, a $\bQ$-ideal is the same thing as a maximal ideal, $\Spec_\bQ\left(\A\right)$ is the maximal spectrum, and the $\bQ$-radical of $\A$ is the same thing as the Jacobson radical. Another example is the $\bW$-radical considered in \cite{sdg}, Section III.9, and is closely related to the concept of near-point determined algebras discussed in Section \ref{sec:subnpd} of this paper.
\end{rem}

\begin{prop}
For any ideal $I$ of a $\bE$-algebra $\A,$ we have
\begin{itemize}
\item[1)] $\Rad_\bQ\left(I\right)=\pi_I^{-1}\left(\RAD_\bQ\left(\A/I\right)\right),$ where $$\pi_I:\A \to \A/I$$ is the canonical projection.
\item[2)] $\Rad_\bQ\left(\Rad_\bQ\left(I\right)\right)=\Rad_\bQ\left(I\right).$
\end{itemize}
\end{prop}
\begin{proof}
Condition $1)$ follows immediately from the lattice theory of ideals. For $2),$ observe that for any $\bQ$-point $p$ of $\A,$ such that $\Ker\left(p\right) \supseteq I,$ $p\left(a\right)=0$ for all $a \in \Rad_\bQ\left(I\right),$ by definition.
\end{proof}

\begin{cor}
For every $\bE$-algebra $\A,$ $\Rad_\bQ\left(\RAD_\bQ\left(\A\right)\right)=\RAD_\bQ\left(\A\right).$
\end{cor}

\begin{prop}\label{prop:sptdet}
For a $\bE$-algebra $\A$, the following conditions are equivalent:

\begin{itemize}
\item[1)] $\RAD_\bQ\left(\A\right)=0.$
\item[2)] There is an embedding $$\A \hookrightarrow \prod\limits_\alpha \Q_\alpha$$ of $\A$ into a product of algebras in $\bQ.$
\item[3)] For any pair of maps $f,g:\B\to\A$ we have
\[
\forall\Q\in\bQ,\;\forall\: p:\A\to\Q,\;p\circ f=p\circ g\implies f=g.
\]
\item[4)] For any element $a \in \A,$ if $p\left(a\right)=0$ for all $\bQ$-points $p,$ then $a=0.$
\end{itemize}
\end{prop}
\begin{proof}
Suppose $1)$ holds. Choose for each $\bQ$-ideal $P$ a homomorphism $\A \to \Q_P$ with $\Q_P \in \bQ$ whose kernel is $P,$ and denote the associated embedding $$\A/P \hookrightarrow \Q_P$$ by $\varphi_P.$ Denote $$\varphi:=\prod\limits_{P\in\Spec_\bQ\left(\A\right)} \varphi_P:\prod\limits_{P\in\Spec_\bQ\left(\A\right)} \A/P \hookrightarrow \prod\limits_{P\in\Spec_\bQ\left(\A\right)} \Q_P,$$ and consider the canonical composite
$$\A \stackrel{\theta}{\longrightarrow} \prod\limits_{P\in\Spec_\bQ\left(\A\right)} \A/P \stackrel{\varphi}{\longrightarrow} \prod\limits_{P\in\Spec_\bQ\left(\A\right)} \Q_P.$$
The kernel of $\theta$ is $\Rad_\bQ\left(\A\right)=0,$ hence the composite is an embedding of $\A$ into a product of algebras in $\Q.$ $2) \implies 3)$ is obvious. Now, suppose that $3)$ holds.
For simplicity, we will assume that $\bE$ is Fermat as opposed to super Fermat, however the proof for the super case is nearly identical. Suppose that $a \in \A$ has $p\left(a\right)=0$ for every $\bQ$-point of $\A.$ Consider the free $\bE$-algebra on one generator, $\KK\{x\}.$ There is a natural bijection $$\Hom\left(\KK\{x\},\A\right) \cong \underline \A,$$ where $\underline \A$ is the underlying set of $\A$. Each object $t$ of $\A$ corresponds to a unique morphism $\lambda^\A_t:\KK\{x\} \to \A$ sending $x$ to $t.$ Now, $p\circ \lambda^\A_a = p\circ \lambda^\A_0,$  for all $p:\A \to \Q,$ with $\Q \in \bQ,$ since both expressions are equal to $\lambda^\Q_0,$ the morphism $$\KK\{x\} \to S$$ classifying the element $0 \in \Q.$ Assuming $3),$ it follows that $\lambda^\A_a=\lambda^\A_0,$ hence $a=0.$ $4) \implies 1)$ is obvious.
\end{proof}

\begin{defn}
If an $\bE$-algebra $\A$ satisfies either of the equivalent conditions of Proposition \ref{prop:sptdet}, it is said to be \emph{$\bQ$-point determined.} Denote the full subcategory of $\bQ$-point determined algebras by $\bE\Alg_{\bQ\mathrm{det}}$. An ideal $I$ is said to be \emph{$\bQ$-point determined} (or \emph{$\bQ$-radical}) if $\Rad_{\bQ}\left(I\right)=I.$
\end{defn}

\begin{rem}\label{rem:subpt}
If $\A$ is $\bQ$-point determined, then any sub-$\bE$-algebra $\B$ of $\A$ is also $\bQ$-point determined.
\end{rem}

\begin{rem}
Assume that $\bE=\bcom.$ If $\bQ$ is the subcategory of integral domains, then a commutative ring is $\bQ$-point determined if and only if it is reduced. When $\bQ$ is the subcategory of fields, a commutative ring is $\bQ$-point determined if and only if it is Jacobson semisimple (a.k.a semiprimitive).
\end{rem}

\begin{prop}\label{prop:pdtdideal}
Let $\A$ be any $\bE$-algebra, and $I$ an ideal. Then $\A/I$ is $\bQ$-point determined if and only if $I$ is.
\end{prop}
\begin{proof}
For any ideal $J$ of $\A$ such that $J\supseteq I,$ $\pi_I\left(J\right)=0$ if and only if $J=I,$ where $$\pi_I:\A \to \A/I$$ is the canonical projection. Hence
$$\RAD_\bQ\left(\A/I\right)=\pi_I\left(\Rad_\bQ\left(I\right)\right)=0$$ if and only if $\Rad_\bQ\left(I\right)=I.$
\end{proof}

\begin{cor}
For any $\A,$ $\A/\RAD_\bQ\left(\A\right)$ is $\bQ$-point determined.
\end{cor}

\begin{prop}
Let $I$ be an ideal of $\A.$ Then $I$ is $\bQ$-point determined if and only if $I$ is the kernel of a homomorphism $f:\A \to \B,$ with $\B$ a $\bQ$-point determined algebra.
\end{prop}
\begin{proof}
Suppose  that $I$ is $\bQ$-point determined. Then by Proposition \ref{prop:pdtdideal}, $\A/I$ is $\bQ$-point determined, and $I$ is the kernel of $$\A \to \A/I.$$ Conversely, suppose that $f:\A \to \B$ and $\B$ is $\bQ$-point determined. Then $\A/\Ker\left(f\right)$ is a sub-$\bE$-algebra of $\B,$ so is $\bQ$-point determined. So by Proposition \ref{prop:pdtdideal}, $\Ker\left(f\right)$ is $\bQ$-point determined.
\end{proof}

\begin{prop}\label{prop:npdreflects}
The assignment $\A \mapsto \A/\RAD_\bQ\left(\A\right)$ extends to a functor
$$L_\bQ:\bE\Alg \to \bE\Alg_{\bQ\mathrm{det}}$$ which is left adjoint to the inclusion, with the unit $\eta_\A:\A \to \A/\RAD_\bQ\left(\A\right)$ given by the canonical projection.
\end{prop}
\begin{proof}
If $\phi:\A\to\B$ is a map and $a\in \RAD_\bQ\left(\A\right)$, then for every $\Q\in\bQ$ and every $g:\B\to\Q$ we have $g(\phi(a))=(g\circ\phi)(a)=0$, so $\phi(a)\in \RAD_\bQ\left(\B\right)$, and thus $L_\bQ$ is indeed a functor.  It is left adjoint to the inclusion with the described unit since any map from an arbitrary $\A$ to a $\bQ$-point determined $\B$ must send the elements of $\RAD_\bQ\left(\A\right)$ to $0$, hence factors uniquely through $\eta_\A.$
\end{proof}

\begin{rem}
Given a subcategory $\bQ$ of $\bE\Alg,$ we may consider its \emph{saturation} $\overline{\bQ}$ with respect to arbitrary products and subobjects, i.e. the smallest subcategory of $\bE\Alg$ which is closed under arbitrary products and subobjects, which contains $\bQ$. On one hand, by Proposition \ref{prop:sptdet}, the category of $\bQ$-point determined algebras is contained in $\overline{\bQ}.$ On the other hand, by Proposition \ref{prop:npdreflects} the full subcategory $\bE\Alg_{\bQ\mathrm{det}}$ is reflective, hence closed under arbitrary limits. In particular, it is closed under arbitrary products. Moreover, by Remark \ref{rem:subpt}, $\bE\Alg_{\bQ\mathrm{det}}$ is closed under subobjects. Hence, $\overline{\bQ}$ is contained in $\bE\Alg_{\bQ\mathrm{det}}.$ Therefore the subcategory $\bE\Alg_{\bQ\mathrm{det}}$ of $\bE\Alg$ may be identified with the saturation $\overline{\bQ}.$ Notice that this notion of saturation makes sense in a much more general context, even where many of the various equivalent conditions in Proposition \ref{prop:sptdet} do not make sense.
\end{rem}

Recall that the radical of an ideal $I$ in a commutative ring $\A$ is given by $$\sqrt{I}=\{a \in \A\mspace{2mu}|\mspace{3mu}a^n \in I\mbox{ for some }n\in\ZZ_+\}.$$ This is the same as $\Rad_{\bQ}\left(I\right),$ when $\bQ$ is the subcategory of $\bcom\Alg$ consisting of integral domains. An important property of the radical is that for any two ideals $I$ and $J$ of $\A,$ \begin{equation}\label{eq:sqrt}
\sqrt{I \cap J}=\sqrt{I} \cap \sqrt{J}.
\end{equation}
An analogous equation holds for the Jacobson radical of ideals. A natural question is, for $\bQ$ any subcategory of $\bE\Alg,$ when does (\ref{eq:sqrt}) hold? The following proposition offers a partial answer:

\begin{prop}
If each $\Q$ in $\bQ$ is an integral domain, then the following equation is satisfied
\begin{equation}\label{eq:satt}
\Rad_\bQ\left(I\cap J\right)=\Rad_\bQ\left(I\right)\cap \Rad_\bQ\left(J\right)
\end{equation}
for all $\A \in \bE\Alg$ and all $I$ and $J$ ideals of $\A.$
\end{prop}
\begin{proof}
Suppose that each $\Q$ in $\bQ$ is an integral domain. Notice that the inclusion $$\Rad_\bQ\left(I\cap J\right) \subseteq \Rad_\bQ\left(I\right)\cap \Rad_\bQ\left(J\right)$$ is always true. It suffices to show the reverse inclusion. Let $$p:\A \to \Q$$ be a $\bQ$-point of $\A,$ such that $p\left(I\cap J\right)=0.$ Notice that $IJ \subseteq I \cap J,$ so
\begin{equation}\label{eq:integ}
p\left(ij\right)=p\left(i\right)p\left(j\right)=0
\end{equation}
for all $i \in I$ and $j \in J.$ Suppose that $p\left(I\right) \ne 0.$  Then there exists $i \in I$ such that $p\left(i\right)\ne 0.$ In this case, equation (\ref{eq:integ}) holds in $\Q,$ which is an integral domain. It follows that $p\left(j\right)=0,$ for all $j\in J,$ i.e. $p\left(J\right)=0.$ So, for every $\bQ$-point $p$ whose kernel contains $I\cap J$, either $\Ker\left(p\right)$ contains $I$ or it contains $J.$  It follows that $$\Rad_\bQ\left(I\right)\cap \Rad_\bQ\left(J\right) \subseteq \Rad_\bQ\left(I\cap J\right).$$
\end{proof}

This last proposition explains why (\ref{eq:satt}) is satisfied in the case of prime and Jacobson radicals. If $\bQ$ does not consist entirely of integral domains, equation (\ref{eq:satt}) may still be satisfied for \emph{coprime} ideals, as the following proposition shows:

\begin{prop}\label{prop:coprimeradical}
If every algebra $\Q$ in $\bQ$ is a local $\bE$-algebra, then for all $\bE$-algebras $\A$ and all coprime ideals $I$ and $J$ of $\A,$ equation (\ref{eq:satt}) holds.
\end{prop}
\begin{proof}
It suffices to show that if $p:\A \to \Q$ is a $\bQ$-point of $\A$ such that $$p\left(I \cap J\right)=0,$$ then either $p\left(I\right)=0$ or $p\left(J\right)=0.$ Since $I$ and $J$ are coprime, there exists $\zeta \in I$ and $\omega \in J$ such that $$\zeta + \omega=1.$$ Hence $$p\left(\zeta\right)+ p\left(\omega\right)=1.$$ Since $\Q$ is local, either $p\left(\zeta\right)$ or $p\left(\omega\right)$ is a unit, otherwise they would both be in the unique maximal ideal $\mathfrak{m}$ of $\Q,$ but this would imply that $1 \in \mathfrak{m},$ which is absurd. Assume without loss of generality that $p\left(\zeta\right)$ is a unit. Then, for all $j \in J,$ $$p\left(\zeta j\right)=p\left(\zeta\right)p\left(j\right)=0,$$ and since $p\left(\zeta\right)$ is a unit, this implies $p\left(j\right)=0,$ for all $j\in J.$
\end{proof}

\begin{cor}\label{cor:reflectoprods}
If every algebra $\Q$ in $\bQ$ is a local $\bE$-algebra, then the reflector $$L_\bQ:\bE\Alg \to \bE\Alg_{\bQ\mathrm{det}}$$ preserves finite products.
\end{cor}

\begin{proof}
The reflector $L_\bQ$ always preserves the terminal algebra. Let $\A$ and $\B$ be $\bE$-algebras. Consider the composite of surjections $$\A \times \B \stackrel{pr_1}{\longlongrightarrow} \A \stackrel{\eta_\A}{\longlongrightarrow} \A/\RAD_\bQ\left(\A\right),$$ and similarly with the role of $\A$ and $\B$ exchanged. Their kernels are $\Rad_\bQ\left(\{0\} \times \B\right)$ and $\Rad_\bQ\left( \A\times \{0\}\right)$ respectively. Notice that these two ideals are coprime since the former contains $\left(0,1\right)$ and the latter contains $\left(1,0\right).$  Hence, by Proposition \ref{prop:coprimeideals}, it follows that the induced map $$\A \times \B \to \A/\RAD_\bQ\left(\A\right) \times \B/\RAD_\bQ\left(\B\right)$$ is surjective, with kernel $\Rad_\bQ\left(\{0\} \times \B\right) \cap \Rad_\bQ\left( \A\times \{0\}\right).$ Since every algebra $\Q$ in $\bQ$ is a local $\bE$-algebra, by Proposition \ref{prop:coprimeradical}, this kernel is equal to $$\Rad_\bQ\left(\left(\{0\} \times \B\right) \cap \left( \A\times \{0\}\right)\right)=\Rad_\bQ\left(0\right)=\RAD_\bQ\left(\A \times \B\right).$$ By the first isomorphism theorem, it follows that
\begin{eqnarray*}
L_\bQ\left(\A \times \B\right)&=&\left(\A\times\B\right)/\RAD_\bQ\left(\A\times\B\right)\\
&\cong& \A/\RAD_\bQ\left(\A\right) \times \B/\RAD_\bQ\left(\B\right)\\
&=&L_\bQ\left(\A\right) \times L_\bQ\left(\B\right).
\end{eqnarray*}
\end{proof}

\begin{rem}
By the same proof, $L_\bQ$ also preserves finite products if every algebra $\Q$ in $\bQ$ is an integral domain.
\end{rem}

\subsubsection{Relative Reduction}
Let $\bE$ be a super Fermat theory, and let $\bQ$ be a full subcategory of $\bE\Alg$. Denote by
$$j_\bQ:\bQ \hookrightarrow \bE\Alg$$ the full and faithful inclusion. Consider the forgetful functor
$$U_\bE:\bE\Alg \to \Set^{\set{\even,\odd}}.$$ Denote the composite by $$\underline{\KK_\bQ}:= U_\bE \circ j_\bQ.$$ The $\ZZ_2$-graded object $\set{\underline{\KK_\bQ}_\even,\underline{\KK_\bQ}_\odd}$ of $\Set^\bQ,$ may be regarded as a $\bE$-algebra in the topos $\Set^\bQ.$ For each $\Q \in \bQ,$ there is universal map of $2$-sorted Lawvere theories $$\chi_\Q:\bE \to \End_{U_\bE\left(\Q\right)}$$ classifying $\Q.$ (See Example \ref{eg:end}.) Hence, for all $\left(n|m\right)$ and $\left(p|q\right),$ the functor $\chi_\Q$ provides natural maps $$\bE\left(\left(n|m\right),\left(p|q\right)\right) \to \Set\left(\Q^n_\even \times \Q^m_\odd,\Q^p_\even \times \Q^q_\odd\right).$$ Pick $f \in \bE\left(\left(n|m\right),\left(p|q\right)\right),$ then the maps $\left(\chi_\Q\left(f\right)\right)_{\Q \in \bQ}$ assemble into a natural transformation (i.e. a map in $\Set^{\bQ}$) $$\ev\left(f\right):\underline{\KK_\bQ}^n_\even \times \underline{\KK_\bQ}^m_\odd \to \underline{\KK_\bQ}^p_\even \times \underline{\KK_\bQ}^q_\odd.$$ This yields functions
\begin{equation}\label{eq:evals}
\ev_{\left(n|m\right),\left(p|q\right)}:\bE\left(\left(n|m\right),\left(p|q\right)\right) \to \Set^{\bQ}\left(\underline{\KK_\bQ}^n_\even \times \underline{\KK_\bQ}^m_\odd,\underline{\KK_\bQ}^p_\even \times \underline{\KK_\bQ}^q_\odd\right).
\end{equation}
Note that $$\bE\left(\left(n|m\right),\left(1|0\right)\right)\cong \KK\set{x^1,\cdots,x^n,\xi^1,\cdots,\xi^m}_\even$$ and $$\bE\left(\left(n|m\right),\left(0|1\right)\right)\cong \KK\set{x^1,\cdots,x^n,\xi^1,\cdots,\xi^m}_\odd.$$ Hence, we get even and odd evaluation maps:
$$\ev^{\left(n|m\right)}_\even:\KK\set{x^1,\cdots,x^n,\xi^1,\cdots,\xi^m}_\even \to \Set^{\bQ}\left(\underline{\KK_\bQ}^n_\even \times \underline{\KK_\bQ}^m_\odd,\underline{\KK_\bQ}_\even\right)$$ and $$\ev^{\left(n|m\right)}_\odd:\KK\set{x^1,\cdots,x^n,\xi^1,\cdots,\xi^m}_\odd \to \Set^{\bQ}\left(\underline{\KK_\bQ}^n_\even \times \underline{\KK_\bQ}^m_\odd,\underline{\KK_\bQ}_\odd\right).$$

\begin{rem}
If $\bE$ is a non-super Fermat theory, all of this construction carries through; however, one only needs to use one sort.
\end{rem}

\begin{defn}
The super Fermat theory $\bE$ is \emph{$\bQ$-reduced} if for all $\left(n|m\right)$ the evaluation maps $\ev^{\left(n|m\right)}_\even$ and $\ev^{\left(n|m\right)}_\odd$ are injective.
\end{defn}

\begin{prop}\label{prop:whatreducedmeans}
A super Fermat theory $\bE$ is $\bQ$-reduced if and only if each finitely generated free $\bE$-algebra is $\bQ$-point determined.
\end{prop}
\begin{proof}
Suppose that for some $n$ and $m$, $\KK\set{x^1,\cdots,x^n,\xi^1,\cdots,\xi^m}$ is not $\bQ$-point determined. Then there is a non-zero $f \in \KK\set{x^1,\cdots,x^n,\xi^1,\cdots,\xi^m}$ such that $\varphi\left(f\right)=0$ for all
$$\varphi:\KK\set{x^1,\cdots,x^n,\xi^1,\cdots,\xi^m} \to \Q,$$ with $\Q \in \bQ.$ So, for all $\bQ$-points $\varphi,$ we have
\begin{equation}\label{eq:redq}
\varphi\left(f\right)=\Q\left(f\right)\left(\left(\varphi\left(x^1\right),\cdots,\varphi\left(x^n\right)\right),\left(\varphi\left(\xi^1\right),\cdots,\varphi\left(\xi^m\right)\right)\right)=0.
\end{equation}
Since $\KK\set{x^1,\cdots,x^n,\xi^1,\cdots,\xi^m}$ is free, this implies for all $\Q$ and any collection $$a_1,\cdots,a_n$$ of even elements of $\Q$ and $$b_1,\cdots,b_m$$ odd elements, $$\Q\left(f\right)\left(\left(a_1,\cdots,a_n\right),\left(b_1,\cdots,b_m\right)\right)=0.$$ Hence, for all $\Q,$ $$\chi_\Q\left(f\right)=\chi_\Q\left(0\right).$$ In particular, this implies that the evaluation map $\ev^{\left(n|m\right)}$ of the same parity as $f$ is not injective.

Conversely, suppose that each finitely generated free $\bE$-algebra is $\bQ$-point determined. Suppose that $f$ and $g$ are in $\KK\set{x^1,\cdots,x^n,\xi^1,\cdots,\xi^m},$ have the same parity, and $$\ev^{\left(n|m\right)}\left(f\right)=\ev^{\left(n|m\right)}\left(g\right).$$ By equation (\ref{eq:redq}), this implies that for all $\bQ$-points $\varphi,$ $\varphi\left(f-g\right)=0.$ By Proposition \ref{prop:sptdet}, this implies that $f=g,$ so that each evaluation map is injective.
\end{proof}

\begin{cor}\label{cor:rednpdis}
A super Fermat theory $\bE$ is $\bQ$-reduced if and only if each free $\bE$-algebra is $\bQ$-point determined.
\end{cor}

\begin{proof}
If every free $\bE$-algebra is $\bQ$-point determined, then $\bE$ is $\bQ$-reduced by Proposition \ref{prop:whatreducedmeans}. Suppose that $\bE$ is $\bQ$-reduced. By the same proposition, every finitely generated $\bE$-algebra is $\bQ$-point determined. Let $\bT$ be some $\ZZ_2$-graded set and let $\bE\left(\bT\right)$ be the free $\bE$-algebra on $\bT$. Suppose that $f$ and $g$ are two elements thereof and that for every $\bQ$-point $$p:\bE\left(\bT\right) \to \Q,$$ $p\left(f\right)=p\left(g\right).$ Since $\bE\left(\bT\right)$ is a filtered colimit of finitely generated free algebras, there exists a finite $\ZZ_2$-graded subset $\bT_0$ of $\bT$ such that $f$ and $g$ are in the image of $$i:\bE\left(\bT_0\right) \to \bE\left(\bT\right).$$ Say $i\left(\tilde f\right)=f$ and $i\left(\tilde g\right)=g.$ Let $$q:\bE\left(\bT_0\right) \to \Q$$ be any $\bQ$-point. Then $q$ can be extended along $i$ to a $\bQ$-point $p$, for example, by setting $$p\left(t\right)=q\left(t\right)$$ for all $t \in \bT_0,$ and by letting $p$ be zero on all other generators. This implies that $$q\left(\tilde f\right)=pi\left(\tilde f\right)=p\left(f\right).$$ Hence $q\left(\tilde f\right)=q\left(\tilde g\right)$ for every $\bQ$-point $q,$ and hence $\tilde f = \tilde g,$ since $\bE\left(\bT_0\right)$ is finitely generated, and hence $\bQ$-point determined. Therefore, $f=g,$ and $\bE\left(\bT\right)$ is also $\bQ$-point determined.
\end{proof}

\begin{rem}
A non-super Fermat theory $\bE$ is reduced if and only if it is $\bK$-reduced, where $\bK$ is the full subcategory of $\bE$ spanned by the initial $\bE$-algebra $\KK.$
\end{rem}

\begin{defn}
Let $\bE$ be a super Fermat theory, and let $\mathbf{\Lambda}$ denote the subcategory of $\bE$ consisting of the Grassman algebras (Definition \ref{dfn:grass}.) The super Fermat theory $\bE$ is \emph{super reduced} if it is $\mathbf{\Lambda}$-reduced.
\end{defn}

\begin{prop}
If $\bE$ is a reduced Fermat theory, $\bS\bE$ is a super reduced super Fermat theory.
\end{prop}

\begin{proof}
By Proposition \ref{prop:superization}, $\bS\bE$ is a super Fermat theory. It suffices to show that it is super reduced. By Proposition \ref{prop:whatreducedmeans}, it suffices to show that for all $n$ and $m,$ $E\left(n\right) \otimes \Lambda^m$ is $\mathbf{\Lambda}$-point determined. Since $\bE$ is reduced, by Proposition \ref{prop:whatreducedmeans}, each $E\left(n\right)$ is $\bK$-point determined, hence by Proposition \ref{prop:sptdet}, there is an embedding $$\psi:E\left(n\right) \hookrightarrow \prod\limits_\alpha \KK.$$ Consider now the composite $$E\left(n\right) \otimes \Lambda^m \to \left(\prod\limits_\alpha \KK\right) \otimes \Lambda^m \to \prod\limits_\alpha \Lambda^m.$$ The first morphism is injective since $\Lambda^m$ is free, hence flat as a $\KK$-module. The second is always injective. Hence, we have an embedding of $E\left(n\right) \otimes \Lambda^m$ into a copy of algebras in $\mathbf{\Lambda},$ so by Proposition \ref{prop:sptdet}, we are done.
\end{proof}

Define a $2$-sorted Lawvere theory $\End_{\underline{\KK_\bQ}}$ by setting $$\End_{\underline{\KK_\bQ}}\left(\left(n|m\right),\left(p|q\right)\right)=\Set^{\bQ}\left(\underline{\KK_\bQ}^n_\even \times \underline{\KK_\bQ}^m_\odd,\underline{\KK_\bQ}^p_\even \times \underline{\KK_\bQ}^q_\odd\right).$$
Notice that (\ref{eq:evals}) yields a canonical map of theories $\ev_\bQ:\bE \to \End_{\underline{\KK_\bQ}}.$ It is easy to see that $\bE$ is $\bQ$-reduced if and only if this map is faithful. Moreover, the $\ZZ_2$-graded set $$\Set^{\bQ}\left(\underline{\KK_\bQ}^n_\even \times \underline{\KK_\bQ}^m_\odd,\underline{\KK_\bQ}\right):=\set{\Set^{\bQ}\left(\underline{\KK_\bQ}^n_\even \times \underline{\KK_\bQ}^m_\odd,\underline{\KK_\bQ}_\even\right),\Set^{\bQ}\left(\underline{\KK_\bQ}^n_\even \times \underline{\KK_\bQ}^m_\odd,\underline{\KK_\bQ}_\odd\right)}$$ has the point-wise structure of an $\bE$-algebra, and the morphisms
$\ev^{\left(n|m\right)}_\even$ and $\ev^{\left(n|m\right)}_\odd$ define an $\bE$-algebra map $$\ev^{\left(n|m\right)}_\bQ:\KK\set{x^1,\cdots,x^n,\xi^1,\cdots,\xi^m} \to \Set^{\bQ}\left(\underline{\KK_\bQ}^n_\even \times \underline{\KK_\bQ}^m_\odd,\underline{\KK_\bQ}\right).$$
\begin{defn}
Given a super Fermat theory $\bE,$ we define its $\bQ$-reduction $\bE_{\bQ\red}$ to be the image of $\ev_\bQ.$ Explicitly, the finitely generated $\bE_{\bQ\red}$-algebra on the sort $\left(n|m\right)$ is given by $$\bE_{\bQ\red}\left(n|m\right)=\Im\left(\ev^{\left(n|m\right)}_\bQ\right)=\bE\left(n|m\right)/\Ker\left(\ev^{\left(n|m\right)}_\bQ\right).$$
\end{defn}

\begin{rem}
By the proof of Proposition \ref{prop:whatreducedmeans} one sees that $$\Ker\left(\ev^{\left(n|m\right)}_\bQ\right)=\RAD_\bQ\left(\bE\left(n|m\right)\right),$$ so that $\bE_{\bQ\red}\left(n|m\right)=L_\bQ\left(\bE\left(n|m\right)\right).$
\end{rem}

The proof of Proposition \ref{prop:reduction} readily generalizes:

\begin{prop}
If $\bE$ is a super Fermat theory, $\bE_{\bQ\red}$ is a $\bQ$-reduced super Fermat theory. Moreover, the assignment $\bE\mapsto\bE_{\bQ\red}$ is functorial and is left adjoint to the inclusion
\[
\mathbf{SFTh}_{\bQ\red}\hookrightarrow\mathbf{SFTh}
\]
of the full subcategory of $\bQ$-reduced super Fermat theories. In particular, super reduced super Fermat theories are a reflective subcategory of super Fermat theories.
\end{prop}

\subsection{Near-point determined superalgebras}

\subsubsection{Near-point determined superalgebras.}\label{sec:subnpd}

\begin{defn}\label{dfn:npd}
Let $\bN$ denote the class of $\fweil$ $\KK$-algebras. An $\bE$-algebra which is $\bN$-point determined is said to be \emph{near-point determined.} Denote the associated subcategory by $\Np.$
\end{defn}

\begin{rem}
If one replaces the role of nilpotent extensions with that of locally nilpotent extensions, (as in Remark \ref{rem:locallynil}), one arrives at an equivalent definition of near-point determined. The reason is that if $$\A \to \KK$$ is a locally nilpotent extension with kernel $N,$ the natural map $$\A \to \prod\limits_{n=0}^{\infty} \A/N^{n+1}$$ is an embedding into a product of formal Weil algebras, hence $\A$ is near-point determined.
\end{rem}

\begin{rem}\label{rem:reflsprods}
In light of Remark \ref{rem:local}, from Corollary \ref{cor:reflectoprods} it follows that, if $\KK$ is a field, each $\fweilk$ is local, so the reflector  $$L_\bN:\bE\Alg \to \Np$$ preserves finite products.
\end{rem}

\begin{prop}\label{prop:redisnpdfree}
If $\bE$ is a super reduced super Fermat theory, then each free $\bE$-algebra is near-point determined.
\end{prop}
\begin{proof}
Since $\Lambda \subset \bN,$ the result follows from Corollary \ref{cor:rednpdis}.
\end{proof}

We will now give an alternate characterization of what it means to be near-point determined. First, we will make some basic observations. Suppose that $$p:\A \to \KK$$ is a $\KK$-point of an $\bE$-algebra $\A$, where $\KK$ is the ground ring. Let $M$ denote the kernel of $p.$ For any $k \ge 0,$ there is a canonical factorization
$$\xymatrix@R=0.5cm{\A \ar[dr]_-{\pi} \ar[rr]^-{p} & & \KK\\
& \A/M^{k+1}_p. \ar[ru]_-{\tilde p} & }$$
If $\pi\left(a_1\right),\ldots\pi\left(a_{k+1}\right)$ are arbitrary elements of $\Ker\left(\tilde p\right),$ then each $a_i \in M_p,$ so $$\pi\left(a_1\right)\pi\left(a_2\right)\cdots\pi\left(a_{k+1}\right)=\pi\left(a_1a_2\cdots a_{k+1}\right)=0.$$ So, $\Ker\left(\tilde p\right)$ is nilpotent of degree $k$ and therefore $\A/M^{k+1}_p$ is a $\fweilk.$ We introduce the notation $$\A_p^{(k)}:=\A/M^{k+1}_p.$$
We note that $\A_p^{(k)}$ is universal among $\fweilk\mbox{s}$ of nilpotency degree $k$ covering the $\KK$-point $p.$ I.e., if $$\rho:\W \to \KK$$ is a $\fweilk$  with $\Ker\left(\rho\right)^{k+1}=0,$ and $$\varphi:\A \to \W$$ is such that $$\rho \circ \varphi = p,$$ then there is a unique factorization $$\xymatrix{\A \ar[d]_-{\pi} \ar[r]^-{\varphi} \ar[r] & \W \ar[d]^-{\rho} \\ \A_p^{(k)} \ar@{-->}[ru]^{\tilde \varphi} \ar[r]^-{\tilde p} & \KK.}$$


\begin{prop}\label{prop:altnpd}
An $\bE$-algebra $\A$ is near-point determined if and only if the canonical map
\begin{equation}\label{eq:npdinj}
\A \to \prod\limits_{p:\A \to \KK} \prod\limits_{k=0}^{\infty} \A_p^{(k)}
\end{equation}
is injective.
\end{prop}
\begin{proof}
Since each $\A_p^{(k)}$ is a $\fweilk,$ if (\ref{eq:npdinj}) is injective, it is an embedding of $\A$ into a product of $\fweil$ $\KK$-algebras, so $\A$ is near-point determined. Conversely, suppose that $\A$ is near-point determined. Notice that the kernel of (\ref{eq:npdinj}) is the intersection over all $\KK$-points $p$ of $\A$ and all $k \ge 1,$ of $\Ker(p)^k.$ Let $a$ be a non-zero element of $\A.$ Then there exists a morphism $\phi:\A \to \W$ to a $\fweilk$ such that $\phi\left(a\right) \ne 0.$ The algebra $\W$ comes equipped with a surjection $$\rho:\W \to \KK.$$ Let $p:=\rho \circ \phi$ and let $M_p$ denote the kernel of $p.$ Notice that $$\phi\left(M_p\right) \subset \Ker\left(\rho\right).$$ Let $n$ be the nilpotency degree of $\Ker\left(\rho\right).$ Then there is a unique factorization of $\phi$ of the form $$\A \stackrel{\lambda}{\longrightarrow} A/M^{n+1}_p=\A^{(n)}_p \stackrel{\tilde \phi}{\longrightarrow} \W.$$ Since $\phi\left(a\right) \ne 0,$ $\lambda\left(a\right) \ne 0,$ so $a$ is not in $M^{n+1}_p.$ Hence, $a$ is not in the kernel of (\ref{eq:npdinj}). It follows  that (\ref{eq:npdinj}) is injective.
\end{proof}

\begin{rem}
It follows that a $\cinf$-algebra is near-point determined in the sense of Definition \ref{dfn:npd} if and only if it is near-point determined in the sense of \cite{Borisov}.
\end{rem}

\begin{lem}\label{lem:algmapsok}
Suppose that $F:\bT' \to \bT$ is a morphism of $\SS$-sorted Lawvere theories. Let $\bD$ be the full subcategory of $\bT$-algebras on those algebras $\A$ with the property that for any $\bT$-algebra $\B,$ any $\bT'$-algebra morphism $$f:F^*\left(\B\right) \to F^*\left(\A\right)$$ is of the form $F^*\left(g\right)$ for a unique $\bT$-algebra morphism $$g:\B \to \A.$$ Then $\bD$ is closed under subobjects and arbitrary limits in $\bT\Alg.$
\end{lem}

\begin{proof}
The fact that $\bD$ is closed under arbitrary limits is clear by universal properties, since $F^*$ is a right adjoint. Suppose that $$j:\C \hookrightarrow \A$$ is a sub-$\bT$-algebra of $\A,$ with $\A \in \bD.$ We wish to show that $\C$ is in $\bD.$ Let $$\varphi:F^*\left(\B\right) \to F^*\left(\A\right)$$ be a $\bT'$-algebra morphism. We wish to show that for all $$f \in \bT\left(\left(n_s\right),\left(m_s\right)\right),$$ the following diagram commutes
\begin{equation}\label{eq:injjjj}
\xymatrix{\prod\limits_s \B^{n_s}_s \ar[d]_-{\prod\limits_s \varphi^{n_s}_s} \ar[r]^-{\B\left(f\right)} & \prod\limits_s \B^{m_s}_s \ar[d]^-{\prod\limits_s \varphi^{m_s}_s}\\
\prod\limits_s \C^{n_s}_s \ar[r]^-{\C\left(f\right)} & \prod\limits_s \C^{m_s}_s.}
\end{equation}
Notice, however, that the following two diagrams commute since $j$ and $F^*\left(j\right) \circ f$ are $\bT$-algebra maps:
$$\xymatrix{\prod\limits_s \C^{n_s}_s \ar[d]_-{\prod\limits_s j^{n_s}_s} \ar[r]^-{\C\left(f\right)} & \prod\limits_s \C^{m_s}_s \ar[d]^-{\prod\limits_s j^{m_s}_s}\\
\prod\limits_s \A^{n_s}_s \ar[r]^-{\A\left(f\right)} & \prod\limits_s \A^{m_s}_s}$$
$$\xymatrix{\prod\limits_s \B^{n_s}_s \ar[d]_-{\prod\limits_s \varphi^{n_s}_s} \ar[r]^-{\B\left(f\right)} & \prod\limits_s \B^{m_s}_s \ar[d]^-{\prod\limits_s \varphi^{m_s}_s}\\
\prod\limits_s \C^{n_s}_s \ar[d]_-{\prod\limits_s j^{n_s}_s}  & \prod\limits_s \C^{m_s}_s \ar[d]^-{\prod\limits_s j^{m_s}_s}\\
\prod\limits_s \A^{n_s}_s \ar[r]^-{\A\left(f\right)} & \prod\limits_s \A^{m_s}_s.}$$
Since $j$ is a monomorphism, this implies that diagram (\ref{eq:injjjj}) commutes, so we are done.
\end{proof}

\begin{cor}\label{cor:npdralgs}
If $\A$ and $\B$ are $\bE$-algebras and $\B$ is near-point determined, then any $\KK$-algebra morphism $\varphi:\A_\sharp \to \B_\sharp$ is a map of $\bE$-algebras.
\end{cor}

\begin{proof}
This is true when $\B$ is a formal Weil algebra by Corollary \ref{cor:weilisgoodsup}. The result now follows from Lemma \ref{lem:algmapsok}, since by definition, any near-point determined $\bE$-algebra is a subalgebra of a product of formal Weil algebras.

\end{proof}

\begin{rem}
In case that $\bE=\bcinf,$ Corollary \ref{cor:npdralgs} gives a completely algebraic proof of \cite{Borisov}, Proposition $8$. (The proof in \cite{Borisov} uses a topological methods.)
\end{rem}

\begin{rem}
The near point determined condition in Corollary \ref{cor:npdralgs} is necessary. It was shown in \cite{Reichard}, that there is a counterexample in the case of $\bcinf$-algebras. In slightly more detail, by Borel's theorem (c.f. \cite{msia}), the canonical $\RR$-algebra map $$T:\bcinf\left(\RR\right)_0 \to \RR[[x]],$$ from the algebra of germs of smooth functions at the origin, to the algebra of formal power series, assigning the germ of a function $f$ its Taylor polynomial, is surjective. By Corollary \ref{cor:idealsarecongs}, this endows $\RR[[x]]$ with the canonical structure of a $\bcinf$-algebra, making $T$ a $\bcinf$-map. Reichard proves (assuming the axiom of choice) in \cite{Reichard} that there exists an $\RR$-algebra map $$\phi:\RR[[x]] \to \bcinf\left(\RR\right)_0$$ splitting $T$, which sends $x$ to the germ of the identity function. If $\phi$ were a $\bcinf$-algebra map, $\phi \circ T$ would be too, and since the latter sends the generator $x$ to itself, one would have to have that $\phi \circ T=id_{\bcinf\left(\RR\right)_0}.$ This is not possible, since the existence of non-zero flat functions imply $T$ is clearly not an isomorphism.
\end{rem}

Suppose that $\bE$ is super Fermat. Consider the composite of adjunctions
\begin{equation}\label{eq:npdadj}
\xymatrix@C=2cm{\Np \ar@<-0.65ex>[r]_-{i_\bN} & \bE\Alg \ar@<-0.65ex>[l]_-{L_\bN}  \ar@<-0.65ex>[r]_-{\left(\quad\right)_\sharp} & \comsalg_\KK . \ar@<-0.65ex>[l]_-{\widehat{\left(\quad\right)}}}
\end{equation}
By Corollary \ref{cor:npdralgs}, the composite $\left(\quad\right)_\sharp \circ i_\bN$ is full and faithful. Similarly for $\bE$ Fermat. Hence we have the following corollary:

\begin{cor}\label{cor:npdffinalgs}
Suppose that $\bE$ is super Fermat. The category $\Np$ of near-point determined $\bE$ algebras is a reflective subcategory of $\comsalg_\KK.$ In particular, for each near-point determined $\bE$-algebra $\A,$ $\A$ is isomorphic to $L_\bN$ applied to $\widehat{\A_\sharp}.$ Similarly for $\bE$ Fermat.
\end{cor}

\subsubsection{Finitely generated near-point determined superalgebras}
\begin{defn}\label{dfn:jet}
Suppose that $\bE$ is super Fermat. For each $k\geq 0$, $m,n\geq 0$ the \emph{$\left(n|m\right)$-dimensional $k^{th}$ jet algebra} is defined to be the supercommutative $\KK$-algebra $\J_{n|m}^k=\KK[x^1,\ldots,x^n;\xi^1,\ldots,\xi^m]/m_0^{k+1}$, where $m_0=(x^1,\ldots,x^n;\xi^1,\ldots,\xi^m).$ Similarly for $\bE$ Fermat, one has jet algebras $\J_n^k$.
\end{defn}

\begin{rem}
Clearly, each $\J_{n|m}^k$ is a Weil algebra.
\end{rem}

\begin{prop}\label{prop:jete}
Each $$\J_{n|m}^k \cong\KK\set{x^1,\ldots,x^n;\xi^1,\ldots,\xi^m}/m_0^{k+1}$$ as an $\bE$-algebra for any super Fermat theory with ground ring $\KK.$ (And similarly for $\J_{n}^k$ when $\bE$ is not super.)
\end{prop}

\begin{proof}
Notice that by Proposition \ref{prop:completionofpresentations}, $\KK\set{x^1,\ldots,x^n;\xi^1,\ldots,\xi^m}/m_0^{k+1}$ can be identified with the $\bE$-completion of the Weil algebra $\KK[x^1,\ldots,x^n;\xi^1,\ldots,\xi^m]/m_0^{k+1},$ so we are done by Corollary \ref{cor:weilcomp}.
\end{proof}

\begin{prop}\label{weilquot}
Every Weil algebra is a quotient of some $\J_{n|m}^k$.
\end{prop}

\begin{proof}
Let $\pi:\W\to\KK$ be a Weil $\KK$-algebra. As a $\KK$-module, $$\W\cong \KK \oplus \N,$$ with $\N$ finitely generated. Let $a_1,\ldots,a_n$ be generators of $\N_\even$ and $b_1,\ldots,b_m$ be generators of $\N_\odd.$ Then there exists a surjective $\KK$-algebra map $$\phi:\KK[x^1,\ldots,x^n;\xi^1,\ldots,\xi^m] \to \W$$ sending each $x^i$ to $a_i$ and each $\xi^j$ to $b_j.$ Let $k$ be the nilpotency degree of $\N=\Ker\left(\pi\right)$. Then $I:=\Ker\left(\phi\right) \supseteq m^{k+1}_0,$ hence $\W$ is a quotient of $\J_{n|m}^k.$
\end{proof}

\begin{rem}
Both Proposition \ref{prop:jete} and Proposition \ref{weilquot} have non-finitely generated analogues; one may introduce jet algebras with infinitely many generators, and then Proposition \ref{prop:jete} holds and Proposition \ref{weilquot} holds for $\fweilk\mbox{s}.$ The proofs are the same.
\end{rem}

\begin{prop}\label{prop:weilnpd}
If $\A$ is a finitely generated $\bE$-algebra and $p:\A \to \KK$ is a $\KK$-point of $\A,$ then each $\A^{(k)}_p,$ is a Weil algebra.
\end{prop}
\begin{proof}
Since $\A$ is finitely generated, there exists a surjection $$\varphi:\KK\set{x^1,\ldots,x^n;\xi^1,\ldots,\xi^m} \to \A.$$ If $p:\A \to \KK$ is a $\KK$-point of $\A,$ then by composition there is an induced $\KK$-point $q$ of $\KK\set{x^1,\ldots,x^n;\xi^1,\ldots,\xi^m}.$ Denote by $$u:\KK\set{x^1,\ldots,x^n;\xi^1,\ldots,\xi^m} \to \KK$$ the unique $\KK$-point sending each of the generators to zero. For every $\KK$-point $t$ of $\KK\set{x^1,\ldots,x^n;\xi^1,\ldots,\xi^m},$ consider the automorphism
\begin{eqnarray*}
\theta_t:\KK\set{x^1,\ldots,x^n;\xi^1,\ldots,\xi^m} &\to& \KK\set{x^1,\ldots,x^n;\xi^1,\ldots,\xi^m}\\
x^i &\mapsto& x^i-t\left(x^i\right)\\
\xi^j &\mapsto& \xi^j-t\left(\xi^j\right),
\end{eqnarray*}
with inverse $\theta_{-t}.$ Let $M_p$ denote the kernel of $p.$ Notice that the following diagram commutes
$$\xymatrix{\KK\set{x^1,\ldots,x^n;\xi^1,\ldots,\xi^m} \ar[r]^-{\varphi \theta_q} \ar[d]_-{u} & \A \ar[ld]_-{p} \ar[d]^-{\pi}\\
\KK & \A/M^{k+1}_p \ar[l]^-{\tilde p}.}$$
Hence, the image of each generator under $\pi \varphi \theta_q$ is in $\Ker\left(\tilde p\right),$ which has nilpotency degree $k.$ It follows from Proposition \ref{prop:jete} that there is a commutative diagram
$$\xymatrix@C=0.5cm{\J_{n|m}^{k} \ar[rr]^{\tilde \varphi} \ar[rd]_-{\tilde u}& & \A/M^{k+1}_p \ar[ld]^-{\tilde p}\\
& \KK &}$$ such that $\tilde \varphi$ is surjective. Since $\J_{n|m}^{k}$ is a Weil algebra, $\Ker\left(\tilde u\right)$ is finitely generated as a $\KK$-module. Notice that we have a canonical isomorphism of $\KK$-modules $$\Ker\left(\tilde p\right) \cong \Ker\left(\tilde u\right)/\Ker\left(\tilde \varphi\right),$$ so hence $\Ker\left(\tilde p\right)$ is also finitely generated and $\A^{(k)}_p$ is a Weil algebra.
\end{proof}

\begin{notation}
Let $\bW$ denote the full subcategory of $\bE\Alg$ spanned by Weil $\KK$-algebras and their homomorphisms.
\end{notation}

\begin{cor}\label{cor:weilptnpdfin}
If $\A$ is a finitely generated $\bE$-algebra which is near-point determined, it is $\bW$-point determined.
\end{cor}

\begin{proof}
This follows immediately from Proposition \ref{prop:altnpd} and Proposition \ref{prop:weilnpd}.
\end{proof}

\begin{rem}
From Corollary \ref{cor:weilptnpdfin}, it follows that if $\A$ is a finitely generated $\cinf$-algebra, then it is near-point determined in the sense of Definition \ref{dfn:npd} if and only if it is near-point determined in the sense of \cite{msia}.
\end{rem}



\begin{rem}
Since $\bW \subset \bN,$ if $\A$ is a $\bW$-point determined $\bE$-algebra, then it is also near-point determined. By Corollary \ref{cor:weilptnpdfin}, the converse is true provided that $\A$ is finitely generated. However, the converse is \emph{not true} in general, as the following example shows. It was suggested to us by Pierre Lairez:

Let $\bE=\bcom_\KK$, $\KK$ a \emph{field}. Let $$\A:=\KK\left[u,x_1,x_2,\cdots\right]/\left(u^2,\left(x_ix_j-\delta_{ij}u\right)|_{i,j}\right).$$ Let $\bar u$ and $\bar x_i$ denote the images of $u$ and each $x_i$ in $\A,$ respectively. Consider the canonical projection to $\KK$ with kernel $$\mathfrak{m}=\left(\bar u,\bar x_1,\bar x_2,\cdots\right).$$ Notice that $\mathfrak{m}^4=0,$ so that $\A$ is a nilpotent extension of $\KK.$

Suppose that $\varphi:\A \to \W$ is a morphism to a Weil algebra which does not annihilate $\bar u.$ The image of $\varphi$ is a subalgebra of a Weil algebra. By \cite{msia}, Corollary $3.21$ $b),$ the image of $\varphi$ is also a Weil algebra. Hence, we may assume without loss of generality that $\varphi$ is surjective. Hence $$\A \to \W \to \KK$$ (where the latter map
is the one defining $\W$ as an extension of $\KK$) is also surjective, and its kernel must be $\mathfrak{m}.$ This means, that $A \to \W \to \KK$ and $\pi:\A \to \KK$ can only differ by an automorphism of $\KK$- but this automorphism must be an automorphism of $\KK$-algebras, so must be the identity. Hence, we have that the natural diagram commutes. This implies that $\Ker(\varphi)=I$ must be a subideal of $\mathfrak{m}$ which does not contain $\bar u,$ and such that $\mathfrak{m}/I$ is a finitely generated $\KK$-module.

If no finite $\KK$-linear combination of $\bar x_1, \bar x_2,\ldots$ is in $I,$ then the images of all the $\bar x_i$ in $\mathfrak{m}/I$ would be linearly independent, so $\mathfrak{m}/I$ would be infinite dimensional. Hence, there exists some $$y=\sum\limits_{i=1}^{\infty} a_i \bar x_i \in I,$$ with all but finitely many $a_i$ zero, and at least one $a_i$ non-zero. Without loss of generality, assume that $a_1$ is nonzero. Then
\begin{eqnarray*}
x_1y&=&\sum\limits_{i=1}^{\infty} a_i \bar x_1\bar x_i\\
&=&\sum\limits_{i=1}^{\infty} a_i  \delta_{i1} \bar u\\
&=&a_1 \bar u \in I,
\end{eqnarray*}
and hence $\bar u \in I,$ which is a contradiction.
\end{rem}

Finally, we note that Weil (super) algebras enjoy the following useful property (\cite{1forms}, remark directly preceding Section 2.) 

\begin{prop}\label{prop:coexp}
For any $\W\in\bW$, the functor $\W\odot-:\bE\Alg\to\bE\Alg$ has a left adjoint, $(\quad)_\W:\bE\Alg\to\bE\Alg$. The same holds in $\comalg_{\KK}$. In particular, $\W$-points of an algebra $\A$ are in bijection with $\KK$-points of $\A_\W$. Moreover (\cite{sgm}, Theorem 9.3.1) if $\A$ is free (respectively finitely generated) so is $\A_\W.$
\end{prop}

\begin{rem}
This property corresponds to exponentiability in $\bE\Alg^\op$: if $\overline{\A}$ corresponds to $\A$, $\overline{\W}$ to $\W$, then $\overline{\A_\W}=\overline{\A}^{\overline{\W}}$ in $\bE\Alg^\op$.
\end{rem}

\subsection{Flatness of near-point determined superalgebras}
In this subsection, $\bE$ will be a fixed (possibly super) Fermat theory whose ground ring $\KK$ is a \emph{field}. We will investigate the properties of the intrinsic tensor product of near-point determined $\bE$-algebras and derive some important properties of it. These properties seem to suggest that every near-point determined $\bE$-algebra is flat with respect to this intrinsic tensor product, which we state as a conjecture at the end of this subsection.

Recall that the inclusion $$\Np \hookrightarrow \bE\Alg$$ admits a left-adjoint $L_\bN$ (Proposition \ref{prop:npdreflects}). In particular, since $\bE\Alg$ is both complete and cocomplete, $\Np$ also enjoys both of these properties. However, as with all reflective subcategories, colimits are computed by first computing the colimit in the larger category, in this case $\bE\Alg,$ and then applying the reflector $L_\bN.$ In particular, from this observation, there is an intrinsic notion of tensor product (i.e. coproduct) of near-point determined $\bE$-algebras, which we shall denote by the binary operation $\oinftpy\!,$ and may not agree a priori with $\odot.$

\begin{defn}
A near-point determined $\bE$-algebra is \emph{flat} if the endofunctor $$\A \oinftpy \left(\quad\right): \Np \to \Np$$ preserves finite limits (i.e. is left exact.)
\end{defn}

\begin{lem}\label{lem:npfinp}
For any near-point determined $\bE$-algebra $A,$ the endofunctor
$$\A \oinftpy \left(\quad\right): \Np \to \Np$$ preserves finite products.
\end{lem}

\begin{proof}
The endofunctor in question factors as
$$\Np \hookrightarrow \bE\Alg \stackrel{\A\odot\left(\quad\right)}{\longlongrightarrow} \bE_\A\Alg \stackrel{\left(\quad\right) \circ u}{\longlongrightarrow} \bE\Alg \stackrel{L_\bN}{\longlongrightarrow} \Np.$$
By Corollary \ref{cor:complepresprod}, $\A\odot\left(\quad\right)$ preserves finite products, and by Remark \ref{rem:reflsprods}, $L_\bN$ does as well. The rest of the functors are right adjoints, so it follows that the  composite preserves finite products.
\end{proof}

We have the following generalization of Corollary \ref{cor:weilisgood}:

\begin{prop}\label{prop:tensweilnpd}
Let $\W$ be a $\fweilk$. Then $\W$ viewed as an $\bE$-algebra (as in Corollary \ref{cor:weilisgood}) is near-point determined. Moreover, if $\A$ is any other near-point determined $\bE$-algebra, the natural map $\A_\sharp\otimes\W\to(\A\oinftpy\W)_\sharp$ is an isomorphism. In particular, $A \otimes W$ is near-point determined as an $\bE$-algebra.
\end{prop}

\begin{proof}
By Proposition \ref{prop:sptdet}, $\W$ is clearly near-point determined. It suffices to show that $\A_\sharp\otimes\W$ with its canonical structure of an $\bE$-algebra is near-point determined. Let $$j:\A\hookrightarrow\prod_\alpha \W_\alpha$$ be an embedding into a product of $\fweil\mbox{ algebras}$, whose existence is ensured by Proposition \ref{prop:sptdet}. Then, since $\KK$ is a field, the tensor product of $\KK$-algebras is left exact (and hence also preserves monomorphisms), and we have an embedding $$j\otimes\mathrm{id_W}:\A\otimes \W \hookrightarrow (\prod_\alpha \W_\alpha)\otimes \W.$$ The canonical map $$(\prod_\alpha \W_\alpha)\otimes \W \to \prod_\alpha(\W_\alpha\otimes \W)$$ is also a monomorphism, as this is a property of vector spaces over $\KK.$
By Corollary \ref{cor:weilisgood}, each $\W_\alpha\otimes \W$ is a $\fweil\mbox{ algebra}$. So, $\A \otimes \W$ embeds into a product of $\fweil\mbox{ algebras}$, and hence is near-point determined, again by Proposition \ref{prop:sptdet}.
\end{proof}

\begin{lem}\label{lem:monoprev}
Let $\A$ be a near-point determined $\bE$-algebra. Then the endofunctor
$$\A \oinftpy \left(\quad\right): \Np \to \Np$$ preserves monomorphisms.
\end{lem}
\begin{proof}
Let $i:\B \hookrightarrow \C$ be a monomorphism of near-point determined $\bE$-algebras. Suppose that $$id_\A \oinftpy i:\A \oinftpy \B \to \A \oinftpy \C$$ is not a monomorphism. Then there exists non-zero element $k \in \A \oinftpy \B$ in its kernel. So, there exists a homomorphism $\phi:\A\oinftpy\B \to \W$ to a $\fweil\mbox{ algebra}$, such that $\phi\left(k\right) \ne 0.$ Denote by $$i_\A:\A \to \A \oinftpy \B$$ the canonical morphism, and similarly for $\B.$ Notice that the following diagram commutes:
$$\xymatrix@C=2cm@R=1.5cm{\A \oinftpy \B \ar[r]^-{\phi i_\A \oinftpy id_\B} \ar[d]_-{id_\A \oinftpy i} & \W \oinftpy \B \ar[r]^-{=} \ar[d]_-{id_{\W} \oinftpy i} & \W \otimes \B \ar[d]^-{id_{\W} \otimes i}\\
\A \oinftpy \C \ar[r]^-{\phi i_\A \oinftpy id_\C} & \W \oinftpy \C \ar[r]^-{=} & \W \otimes \C.}$$
The homomorphism $$id_{\W} \otimes i:\W \otimes \B \to \W \otimes \C$$ is a monomorphism since tensoring with $\KK$-algebras preserves monomorphisms, since $\KK$ is a field. Notice that $\phi$ factors as
$$\A \oinftpy \B \stackrel{\phi i_\A \oinftpy id_\B}{\longlongrightarrow} \W \oinftpy \B \stackrel{id_\W \oinftpy \phi i_\A}{\longlongrightarrow} \W \oinftpy \W \stackrel{\nabla}{\longlongrightarrow} \W,$$
hence $$\left(\phi i_\A \oinftpy id_\B\right)\left(k\right) \ne 0.$$ It follows that $$\left(id_{\W} \otimes i\right)\left(\phi i_\A \oinftpy id_\B\right)\left(k\right) \ne 0,$$ whereas $$\left(\phi i_\A \oinftpy id_\C\right)\left(id_\A \oinftpy i\right)\left(k\right)=0$$ since $k$ is in the kernel. This is a contradiction.
\end{proof}

Note that any left exact functor preserves monomorphisms. In particular, in order for taking the $\oinftpy$-tensor product with a near-point determined $\bE$-algebra to be left exact, it is necessary that it preserve monomorphisms. In fact, for commutative rings, the converse is true. That is, for $\R$ a commutative ring and $\A$ an $\R$-algebra, one has the following result:\\

\emph{$\A$ is flat if and only if the endofunctor $\A \otimes_\R \left(\quad\right) \bcom_\R \to \bcom_\R$ preserves monomorphisms.}
\newline

\noindent (This can be proven by using square-zero extensions to show that if tensoring with $\A$ preserves monomorphisms of $\R$-algebras, it also preserves monomorphisms of $\R$-modules.) In light of Lemma \ref{lem:monoprev}, the following conjecture seems plausible:

\begin{conj}\label{conj}
Every near-point determined $\bE$-algebra is flat.
\end{conj}

In light of Lemma \ref{lem:npfinp}, to show that Conjecture \ref{conj} is true, it suffices to show that for any near-point determined $\bE$-algebra $\A,$ the endofunctor $$\A \oinftpy \left(\quad\right): \Np \to \Np$$ preserves pullbacks, or equalizers; either would suffice\footnote{In fact, one would only need to show it preserves coreflexive equalizers.}. We offer the following partial result:

\begin{lem}
Let $$\xymatrix{\P \ar[d] \ar[r] & \C \ar[d]^-{g}\\ \B \ar[r]^-{f} & \D}$$ be a pullback diagram of near-point determined $\bE$-algebras, and consider the pullback diagram $$\xymatrix@C=2.5cm{\P' \ar[d] \ar[r] & \A \oinftpy \C \ar[d]^-{id_\A \oinftpy  g}\\ \A \oinftpy \B  \ar[r]^-{id_\A \oinftpy f} & \A \oinftpy \D.}$$ The canonical map $$\A \oinftpy P \to P'$$ is a monomorphism.
\end{lem}

\begin{proof}
Consider the canonical monomorphism $$\P \to \B \times \C.$$ By Lemma \ref{lem:monoprev} and Lemma \ref{lem:npfinp}, the induced morphism $$\A \oinftpy P \to \left(\A \oinftpy \B\right) \times \left(A \oinftpy \C\right)$$ is a monomorphism, and this map factors through $\A \oinftpy P \to P'.$
\end{proof}

\appendix
\section{Algebraic theories}\label{sec:theories}
In this appendix, we give a rapid introduction to the formalism of abstract algebraic theories. Nearly all the material may be found in \cite{algthy} and we claim no originality for it. This appendix is included merely as a convenience to the reader.
\begin{defn}\cite{borceux2,algthy}
An \emph{(abstract) algebraic theory} is a small category $\bT$ with finite products.
\end{defn}
\begin{rem}
Any algebraic theory $\bT$ has a terminal object; it is the empty product.
\end{rem}
We adjoined the parenthetical adjective \emph{abstract} since we have not provided the data of a chosen set of generators. Much of the theory of algebraic theories works well at this level of generality, but for many applications, it is important to consider the generators as part of the data. This is precisely what one needs to consider algebras as a (family of) sets with extra structure. We discuss this in Section \ref{sec:concrete}. For now, we will simply give the following definition:
\begin{defn}\label{dfn:generate}
A set $$\SS\subset\mathrm{Ob}(\bT)$$ of objects of $\bT$ is said to \emph{generate} $\bT$ as an algebraic theory if very object of $\bT$ is isomorphic to the product of finitely many objects from $\SS.$
\end{defn}

\begin{rem}
Since any algebraic theory $\bT$ is small, the set of all objects of $\bT$ is in particular (a very redundant) set of generators.
\end{rem}

\begin{defn}\cite{borceux2,algthy}
Given an algebraic theory $\bT$, a $\bT$-\emph{algebra} (in $\Set$) is defined to be a finite product preserving functor $$A:\bT\to\Set.$$  $\bT$-algebras form a category $\bT\Alg$, with natural transformations as morphisms; they span a full subcategory of the functor category $\Set^\bT$. A category $\bC$ is said to be \emph{algebraic} if it is equivalent to $\bT\Alg$ for some algebraic theory.
\end{defn}

\begin{rem}
A $\bT$-algebra $$A:\bT\to\Set$$ must send the terminal object in $\bT$ to the singleton set.
\end{rem}

\begin{rem}
Suppose that $\bT$ has a chosen generating set $\SS$. This enables us to describe an algebra $A$ by a collection of sets $$\A=\set{\A_s=A(s)|s\in\SS}$$ together with finitary operations
$$A(f):\prod_{i=1}^N\A_{s_i}^{n_i}\To\A_s,\quad s_1,\ldots,s_N,s\in\SS,\;n_i\in\NN$$
for each morphism $f$ of $\bT$, satisfying coherence conditions following from functoriality and preservation of products. Indeed, since $\SS$ generates $\bT,$ then every object $t$ is of the form $$t \cong \prod_{s \in \SS} s^{n_s}$$ for some integers $n_s \ge 0,$ with only finitely many non-zero. It follows that $$\A\left(t\right) \cong \prod_{s \in \SS} \A\left(s\right)^{n_s}.$$ When the set $\SS$ is a singleton, up to equivalence, we can assume that $\bT$ has the non-negative integers as objects, with product given by addition, $0$ as the terminal object, and with $1$ as the generator. In this case a $\bT$-algebra is the same thing as a set $\A=A(1)$ together with finitary operations $A(f):\A^n\to\A$  corresponding to the morphisms $f\in\bT(n,1)$ and satisfying structure equations coming from the composition of morphisms. In general, we may replace $\bT$, up to equivalence, by a category whose objects are $\SS$-indexed families of non-negative integers. We will discuss this in detail in Section \ref{sec:concrete}.
\end{rem}

\begin{rem}
The Yoneda embedding
$$Y_{\bT^{\op}}:\bT^\op\To\Set^\bT$$
actually factors through $\bT\Alg$ (since representable functors preserve all limits, hence in particular, finite products) and identifies $\bT^\op$ with the full subcategory of finitely-generated free $\bT$-algebras (see Remark \ref{rem:underset}).
\end{rem}

\begin{rem}
Let $Z$ be any set, and let $A:\bT \to \Set$ be an algebra for an algebraic theory $\bT$. Notice that the functor
\begin{eqnarray*}
\left(\quad\right)^{Z}:\Set &\to& \Set\\
X &\mapsto& X^{Z}=\Hom\left(Z,X\right)
\end{eqnarray*}
is right adjoint to the functor $X \mapsto X \times Z,$ so preserves all limits. In particular, the functor $$A^Z:=\left(\quad\right)^{Z} \circ A:\bT \to \Set$$ is a $\bT$-algebra. If $k$ is a generator of $\bT$, we have
\[
\A^Z_k=A^Z(k)=\A_k^Z=\Hom(Z,\A_k),
\]
with the operations applied ``pointwise''.
\end{rem}

\subsection{Sifted Colimits}
\begin{defn}
A category $\bD$ is said to be \emph{sifted} if for every finite set $X,$ (regarded as a discrete category) and every functor
$$F : \bD \times X \to \Set,$$
the canonical morphism
$$\left({\varinjlim} \prod_{x \in X} F\left(d,x\right)\right)
   \longrightarrow  \prod_{x \in X} {\varinjlim} F\left(d,x\right)$$
is an isomorphism. A \emph{sifted colimit} in $\bC$ is a colimit of a diagram $$F:\bD \to \bC,$$ with $\bD$ a sifted category.
\end{defn}

\begin{rem}
Sifted colimits commute with finite products in $\Set$ (by definition). Notice the similarity between sifted colimits, and \emph{filtered colimits}, which commute with all finite limits in $\Set.$ In particular, filtered colimits are a special case of sifted colimits.
\end{rem}

\begin{defn}
Let $\varinjlim \left(R \rightrightarrows A\right)$ be a coequalizer in a category $\bC$. It is a \emph{reflexive coequalizer} if for all objects $C,$  the induced map $$\Hom_{\bC}\left(C,R\right) \to \Hom_{\bC}\left(C,A\right) \times \Hom_{\bC}\left(C,A\right)$$ is injective, and hence determines a relation on the set $\Hom_{\bC}\left(C,A\right),$ and moreover, this relation is reflexive.
\end{defn}

In $\Set,$ one can easily check that reflexive coequalizers commute with finite products. The following proposition follows:

\begin{prop}
Reflexive coequalizers are sifted colimits.
\end{prop}

\begin{prop}\cite{algthy}\label{prop:cocompl}
A category $\bC$ is cocomplete if and only if it has all sifted colimits, and binary coproducts. Similarly, a functor preserves all small colimits if and only if it preserves all sifted colimits and binary coproducts.
\end{prop}
\begin{proof}
The standard proof that all colimits can be constructed out of arbitrary coproducts and coequalizers only uses reflexive coequalizers. The result now follows since any coproduct is a filtered colimit of finite coproducts.
\end{proof}

The following proposition is standard:

\begin{prop}\cite{sifted}
A category $\bD$ is sifted if and only if its diagonal functor is final.
\end{prop}

\begin{cor}\cite{sifted}\label{cor:sif}
Any category $\bD$ with finite coproducts is sifted.
\end{cor}

Recall that for a  small category $\bC,$ one can construct a category $\Ind\left(\bC\right)$ of $\Ind$-objects of $\bC$, that is formal filtered colimits of objects of $\bC.$ Formally, $\Ind\left(\bC\right)$ is the free cocompletion of $\bC$ with respect to filtered colimits. One can also do the analogous thing for sifted colimits:

\begin{defn}
Let $\bC$ be a small category. We define $\Sind\left(\bC\right)$ to be the free cocompletion of $\bC$ with respect to sifted colimits. It is a category under $\bC,$ $$Y_{\Sind}:\bC \to \Sind\left(\bC\right)$$ determined uniquely up to equivalence by the property that the functor $Y_{\Sind}$ satisfies the following universal property:

For all categories $\bB$ with sifted colimits, composition with $Y_{\Sind}$ induces an equivalence of categories $$\Fun_{\mathit{sift}}\left(\Sind\left(\bC\right),\bB\right) \stackrel{\sim}{\longlongrightarrow} \Fun\left(\bC,\bB\right),$$
where $\Fun_{\mathit{sift}}\left(\Sind\left(\bC\right),\bB\right)$ is the full subcategory of the functor category\\ $\Fun\left(\Sind\left(\bC\right),\bB\right)$ spanned by those functors which preserve sifted colimits.
\end{defn}

\begin{prop}\cite{algthy}\label{prop:presind}
For a small category $\bC,$ $\Sind\left(\bC\right)$ may be constructed as the full subcategory of the presheaf category $\Set^{\bC^{\op}}$- the free cocompletion of $\bC$- spanned by those presheaves which are sifted colimits of representables, and $Y_{\Sind}$ may be taken as the codomain-restricted Yoneda embedding.
\end{prop}

\begin{prop}\label{prop:cops}
The functor $$Y_{\Sind}:\bC \to \Sind\left(\bC\right)$$ preserves any finite coproducts that exist in $\bC$.
\end{prop}

\begin{proof}
Let $C$ and $D$ be objects of $\bC$ for which $C \coprod D$ exists. By Proposition \ref{prop:presind}, we can identify $\Sind\left(\bC\right)$ with a subcategory of $\Set^{\bC^{\op}}$ and $Y_{\Sind}$ with the Yoneda embedding. Let $$X=\varinjlim Y\left(E_\gamma\right)$$ be a sifted colimit. This represents an arbitrary object of $\Sind\left(\bC\right).$
We have the following chain of natural isomorphisms:
\begin{eqnarray*}
\Hom\left(Y\left(C\coprod D\right),\varinjlim Y\left(E_\gamma\right)\right) &\cong& \left(\varinjlim Y\left(E_\gamma\right)\right)\left(C\coprod D\right)\\
&\cong& \varinjlim \Hom\left(C \coprod D,E_\gamma\right)\\
&\cong& \varinjlim \left(\Hom\left(C,E_\gamma\right) \times \Hom\left(D,E_\gamma\right)\right)\\
&\cong& \left(\varinjlim \Hom\left(C,E_\gamma\right)\right) \times \left(\varinjlim \Hom\left(D,E_\gamma\right)\right)\\
&\cong& \varinjlim Y\left(E_\gamma\right)\left(C\right) \times \varinjlim Y\left(E_\gamma\right)\left(D\right)\\
&\cong& \Hom\left(C,\varinjlim Y\left(E_\gamma\right)\right) \times \Hom\left(D,\varinjlim Y\left(E_\gamma\right)\right).\\
\end{eqnarray*}
\end{proof}

\begin{cor}
If $\bC$ has binary coproducts, then $\Sind\left(\bC\right)$ is cocomplete.
\end{cor}

\begin{proof}
By Proposition \ref{prop:cocompl}, it suffices to show that $\Sind\left(\bC\right)$ has binary coproducts. However, by Proposition \ref{prop:cops}, coproducts of objects in the essential image of $Y_{\Sind}$ exist in $\Sind\left(\bC\right)$. Since every object of $\Sind\left(\bC\right)$ is a sifted colimit of representables, the result follows.
\end{proof}

\begin{cor}
If $\bC$ has binary coproducts, then for any cocomplete category $\bB$, composition with $Y_{\Sind}$ induces an equivalence of categories $$\Fun_{\mathit{cocont.}}\left(\Sind\left(\bC\right),\bB\right) \stackrel{\sim}{\longlongrightarrow} \Fun_{\coprod}\left(\bC,\bB\right),$$
where $\Fun_{\mathit{cocont.}}\left(\Sind\left(\bC\right),\bB\right)$ is the full subcategory of the functor category\\ $\Fun\left(\Sind\left(\bC\right),\bB\right)$ spanned by those functors which preserve all colimits, and\\ $\Fun_{\coprod}\left(\bC,\bB\right)$ is the full subcategory of $\Fun\left(\bC,\bB\right)$ spanned by those functor which preserve binary coproducts.
\end{cor}

\begin{cor}
If $\bC$ has binary coproducts, $\Sind\left(\bC\right)$ is reflective in $\Set^{\bC^{\op}}.$
\end{cor}

\begin{proof}
It is easily checked that the left Kan extension $\Lan_{Y}\left(Y_{\Sind}\right)$ of $Y_{\Sind}$ along the Yoneda embedding, which exists by virtue of the cocompleteness of $\Sind\left(\bC\right)$, is a left adjoint to the inclusion $$\Sind\left(\bC\right) \hookrightarrow \Set^{\bC^{\op}}.$$
\end{proof}

\begin{cor}
If $\bC$ has binary coproducts, $\Sind\left(\bC\right)$ is locally finitely presentable, so in particular is complete and cocomplete. Moreover, limits and sifted (and hence filtered) colimits are computed pointwise.
\end{cor}

\begin{proof}
The inclusion $$\Sind\left(\bC\right) \hookrightarrow \Set^{\bC^{\op}}$$ preserves sifted colimits by construction, hence in particular, filtered colimits, so is accessible. Since $\Sind\left(\bC\right)$ is fully reflective in $\Set^{\bC^{\op}},$ it follows from \cite{locpres}, Proposition 1.46, that $\Sind\left(\bC\right)$ is locally finitely presentable. The final statement is true by construction, from Proposition \ref{prop:presind}.
\end{proof}

\begin{thm}\cite{algthy}
Let $\bT$ be an algebraic theory. Then its category of algebras, $\bT\Alg$, is equivalent to $\Sind\left(\bT^{\op}\right).$
\end{thm}

\begin{proof}
Any representable presheaf is clearly an algebra, and therefore so is any sifted colimit of representables. Hence every functor $$F:\bT= \left(\bT^{\op}\right)^{\op}\to\Set$$ in $$\Sind\left(\bT^{\op}\right) \subseteq \Set^\bT$$ is a $\bT$-algebra. It suffices to show that if $\A$ is a $\bT$-algebra, then it is a sifted colimit of representable presheaves. The functor $$\A:\left(\bT^{\op}\right)^{\op} \to \Set$$ is canonically the colimit of $$\int\limits_{\bT^{\op}} \A \to \Set^{\left(\bT^{\op}\right)^{\op}}=\Set^\bT,$$ where $\left(\int\limits_{\bT^{\op}} \A\right)$ is Grothendieck construction of the presheaf $\A.$ It therefore suffices to show that $\left(\int\limits_{\bT^{\op}} \A\right)$ is sifted. By Corollary \ref{cor:sif}, it suffices to show that $\left(\int\limits_{\bT^{\op}} \A\right)^{\op}$ has binary products. The objects of $\left(\int\limits_{\bT^{\op}} \A\right)^{\op}$ can be described as pairs $\left(t,\alpha\right)$ such that $t \in \bT$ and $\alpha \in \A\left(t\right).$ Arrows $$\left(t,\alpha\right) \to \left(t',\alpha'\right)$$ are morphisms $$g:t' \to t$$ such that $$\A\left(g\right)\left(\alpha'\right)=\alpha.$$ It follows that $$\left(t,\alpha\right) \times \left(t',\alpha'\right)=\left(t \times t',\left(\alpha,\alpha'\right)\in \A\left(t\right) \times \A\left(t'\right)=\A \left(t \times t'\right)\right).$$
\end{proof}

\begin{cor}
For an algebraic theory $\bT$, its category of algebras $\bT\Alg$ is locally finitely presentable, so in particular is complete and cocomplete. Moreover, limits and sifted (and hence filtered) colimits are computed pointwise.
\end{cor}

\subsection{Morphisms of Theories}

\begin{defn}
Algebraic theories naturally form a $2$-category $\ath$. A \emph{morphism of algebraic theories} $$\bT \to \bT'$$ is a finite product preserving functor. A $2$-morphism is simply a natural transformation of functors. We will mostly be concerned only with truncation $\ath$ to a $1$-category in this paper, and denote it by $\ATh.$ 
\end{defn}

\begin{rem}
Any morphism of algebraic theories must preserve the terminal object, since it is the empty product.
\end{rem}

We will now show that any morphism $F:\bT\to\bT'$ of algebraic theories induces an adjunction $F_! \dashv F^*$ between the corresponding categories of algebras:

\begin{equation}\label{eqadj}
\mbox{$$\Adj{F^*}{\bT\Alg}{\bT'\Alg}{F_!},$$}
\end{equation}
To construct these, observe that the functor $F^{\op}:\bT^{\op} \to \bT'^{\op}$ induces three adjoint functors $F_! \dashv F^* \dashv F_*$
$$\xymatrix{\Set^{\bT} \ar@<-0.8ex>[r] \ar@<0.8ex>[r]  & \Set^{\bT'} \ar[l]}.$$ The adjunction with which we will be concerned is $$F_! \dashv F^*.$$
Indeed, $F_!$ is given as the left Kan extension $\Lan_{Y_{\bT^{\op}}} \left(Y_{\bT'^{\op}} \circ F^\op\right)$
\begin{equation}\label{eq:freealg}
\xymatrix{\Set^{\bT} \ar@{-->}[r]_-{F_!} & \Set^{\bT'} \\
\bT^{\op} \ar@{^{(}->}[u]^-{Y_{\bT^{\op}}} \ar[r]^-{F^{\op}} & \bT'^{\op} \ar@{^{(}->}[u]_-{Y_{\bT'^{\op}}}}
\end{equation}
of $Y_{\bT'^{\op}} \circ F^{\op}$ along the Yoneda embedding $Y_{\bT^{\op}}:\bT^{\op} \hookrightarrow \Set^{\bT},$ so that $F_!$ is the unique colimit preserving functor which agrees with $F^{\op}$ along representables. By the Yoneda Lemma, it follows that if $X \in \Set^{\bT'},$
\begin{eqnarray*}
F^*\left(X\right)\left(t\right) &\cong& \Hom\left(Y_{\bT^\op}\left(t\right),F^*\left(X\right)\right)\\
&\cong& \Hom\left(F_!\left(Y_{\bT^\op}\left(t\right)\right),X\right)\\
&\cong& \Hom\left(Y_{\bT'^\op}\left(F\left(t\right)\right),X\right) \\
&\cong& \left(X \circ F\right)\left(t\right),\\
\end{eqnarray*}
so that $F^*$ is given simply by precomposition with $F.$ It follows that if $X$ preserves finite products, so does $F^*\left(X\right)$. So there is an induced functor $$F^*:\bT'\Alg \to \bT\Alg.$$
The functor $$F^*:\Set^{\bT'} \to \Set^{\bT}$$ has a right adjoint $F_*,$ which by the Yoneda Lemma is given by the formula
$$F_*\left(X\right)\left(t'\right)=\Hom\left(F^*Y_{\bT'^{\op}}\left(t'\right),X\right).$$ If $X$ happened to be a $\bT$-algebra, there is no guarantee that $F_*\left(X\right)$ is a $\bT'$-algebra, so there is in general no right adjoint to $F^*$ at the level of algebras.
\begin{rem}\label{rem:rightadjoint}
Indeed, $$F^*:\Set^{\bT'} \to \Set^\bT,$$ since it is a left adjoint, preserves all colimits, and these colimits are computed pointwise. It follows that $$F^*:\bT'\Alg \to \bT\Alg$$ at least preserves \emph{sifted} colimits, as these are also computed pointwise. In fact, since both $\bT\Alg$ and $\bT'\Alg$ are locally presentable, it follows by the Adjoint Functor Theorem (\cite{locpres} Theorem 1.66), that $F^*$ has a right adjoint at the level of algebras, if and only if it preserves all small colimits. Since $F^*$ preserves sifted colimits, by Proposition \ref{prop:cocompl}, it follows that $F^*$ has a right adjoint if and only if it preserves finite coproducts.
\end{rem}
Since $F_! \vdash F^*,$ by the universal property of left Kan extensions, another characterization of $F_!$ is that $F_!\left(Z\right)$ is itself the left Kan extension of $Z$ along $F$, that is $$F_!=\Lan_{F}\left(\quad\right):Z \mapsto \Lan_{F}\left(Z\right).$$ Notice that if $Z$ is in fact a $\bT$-algebra, then $Z$ is a sifted colimit of representables, $$Z\cong \varinjlim Y_{\bT^\op}\left(t_\alpha\right).$$ It follows that,
$$F_!\left(Z\right)=\Lan_{F}\left(\varinjlim Y_{\bT^\op}\left(t_\alpha\right)\right) \cong \varinjlim Y_{\bT'^\op}\left(F\left(t_\alpha\right)\right)$$ is a sifted colimit of representables, hence a $\bT'$-algebra. Therefore, $F_!$ restricts to a functor $$F_!:\bT\Alg \to \bT'\Alg.$$

\emph{In summary:}

\medskip
From a morphism of theories $F:\bT \to \bT'$ ones gets an adjunction
$$\Adj{F^*}{\bT\Alg}{\bT'\Alg}{F_!},$$
such that $F^*$ preserves sifted colimits.

This suggests the following notion of a morphism of algebraic categories:

\begin{defn}\label{defn:algcat}
An \emph{algebraic morphism} from one algebraic category $\bC$ to another $\bD,$ is an adjunction
$$\Adj{f^*}{\bC}{\bD}{f_!},$$ such that the right adjoint $f^*$ preserves sifted colimits. With this notion of morphism, algebraic categories naturally form a $2$-category $\algcat$, whose $2$-morphisms between $\left(f^*,f_!\right)$ and $\left(g^*,g_!\right)$ are given by natural transformations $$\alpha:f^* \Rightarrow g^*.$$ We similarly denote its $1$-truncation by the $1$-category $\AlgCat.$
\end{defn}

\begin{rem}
This definition of morphism is dual to that of \cite{algthy}.
\end{rem}

\begin{rem}\label{rem:regepi}
By \cite{algthy}, Theorem 8.19, a limit preserving functor $f:\bC \to \bD$ between algebraic categories preserves sifted colimits if and only if it preserves filtered colimits and regular epimorphisms. Hence, one may equivalently say a morphism of algebraic theories $$F:\bT \to \bT'$$ induces an adjunction
$$\Adj{F^*}{\bT\Alg}{\bT'\Alg}{F_!},$$
such that $F^*$ preserves filtered colimits and regular epimorphisms.
\end{rem}

\begin{rem}
There are some size issues with $2$-category $\algcat$ in Definition \ref{defn:algcat}; morphisms may form a proper class. However, there is no cause for concern as $\algcat$ is at least essentially small, as guaranteed by the duality theorem \cite{algthy}, Theorem 8.14. Indeed, $\algcat$ is equivalent to a full subcategory of $\ath$.
\end{rem}

\begin{rem}
The morphisms in $\algcat^{op}$ may be described as limit preserving functors which preserve sifted colimits. The existence (and uniqueness) of a left adjoint follow from the Adjoint Functor Theorem (\cite{locpres} Theorem 1.66).
\end{rem}

\begin{rem}\label{rem:conservative}
If $F:\bC \to \bD$ is an essentially surjective functor, then $$F^*:\Set^{\bD^\op} \to \Set^{\bC^\op}$$ is faithful and conservative. In particular, if $$F:\bT \to \bT'$$ is an essentially surjective morphism of algebraic theories, then $$F^*:\bT'\Alg \to \bT\Alg$$ is faithful and conservative. Moreover, it preserves \emph{and reflects} all limits and sifted colimits.
\end{rem}

\begin{rem}
The construction outlined in the subsection naturally extends to a $2$-functor $$\ath \to \algcat$$ which sends a morphism  $$F:\bT \to \bT'$$ to the algebraic morphism $\left(F^*,F_!\right).$
\end{rem}

\begin{defn}
When $F$ is faithful, one calls $\bT$ as a \emph{sub-theory} of $\bT'$; in this case, $F^*(\A')$ can be thought of as the underlying $\bT$-algebra of a $\bT'$-algebra $\A'$, while $F_!(\A)$ is the free $\bT'$-algebra generated by a $\bT$-algebra $\A$, or its $\bT'$-\emph{completion}.
\end{defn}

\begin{notation}
When $\bT$ is a sub-theory of $\bT',$ we shall often neglect to mention the inclusion functor and denote the underlying $\bT$-algebra functor by $(\quad)_\sharp$ and its left adjoint, the $\bT'$-completion, by $\widehat{(\quad)}$.
\end{notation}

\begin{rem}
By general considerations, it follows that if $F$ is full and faithful, so is $F_!.$
\end{rem}

\begin{rem}\label{rem:rightadj}
It may happen that $F:\bT\to\bT'$ has a right adjoint $G:\bT'\to\bT$ (so $G^\op$ is \emph{left} adjoint to $F^\op$). Since $G$ is a right adjoint, it automatically preserves products. Hence $G$ induces an adjunction $$(G_! \dashv G^*)$$ between the categories of algebras. Notice that for $t' \in \bT'$ and $t \in \bT,$
\begin{eqnarray*}
 \Hom\left(t,G^*Y_{\bT'^{\op}}\left(t'\right)\right) &\cong& \Hom\left(G_!Y_{\bT^{\op}}\left(t\right),Y_{\bT'^{\op}}\left(t'\right)\right)\\
&\cong& \Hom\left(Y_{\bT'^{\op}}\left(G\left(t\right)\right),Y_{\bT'^{\op}}\left(t'\right)\right)\\
&\cong& \Hom_{\bT'^{\op}}\left(G\left(t\right),t'\right) \\
&\cong& \Hom_{\bT'}\left(t',G\left(t\right)\right)\\
&\cong& \Hom_{\bT}\left(F\left(t'\right),t\right)\\
&\cong& \Hom_{\bT^{\op}}\left(t,F\left(t'\right)\right).\\
\end{eqnarray*}
Hence $$G^*:\Set^{\bT'} \to \Set^{\bT}$$ is colimit preserving (since it has a right adjoint $G_*$) and for all $t'$, $$G^* \circ Y_{\bT'^{\op}}\left(t'\right) = Y_{\bT^{\op}}\left(F\left(t'\right)\right).$$ It follows that $G^*=F_!,$ hence $F_!$ acquires a further left adjoint, namely $G_!$. These then restrict to a triple of adjunctions $G_! \dashv F_! \dashv F^*:$
$$\xymatrix{\bT\Alg \ar[r] & \bT'\Alg \ar@<-0.8ex>[l] \ar@<0.8ex>[l] }.$$
\end{rem}

\section{Multisorted Lawvere Theories}\label{sec:concrete}
We now go on to describe the extra data needed to attach to an abstract algebraic theory in order to give a good sense of underlying set (or family of sets) to its algebras. Again, most of this material can be found in \cite{algthy}, however this appendix also contains some examples and notation important in this paper.

\begin{defn}
An \emph{$\SS$-sorted Lawvere theory} is an algebraic theory $\bT$ together with a injection $$\varphi_\bT:\SS \hookrightarrow \bT_0$$ whose image generates $\bT$ in the sense of definition \ref{dfn:generate}. A morphism of $\SS$-sorted Lawvere theories $$\left(\bT,\varphi_\bT\right) \to \left(\bT',\varphi_{\bT'}\right)$$ is a (natural equivalence class of a) morphism of algebraic theories $$F:\bT \to \bT'$$ such that for all $s \in \SS,$ $$F\left(\varphi_\bT\left(s\right)\right)=\varphi_\bT'\left(s\right).$$ We denote the associated category $\SS\LTh.$ When $\SS$ is a singleton set, we call an $\SS$-sorted Lawvere theory simply a Lawvere theory, and denote the corresponding category by $\LTh.$
\end{defn}

\begin{rem}
Up to equivalence, one can regard an $\SS$-sorted Lawvere theory as a category whose objects are $\SS$-indexed families of non-negative integers. This allows us to refer to a theory $\bT$ without reference to its structural map $\varphi_\bT.$
\end{rem}

\begin{rem}
One can expand this definition by removing the injectivity of the map $\varphi_\bT$ and nothing is lost. We will refer to such a pair $\left(\bT,\varphi\right)$ with $\varphi$ not necessarily injective as an \emph{$\SS$-indexed} Lawvere theory.
\end{rem}

\begin{rem}\label{rem:esssurj}
Let $f:\bT \to \bT'$ be a morphism of $\SS$-sorted Lawvere theories. Then, up to equivalence, $f$ is a bijection on objects.
\end{rem}

\begin{rem}\label{rem:ff}
If $F:\bT \to \bT'$ is a full and faithful morphism of $\SS$-sorted Lawvere theories, then it is an equivalence.
\end{rem}

\begin{defn}
An \emph{algebra} for an $\SS$-sorted Lawvere theory $\bT$ is simply an algebra for its underlying algebraic theory.
\end{defn}

\begin{prop}\label{prop:lawvmorph}
A morphism $F:\bT \to \bT'$ of Lawvere theories induces an adjunction
$$\Adj{F^*}{\bT\Alg}{\bT'\Alg}{F_!},$$
such that $F^*$ preserves and reflects all limits and sifted colimits (equivalently all limits, filtered colimits, and regular epimorphisms).
\end{prop}
\begin{proof}
This follows from Remark \ref{rem:conservative} and Remark \ref{rem:esssurj}. The parenthetical remark follows from Remark \ref{rem:regepi}.
\end{proof}

\begin{eg}\cite{algthy}\label{eg:freeprod}
Let $\SS$ be any set viewed as a discrete category. Let $\bT_\SS$ be the free completion of $\SS$ with respect to finite products. Concretely, the objects of $\bT_\SS$ are finite families $$\left(s_i \in \SS\right)_{i \in I}$$ and morphisms $$\left(s_i \in \SS\right)_{i \in I} \to \left(t_j \in \SS\right)_{j \in J}$$ are functions of finite sets $$f:J \to I$$ such that $$s_{f\left(j\right)}=t_j$$ for all $j \in J.$ There is a canonical functor $$Y_{\prod}:\SS \to \bT_\SS,$$ sending each element $s \in \SS$ to $\left(s\right)$ viewed as a finite family with one element. The universal property of this functor is that for any category $\bD$ with finite products, composition with $Y_{\prod}$ induces an equivalence of categories
\begin{equation}\label{eq:freeprod}
\chi:\Fun_{\prod}\left(\bT_\SS,\bD\right) \stackrel{\sim}{\longrightarrow} \Fun\left(\SS,\bD\right)=\bD^{\SS},
\end{equation}
where $\Fun_{\prod}\left(\bT_\SS,\bD\right)$ is the full subcategory of the functor category on those functors which preserve finite products. It follows that $\bT$ is an algebraic theory whose category of algebras is equivalent to $\Set^\SS$. Moreover, the functor $Y_{\prod}$ gives $\bT_\SS$ the canonical structure of an $\SS$-sorted Lawvere theory.
\end{eg}

\begin{rem}
When $\SS$ is a singleton, $\bT_\SS$ is equivalent to $\mathbf{FinSet^\op},$ the opposite category of finite sets.
\end{rem}

By construction, we have the following proposition:

\begin{prop}\label{prop:initial}
The $\SS$-sorted Lawvere theory $\bT_\SS$ is an initial object.
\end{prop}

Let $\left(\bT,\varphi_\bT\right)$ be an $\SS$-sorted Lawvere theory. From Proposition \ref{prop:initial}, we know that there is a unique morphism of $\SS$-sorted Lawvere theories $$\sigma_\bT:\bT_\SS \to \bT.$$ From Proposition \ref{prop:lawvmorph}, we have the following Corollary:

\begin{cor}\cite{algthy}\label{cor:concretess}
With $$\sigma_\bT:\bT_\SS \to \bT$$ as above, the functor $$U_{\bT}:=\left(\sigma_\bT\right)^*:\bT\Alg \to \Set^{\SS}$$ is faithful and conservative. In particular, it preserves and reflects limits, sifted colimits, monomorphisms, and regular epimorphisms.
\end{cor}

\begin{rem}
The above Corollary only needs that $\sigma_\bT$ is essentially surjective, so it holds true for $\SS$-indexed Lawvere theories too.
\end{rem}

\begin{rem}
In particular, for an $\SS$-sorted Lawvere theory $$\sigma_\bT:\bT_\SS \to \bT,$$ the adjunction $$\Adjlong{\left(\sigma_\bT\right)^*}{\Set^{\SS}}{\bT\Alg}{\left(\sigma_\bT\right)_!}$$ is monadic.
\end{rem}

\begin{rem}\label{rem:underset}
The left adjoint $$\left(\sigma_\bT\right)_!:\Set^{\SS} \to \bT\Alg$$ is the functor assigning an $\SS$-indexed family of sets $\bG=\left(G_s\right)_{s \in \SS}$ the free $\bT$-algebra on $\bG,$ whereas the right adjoint $\left(\sigma_\bT\right)^*$ assigns a $\bT$-algebra its underlying $\SS$-indexed set. Explicitly, if $$\A:\bT \to \Set$$ is a $\bT$-algebra, then the underlying $\SS$-indexed set is $$\left(\A_s=\A\left(\sigma_\bT\left(s\right)\right)\right)_{\{s \in  \SS\}}.$$
The diagram \eqref{eq:freealg} becomes in this case
\[
\xymatrix{\Set \ar@{-->}[r]_-{u_!} & \Set^{\bT} \\
\bT_{\SS} \ar@{^{(}->}[u] \ar[r]^-{\sigma_{\bT}^{\op}} & \bT^{\op}. \ar@{^{(}->}[u]_-{Y_{\bT^{\op}}}}
\]
It follows that $Y_{\bT^{\op}}$ establishes an equivalence of categories between $\bT^\op$ and the full subcategory of $\bT\Alg$ consisting of finitely generated free $\bT$-algebras.
\end{rem}

\begin{rem}
Since every algebra is a sifted colimit of representables, which by the previous remark are precisely the finitely generated free algebras, by Remark \ref{rem:rightadjoint}, it follows that for a map of Lawvere theories $$F:\bT \to \bT',$$ the functor $$F^*:\bT'\Alg \to \bT\Alg$$ has a right a adjoint if and only if for each pair $A,B$ of finitely generated free $\bT'$-algebras, $$F^*\left(A\coprod B\right)=F^*\left(A\right) \coprod F^*\left(B\right).$$
\end{rem}

\begin{notation}
For a Lawvere theory $\bT$, denote the free $\bT$-algebra on $n$ generators by $\bT(n)$. As discussed above, $\bT(n)$ can be identified with the representable functor $Y_{\bT^{\op}}\left(n\right).$ Since the underlying set is given by evaluation at the generator, $1,$ it follows that the underlying set is
\begin{equation}\label{eq:notta}
Y_{\bT^{\op}}\left(n\right)\left(1\right)=\Hom_{\bT^{\op}}\left(1,n\right)=\Hom_{\bT}\left(n,1\right).
\end{equation}
We will use the notation $\bT(n,1)$ for this set of morphisms to emphasize that it encodes the $n$-ary operations of the theory $\bT$. In the $\SS$-sorted case, we adopt the notation $\bT(\{n_s\}_{s\in \SS})$ for the free algebra with $n_s$ generators of sort $s,$ for each $s \in \SS.$ The underlying $\SS$-indexed set of such a free algebra is then $$\left(\bT\left(\varphi_{\bT}\left(s\right)^{n_s},\varphi_{\bT}\left(s\right)\right)\right)_{s\in\SS}.$$
\end{notation}

\begin{rem}\label{rem:relative} 
Let $\bT$ be an $\SS$-sorted Lawvere theory. Given an $\A\in\bT\Alg$, an \emph{$\A$-algebra} is, by definition, an object of the undercategory $\A/\bT\Alg$. Given a morphism of theories $F:\bT\to\bT'$, for each $\A\in\bT'\Alg$ we have an induced adjunction of undercategories
$$\Adj{F^\A_!}{\A/\bT'\Alg}{F^*\A/\bT\Alg.}{F_\A^*}$$
To see this, observe that $\A$-algebras are algebras over the theory $\bT_\A$ whose operations are labeled by elements of free finitely generated $\A$-algebras, i.e.
\[
\bT_\A(\{n_s\}_{s\in \SS})=\A\amalg\bT(\{n_s\}_{s\in \SS}).
\]
The above adjunction is induced by the morphism of theories
\[
F_\A:\bT_{F^*\A}\To\bT'_\A.
\]
Notice that there is also canonical morphisms of theories $$u_\A:\bT' \to \bT'_\A$$ and $$u_{F^*\A}:\bT \to \bT_{F^*\A},$$ such that the following diagram commutes: $$\xymatrix{\bT \ar[r]^-{F} \ar[d]_-{u_{F^{*}\A}} & \bT' \ar[d]^-{u_\A}\\ \bT_{F^{*}\A} \ar[r]_-{F_\A} & \bT'_\A.}$$
It follows that $F^*_\A$ maps $\A \to \A'$ to $F^* \A \to \F^*\A'$ and that $F^\A_!$ takes any $F^*\A$-algebra of the form $$F^*\A \to F^*\A \coprod \B,$$ for a $\bT$-algebra $\B,$ to $$\A \to \A \coprod F_!\B.$$
\end{rem}

\begin{rem}\label{rem:image}
Let $f:\bT \to \bT'$ be a morphism of $\SS$-sorted Lawvere theories. The category $\mathbf{Cat}$ of small categories carries a factorization system. Faithful functors are right orthogonal to functors which are both essentially surjective and full. It turns out that if $f$ is a morphism of $\SS$-sorted Lawvere theories, if $$\bT \to \bC \to \bT'$$ is the unique factorization by an essentially surjective and full functor, followed by a faithful one, then $\bC$ is an $\SS$-sorted Lawvere theory, and all the functors are maps of $\SS$-sorted Lawvere theories. In this case, we denote $\bC$ by $\Im\left(f\right),$ and call it the image of $f.$\footnote{The more traditional notion of image of a functor, uses the factorization system determined by essentially surjective functors, and full and faithful functors. This factorization is ill suited for Lawvere theories due to Remark \ref{rem:ff}.}

We now proceed to construct $\Im\left(f\right).$ We may assume without loss of generality that $f$ is a bijection on objects. Consider for each pair of objects $\left(x,y\right) \in \bT_0$ the induced function
$$f_{x,y}:\Hom_\bT\left(x,y\right) \to \Hom_\bT'\left(x,y\right).$$ Then one can define a new category $\Im\left(f\right)$ with the same objects as $\bT'$ but whose morphisms are given by $$\Hom\left(x,y\right):=\Im\left(f_{x,y}\right),$$ the image of the function $f_{x,y}.$ Since $f$ preserves limit diagrams for finite products, it follows that $\bT'$ has finite products and the induced functor $$\Im\left(f\right) \to \bT'$$ preserves them. One hence gets a factorization of $f$
\begin{equation}\label{eq:factt}
\bT \to \Im\left(f\right) \to \bT'
\end{equation}
by algebraic functors, each of which is a bijection on objects. In particular, the structure map for $\bT,$ $\sigma_\bT,$ gives $\Im\left(f\right)$ the canonical structure of a $\SS$-sorted Lawvere theory in such a way that the factorization (\ref{eq:factt}) consists of morphisms of $\SS$-sorted Lawvere theories. We refer to the theory $\Im\left(f\right)$ as the \emph{image} of $f$. It is clear that the induced map $$\Im\left(f\right) \to \bT'$$ is faithful. We conclude, that the factorization system determined by essentially surjective and full functors and faithful ones descends to a factorization system on $\SS\LTh.$
\end{rem}

\begin{eg}\label{eg:end}
Let $A$ be any set. Define the Lawvere theory $\End_A$ to be the full subcategory of $\Set$ generated by the finite Cartesian powers of $A$, i.e. $$\End_A(n,1)=\Set(A^n,A).$$ Picking out $A$ as the generator makes $\End_A$ into a Lawvere theory. Notice that the full and faithful inclusion $$\End_A \hookrightarrow \Set$$ preserves finite products, so that given a morphism of Lawvere theories $$F:\bT\to\End_A,$$ by composition, one gets an $\bT$-algebra in $\Set.$ Since $F$ preserves the generator, this induced $\bT$-algebra will have underlying set $A.$ It follows that a $\bT$-algebra structure on $A$ for a Lawvere theory $\bT$ is the same thing as a morphism of Lawvere theories $$F:\bT\to\End_A.$$ In the same way, given an $\SS$-indexed family of sets $\A=\left(A_s\right)_{s\in \SS} \in \Set^{\SS}$, one obtains a $\SS$-sorted Lawvere theory $\End_\A$, such that for all $\SS$-sorted Lawvere theories $\bT,$ morphisms of $\SS$-sorted Lawvere theories $$\bT\to\End_\A$$ are in bijection with $\bT$-algebras whose underlying $\SS$-indexed set is $\A.$ Explicitly $\End_\A$ is the full subcategory of $\Set$ generated by the collection of sets $\left(A_s\right)$ for each $s$; these are also the sorts.
\end{eg}

\subsection{Congruences}
\begin{defn}
Let $\bC$ be a category with finite products. An \emph{equivalence relation} on an object $A \in \bC$ is a subobject $$R \rightarrowtail A \times A$$ such that for all objects $C$, $$\Hom_{\bC}\left(C,R\right) \rightarrowtail \Hom_{\bC}\left(C,A\right) \times \Hom_{\bC}\left(C,A\right)$$ is an equivalence relation on the set $\Hom_{\bC}\left(C,A\right).$
\end{defn}

\begin{defn}
Given an equivalence relation $$R \rightarrowtail A \times A,$$ one may consider the induced pair of maps
\begin{equation}\label{eq:relation}
R \rightrightarrows A.
\end{equation}
If the coequalizer of this diagram exists, it is called the \emph{quotient object} $A/R$ of $A$ by the equivalence relation $R.$ \end{defn}

\begin{rem}
In the case that $\bC=\Set,$ one recovers the usual notion of the quotient of a set by an equivalence relation.
\end{rem}

\begin{rem}
A coequalizer of the form (\ref{eq:relation}) is a reflexive coequalizer, hence in particular, a sifted colimit.
\end{rem}

\begin{defn}
Let $\bT$ be an $\SS$-sorted Lawvere theory. An equivalence relation in $\bT\Alg$ is called a \emph{congruence}.
\end{defn}

\begin{defn}
A quotient $A \mapsto A/R$ by an equivalence relation is called an \emph{effective quotient} if the canonical map $$A \to A/R$$ is an effective epimorphism, i.e. $$A/R \cong \varinjlim \left(A \times_{A/R} A \rightrightarrows A\right).$$
\end{defn}

\begin{prop}
Suppose that $\bT$ is an $\SS$-sorted Lawvere theory, then every quotient in $\bT\Alg$ is effective.
\end{prop}

\begin{proof}
For any equivalence relation with a quotient, by definition, the map $$A \to A/R$$ is a regular epimorphism. However, since $\bT\Alg$ has pullbacks, every regular epimorphism is an effective epimorphism, so we are done.
\end{proof}

\begin{cor}
Suppose that $\bT$ is an $\SS$-sorted Lawvere theory, then every regular epimorphism is of the form $$A \to A/R$$ for some congruence $R$ on $A.$
\end{cor}
\begin{proof}
Maps of the form $A \to A/R$ are regular by definition. Conversely, suppose that $$A \to B$$ is a regular epimorphism. Then, since $\bT\Alg$ has pullbacks, it is an effective epimorphism, and hence the map is induced by a colimiting cocone witnessing $B$ as the colimit of $$A \times_B A \rightrightarrows B.$$ This coequalizer is the quotient for the equivalence relation $$R:=A \times_B A \rightarrowtail A \times A,$$ hence is of the form $$A \to A/R.$$
\end{proof}

\begin{prop}\label{prop:congrues}
Suppose that $\bT$ is an $\SS$-sorted Lawvere theory, and $R$ is a congruence on a $\bT$-algebra $A.$ Then the quotient $A/R$ exists. In particular, the underlying $\SS$-indexed set of $A/R$ is the quotient of the underlying $\SS$-indexed set of $A$ by the equivalence relation induced by $R.$
\end{prop}

\begin{proof}
From Corollary \ref{cor:concretess}, the functor $$U_{\bT}:\bT\Alg \to \Set^{\SS}$$ assigning an algebra its underlying $\SS$-indexed set, preserves and reflects reflexive coequalizers. In particular, it preserves and reflects quotients by equivalence relations. Since $\Set^{\SS}$ is a topos, it has quotients by equivalence relations, so we are done.
\end{proof}


\bibliographystyle{hplain}
\bibliography{derived}

\begin{thebibliography}{10}

\bibitem{locpres}
Ji\v{r}\'{i} Ad{\'a}mek and Ji\v{r}\'{i} Rosick{\'y}.
\newblock {\em Locally presentable and accessible categories}, volume 189 of
  {\em London Mathematical Society Lecture Note Series}.
\newblock Cambridge University Press, Cambridge, 1994.

\bibitem{sifted}
Ji\v{r}\'{i} Ad{\'a}mek, Ji\v{r}\'{i} Rosick{\'y}, and Enrico~M. Vitale.
\newblock What are sifted colimits?
\newblock {\em Theory Appl. Categ.}, 23:No. 13, 251--260, 2010.

\bibitem{algthy}
Ji\v{r}\'{i} Ad{\'a}mek, Ji\v{r}\'{i} Rosick{\'y}, and Enrico~M. Vitale.
\newblock {\em Algebraic theories}, volume 184 of {\em Cambridge Tracts in
  Mathematics}.
\newblock Cambridge University Press, Cambridge, 2011.
\newblock A categorical introduction to general algebra, With a foreword by F.
  W. Lawvere.

\bibitem{atmcd}
Michael~F. Atiyah and Ian~G. Macdonald.
\newblock {\em Introduction to commutative algebra}.
\newblock Addison-Wesley Publishing Co., Reading, Mass.-London-Don Mills, Ont.,
  1969.

\bibitem{borceux2}
Francis Borceux.
\newblock {\em Handbook of categorical algebra. 2}, volume~51 of {\em
  Encyclopedia of Mathematics and its Applications}.
\newblock Cambridge University Press, Cambridge, 1994.
\newblock Categories and structures.

\bibitem{Borisov}
Dennis Borisov.
\newblock Topological characterization of various types of rings of smooth
  functions.
\newblock \href{http://arxiv.org/abs/1108.5885}{arXiv:1108.5885}, 2011.

\bibitem{borisovnoel}
Dennis Borisov and Justin Noel.
\newblock Simplicial approach to derived differential manifolds.
\newblock \href{http://arxiv.org/abs/1112.0033}{arXiv:1112.003}, 2011.

\bibitem{arch}
Marta Bunge and Eduardo~J. Dubuc.
\newblock Archimedian local {$C\sp \infty$}-rings and models of synthetic
  differential geometry.
\newblock {\em Cahiers Topologie G\'eom. Diff\'erentielle Cat\'eg.},
  27(3):3--22, 1986.

\bibitem{dg2}
David Carchedi and Dmitry Roytenberg.
\newblock Homological algebra of superalgebras of differentiable functions.
\newblock 2012.
\newblock (In Preparation).

\bibitem{dubuc1}
Eduardo~J. Dubuc.
\newblock Sur les mod\`eles de la g\'eom\'etrie diff\'erentielle synth\'etique.
\newblock {\em Cahiers Topologie G\'eom. Diff\'erentielle}, 20(3):231--279,
  1979.

\bibitem{cinfsch}
Eduardo~J. Dubuc.
\newblock {$C\sp{\infty }$}-schemes.
\newblock {\em Amer. J. Math.}, 103(4):683--690, 1981.

\bibitem{germint}
Eduardo~J. Dubuc.
\newblock Germ representability and local integration of vector fields in a
  well adapted model of {SDG}.
\newblock {\em J. Pure Appl. Algebra}, 64(2):131--144, 1990.

\bibitem{1forms}
Eduardo~J. Dubuc and Anders Kock.
\newblock On {$1$}-form classifiers.
\newblock {\em Comm. Algebra}, 12(11-12):1471--1531, 1984.

\bibitem{joyce}
Dominic Joyce.
\newblock An introduction to d-manifolds and derived differential geometry.
\newblock \href{http://arxiv.org/abs/1206.4207}{arXiv:1206.4207}, 2012.

\bibitem{sdg}
Anders Kock.
\newblock {\em Synthetic differential geometry}, volume 333 of {\em London
  Mathematical Society Lecture Note Series}.
\newblock Cambridge University Press, Cambridge, second edition, 2006.

\bibitem{sgm}
Anders Kock.
\newblock {\em Synthetic geometry of manifolds}, volume 180 of {\em Cambridge
  Tracts in Mathematics}.
\newblock Cambridge University Press, Cambridge, 2010.

\bibitem{manin}
Yuri~I. Manin.
\newblock {\em Gauge field theory and complex geometry}, volume 289 of {\em
  Grundlehren der Mathematischen Wissenschaften [Fundamental Principles of
  Mathematical Sciences]}.
\newblock Springer-Verlag, Berlin, second edition, 1997.
\newblock Translated from the 1984 Russian original by N. Koblitz and J. R.
  King, With an appendix by Sergei Merkulov.

\bibitem{loc1}
Ieke Moerdijk and Gonzalo~E. Reyes.
\newblock Rings of smooth functions and their localizations. {I}.
\newblock {\em J. Algebra}, 99(2):324--336, 1986.

\bibitem{smzar}
Ieke Moerdijk and Gonzalo~E. Reyes.
\newblock A smooth version of the {Z}ariski topos.
\newblock {\em Adv. in Math.}, 65(3):229--253, 1987.

\bibitem{msia}
Ieke Moerdijk and Gonzalo~E. Reyes.
\newblock {\em Models for smooth infinitesimal analysis}.
\newblock Springer-Verlag, New York, 1991.

\bibitem{nish}
Hirokazu Nishimura.
\newblock Supersmooth topoi.
\newblock {\em Internat. J. Theoret. Phys.}, 39(5):1221--1231, 2000.

\bibitem{smoothfun}
Ng{\^o}~Van Qu{\^e} and Gonzalo~E. Reyes.
\newblock Smooth functors and synthetic calculus.
\newblock In {\em The {L}. {E}. {J}. {B}rouwer {C}entenary {S}ymposium
  ({N}oordwijkerhout, 1981)}, volume 110 of {\em Stud. Logic Found. Math.},
  pages 377--395. North-Holland, Amsterdam, 1982.

\bibitem{Reichard}
K.~Reichard.
\newblock Nichtdifferenzierbare {M}orphismen differenzierbarer {R}\"aume.
\newblock {\em Manuscripta Math.}, 15:243--250, 1975.

\bibitem{ReyWra}
Gonzalo~E. Reyes and Gavin~C. Wraith.
\newblock A note on tangent bundles in a category with a ring object.
\newblock {\em Math. Scand.}, 42(1):53--63, 1978.

\bibitem{spivak}
David~I. Spivak.
\newblock Derived smooth manifolds.
\newblock {\em Duke Math. J.}, 153(1):55--128, 2010.

\bibitem{yetter}
David~N. Yetter.
\newblock Models for synthetic supergeometry.
\newblock {\em Cahiers Topologie G\'eom. Diff\'erentielle Cat\'eg.},
  29(2):87--108, 1988.

\end{thebibliography}

\end{document}